\documentclass[10pt]{amsart}
\usepackage[dvipdfmx]{graphicx}
\usepackage{amsmath,amscd,amsthm,amsxtra,yhmath,stackrel}
\usepackage{epsfig,graphics,color,colortbl}
\usepackage{amssymb,latexsym} 
\usepackage{mathrsfs}

\usepackage{bm}
\usepackage[poly,all]{xy}
\usepackage{hyperref}
\usepackage{tikz-cd}
\usepackage[draft]{todonotes}
\usepackage{upgreek}
\usepackage{ytableau}
\usepackage{stmaryrd}
\usepackage{verbatim}
\usepackage{mathtools}
\usepackage[title]{appendix}
\usepackage{enumerate}
\usepackage[normalem]{ulem}
\usepackage{datetime}
\usepackage{longtable}
\usepackage[pagewise]{lineno} % for line numbers
\usepackage{kotex} % 배정이 추가했습니다.

\newtheorem{thm}{\bf Theorem}[section]

\newtheorem{prop}[thm]{\bf Proposition}
\newtheorem{cor}[thm]{\bf Corollary}
\newtheorem{lem}[thm]{\bf Lemma}
\newtheorem{rem}[thm]{\bf Remark}
\newtheorem{ex}[thm]{\bf Example}

\numberwithin{equation}{section}

\newcommand{\mc}{\mathcal}
\newcommand{\mf}{\mathfrak}
\newcommand{\ms}{\mathscr}

\newcommand{\pf}{\noindent{\bfseries Proof. }}
\newcommand{\ov}{\overline}

\newcommand{\U}{{\mc U}}

\newcommand{\V}{\mc{V}}
\newcommand{\W}{\mc{W}}
\newcommand{\cP}{\mathscr{P}}

\newcommand{\I}{\mathbb{I}}

\newcommand{\Z}{\mathbb{Z}}

\newcommand{\C}{\mathbb{C}}

\newcommand{\e}{\epsilon}
\newcommand{\de}{\delta}

\newcommand{\td}{\widetilde}

\newcommand{\lang}{\langle}
\newcommand{\rang}{\rangle}
\newcommand{\gl}{\mf{gl}}

\newcommand{\g}{\mf{g}}

\newcommand{\La}{\Lambda}

\newcommand{\la}{\lambda}

\newcommand{\hf}{\frac{1}{2}}

\newcommand{\ot}{\otimes}

\newcommand{\si}{(-1)^{\e_i}}

\newcommand{\bq}{{\bf q}}

\newcommand{\msf}{\mathsf}

\begin{document}
\title
[]{$q$-deformed Howe duality for orthosymplectic Lie superalgebras}

\author{Jeong Bae}
\address{Department of Mathematical Sciences, Seoul National University, Seoul 08826, Korea}
\email{jeongbm@snu.ac.kr}

\author{Jae-Hoon Kwon}
\address{Department of Mathematical Sciences and and Research Institute of Mathematics, Seoul National University, Seoul 08826, Korea}
\email{jaehoonkw@snu.ac.kr}

%\today
%\thanks{This work is supported by the National Research Foundation of Korea(NRF) grant funded by the Korea government(MSIT) (No.2020R1A5A1016126 and RS-2024-00342349).}

\begin{abstract}
We give a $q$-analogue of Howe duality associated to a pair $(\mf{g},G)$, where $\mf{g}$ is an orthosymplectic Lie superalgebra and $G=O_\ell, Sp_{2\ell}$. We define explicitly {commuting actions} of a quantized enveloping algebra of $\mf{g}$ and the $\imath$quantum group of {type AI and AII} on a $q$-deformed supersymmetric space, and describe its semisimple decomposition whose classical limit recovers the $(\mf{g},G)$-duality. As special cases, we obtain $q$-analogues of $(\mf{g},G)$-dualities on symmetric and exterior algebras for $\mf{g}=\mf{so}_{2n}$, $\mf{sp}_{2n}$.

%When $\e=(0^n)$ and $\e=(1^n)$, we obtain $q$-analogues of the Howe dual pairs $(\mf{g},G_\ell)$ acting on symmetric and exterior algebras with $(\mf{g},G_\ell)=(\mf{sp}_{2n},O_\ell), (\mf{so}_n,Sp_{2\ell})$ and $(\mf{so}_{n},O_\ell), (\mf{sp}_{2n},Sp_{2\ell})$, respectively.

\end{abstract}

\maketitle
\setcounter{tocdepth}{1}
%\tableofcontents

\noindent

\section{Introduction}
\subsection{}
A reductive dual pair $(\g,G)$ introduced by Howe \cite{H,H2} is a pair of semisimple Lie (super)algebra $\g$ and a Lie group $G$, which acts on a space $F$ and forms the centralizers of each other. 
It gives a representation theoretic approach to the first fundamental theorem of the classical invariant theory, and has applications in various areas of mathematics.

Let us consider the following well-known four reductive dual pairs $(\g,G)$ acting on $F$:
\begin{align}
&\text{\ $(\mf{sp}_{2m},O_\ell)$\ \  on $S(\C^\ell \ot \C^{m})$},\quad \quad \
\text{$(\mf{so}_{2n},O_\ell)$\ \  on\  $\Lambda(\C^\ell \ot \C^{n})$},\label{eq:O duality}  \\
&\text{$(\mf{so}_{2m},Sp_{2\ell})$ on $S(\C^{2\ell} \ot \C^{m})$},\quad 
\text{\ $(\mf{sp}_{2n},Sp_{2\ell})$ on $\Lambda(\C^{2\ell} \ot \C^{n})$}, \label{eq:Sp duality}
\end{align}
where $S(V)$ (resp. $\Lambda(V)$) is the symmetric (resp. exterior) algebra generated by $V$ (see also \cite{W99}). 
They are obtained from the following dual pairs of type $A$:
\begin{align}
&\text{$(\gl_{m}, GL_l)$ on $S(\C^l \ot \C^{m})$,\quad\  $(\gl_{n}, GL_l)$ on $\Lambda(\C^l \ot \C^{n})$ \quad $(l=\ell, 2\ell)$}, \label{eq:type A duality}
\end{align}
by extending the actions of general linear Lie algebra to those of simple Lie algebras of type $C$ and $D$, while the actions of $O_\ell$ and $Sp_{2\ell}$ are induced from that of $GL_l$. Note that the irreducible representations of $\mf{sp}_{2m}$ and $\mf{so}_{2n}$ appearing in $S(\C^l \ot \C^{m})$ are infinite-dimensional and often called oscillator representations (cf.~\cite[5.6]{GW}).

The dualities \eqref{eq:O duality} and \eqref{eq:Sp duality} can be given in a more uniform way in terms of Lie super algebras \cite{CZ}: Let $S(\C^l \ot \C^{m|n})$ denote the supersymmetric algebra generated by $\C^l \ot \C^{m|n}$ where $\C^{m|n}$ is the complex superspace with $(\C^{m|n})_{\ov 0}=\C^m$ and $(\C^{m|n})_{\ov 1}=\C^n$.
Then each of \eqref{eq:O duality} and \eqref{eq:Sp duality} can be merged into  a single dual pair  
\begin{align}\label{eq:spo duality}
&\text{$(\mf{spo}_{2m|2n},O_\ell)$ on $S(\C^\ell \ot \C^{m|n})$\quad and\quad $(\mf{osp}_{2m|2n},Sp_{2\ell})$ on $S(\C^{2\ell} \ot \C^{m|n})$},  %\\
%&\text{$(\mf{osp}_{2m|2n},O_l)$\ \ on $S(\C^l \ot \C^{m|n})$}, \label{eq:osp duality} 
\end{align}
respectively, where $\mf{spo}_{2m|2n}$ and $\mf{osp}_{2m|2n}$ are orthosymplectic Lie superalgebras. We also remark that there is an $(\mf{spo}_{2m|2n+1},Pin_{2\ell})$-duality on $S(\C^\ell\ot \C^{2m|2n+1})$ \cite[Theorem A.1]{CKW} (see also \cite{CW03}).

\subsection{}
The main result in this paper is to give a $q$-analogue of the dual pair $(\g,G)$ in \eqref{eq:spo duality}, which also yields $q$-analogues of \eqref{eq:O duality} and \eqref{eq:Sp duality} as special cases. 

In our $q$-analogue of \eqref{eq:O duality}, we use a quantized enveloping algebra of $\g$ introduced in \cite{KOS} which arises naturally in the study of three-dimensional Yang-Baxter equation \cite{KO}. We denote it by $\U_X(\e)$ where $\e\in\Z_2^n$ determines the type of Borel subalgebra (or fundamental system) and $X=C,D$ stands for type of Dynkin diagram of $\g$. In particular, we have $\U_X(0^n)=U_q(X_n)$ and $\U_X(1^n)=U_{-q^{-1}}(X_n)$, where $U_q(X_n)$ is the quantum group of type $X_n$. On the other hand, the action of $G$ is replaced by that of the $\imath$quantum group associated to the Lie algebra of $G$ \cite{GK,Le,No}.

Let $\W$ be a $q$-analogue of the supersymmetric algebra $S(\C^{\e})$ generated by $\C^\e$, where $\C^{\e}$ is the superspace of dimension $n$ with $\e$ the parity sequence of standard basis.
We first consider a $\U_X(\e)$-module structure on $\W$ for each $X$ such that it gives a spin representation of type $D_n$ when $(X,\e)=(D,(1^n))$, and $q$-oscillator representations of type $C_n$ when $(X,\e)=(C,(0^n))$ \cite{Ha}. For $\ell\ge 2$, $\W^{\ot\ell}$ has a $(\U_A(\e),{\bf U}(\mf{sl}_\ell))$-module structure, where $\U_A(\e)\subset \U_X(\e)$ is the subalgebra of type $A$ and ${\bf U}(\mf{sl}_\ell)$ is the Drinfeld-Jimbo quantum group of $\mf{sl}_\ell$. Then we show that the action of $\U_X(\e)\supset \U_A(\e)$ on $\W^{\ot\ell}$ commutes with that of the coideal subalgebra or $\imath$quantum group ${\bf U}^\imath(\mf{so}_\ell)\subset {\bf U}(\mf{sl}_\ell)$ corresponding to the symmetric pair $({\mf{sl}}_\ell,\mf{so}_\ell)$ whose Satake diagram is of type AI.

Next, we consider a $\U_X(\e)$-module structure on the space $\W^{\ot 2}$ for each $X$, which is not the $\U_X(\e)$-module induced from the above $\W$ via comultiplication, but is inspired by various dualities on Fock spaces in terms of free field realization with the actions of Clifford and Weyl algebras \cite{W99}. The $\U_X(\e)$-module, which we denote by $\W^{2}$, specializes to a sum of fundamental representations of type $C_n$ when $(X,\e)=(C,(1^n))$ \cite{K14}, and (fundamental type) $q$-oscillator representations of type $D_n$ when $(X,\e)=(D,(0^n))$ \cite{KLO24}.
We show that $(\W^2)^{\ot \ell}$ has a $(\U_X(\e),{\bf U}^\imath(\mf{sp}_{2\ell}))$-module structure, where ${\bf U}^\imath(\mf{sp}_{2\ell})$ is the coideal subalgebra of ${\bf U}(\mf{sl}_{2\ell})$ corresponding to the symmetric pair $(\mf{sl}_{2\ell},\mf{sp}_{2\ell})$ whose Satake diagram is of type AII.

The above two cases can be summarized as follows:
\begin{equation}\label{eq:duality diagram}
\begin{tikzcd}[ampersand replacement=\&] % [row sep=small]
\U_X(\e) \arrow[d, phantom, sloped, "\supseteq"] \& \& \\
\U_A(\e) \arrow[r, phantom, sloped, "\curvearrowright"] \& \W^{\ot \ell} \& {\bf U}(\mf{sl}_\ell) \arrow[l, phantom, sloped, "\curvearrowleft"] \\
 \& \&  {\bf U}^\imath(\mf{so}_\ell) \arrow[u, phantom, sloped, "\subseteq"]
\end{tikzcd}
\quad\quad
\begin{tikzcd}[ampersand replacement=\&] % [row sep=small]
\U_X(\e) \arrow[d, phantom, sloped, "\supseteq"] \& \& \\
\U_A(\e) \arrow[r, phantom, sloped, "\curvearrowright"] \& (\W^2)^{\ot \ell} \& {\bf U}(\mf{sl}_{2\ell}) \arrow[l, phantom, sloped, "\curvearrowleft"] \\
 \& \&  {\bf U}^\imath(\mf{sp}_{2\ell}) \arrow[u, phantom, sloped, "\subseteq"]
\end{tikzcd}
\end{equation}
Finally, we show that the actions of
\begin{equation}\label{eq:commuting action}
\text{$(\U_{X}(\e),{\bf U}^\imath(\mf{so}_\ell))$ on $\W^{\ot \ell}$\quad 
and \quad $(\U_{X}(\e),{\bf U}^\imath(\mf{sp}_{2\ell}))$ on $(\W^2)^{\ot \ell}$}
\end{equation}
are semisimple, and describe its decomposition (Theorem \ref{thm:main result}), whose classical limit recovers \eqref{eq:spo duality} and hence \eqref{eq:O duality}-\eqref{eq:Sp duality} by specializing $\e=(0^n),(1^n)$ in each case. 
We use the notion of classical weight modules of $\imath$quantum groups introduced in \cite{Wat21} to prove the semisimplicity of $\W^{\ot \ell}$ and $(\W^2)^{\ot \ell}$ as a ${\bf U}^\imath(\mf{k})$-module with $\mf{k}=\mf{so}_\ell, \mf{sp}_{2\ell}$.

\subsection{}
Let us give some remarks related to the results in this paper.
We should first remark that $q$-analogues of \eqref{eq:O duality} and \eqref{eq:Sp duality} have  been obtained in \cite{ST}, where the action of $\U_X(\e)$ with $\e=(0^n), (1^n)$ is defined and verified to commute with that of ${\bf U}^\imath(\mf{k})$ by means of diagrammatic calculation generalizing the notion of type $A$ in \cite{CKM}. The result in \cite{ST}, as the authors mentioned, does not seem to yield an immediate generalization to the orthosymplectic case (see \cite{TVW} for type $A$).
The case of $(\mf{so}_{2n},O_\ell)$ on $\Lambda(\C^\ell \ot \C^{n})$ in \eqref{eq:O duality} is also studied in \cite{Ab,Wenzl,Wenzl2}, where a $q$-analogue of the action of $(\mf{so}_{2n},O_\ell)$ is described in terms of $q$-deformed Clifford algebras. To have a multiplicity-free decomposition, the centralizer of the Drinfeld-Jimbo's quantum group of $\mf{so}_{2n}$ on the tensor power of its spin representation is obtained in \cite{Wenzl,Wenzl2}, while the joint highest weight vectors are constructed in \cite{Ab}.
Compared to these results, our description of the action of $(\U_X(\e),{\bf U}^\imath(\mf{k}))$ in \eqref{eq:commuting action} is explicit even in case of \eqref{eq:O duality}-\eqref{eq:Sp duality}, and obtained in a direct and uniform way based on $\U_X(\e)$-modules $\W$ and $\W^{2}$ of fundamental types.

There is an isomorphism between $\U_X(\e)$ and the quantum superalgebra $U_X(\e)$ by Yamane \cite{Ya94} under mild extension (but not an isomorphism of Hopf algebras). This enables us to consider the representations of $\U_X(\e)$ and (super) representations of $U_X(\e)$ in an equivalent way (but not in a tensor category). There are several advantages of using $\U_X(\e)$ instead of $U_X(\e)$ in our context, while we still need to consider $\W^{\ot \ell}$ and $(\W^2)^{\ot \ell}$ as a $U_X(\e)$-module when we take its classical limit. 

One novel point of using $\U_X(\e)$ is that the comultiplication of $\U_X(\e)$ is given by the same formula of the usual quantum group so that the action of $(\U_X(\e),{\bf U}^\imath(\mf{k}))$ in \eqref{eq:commuting action} can be described and verified more easily. 
This may explain why it has been difficult to obtain a $q$-analogue of \eqref{eq:spo duality}, even of \eqref{eq:O duality} and \eqref{eq:Sp duality} for $q$-deformed exterior algebras which can be explained in terms of usual quantum groups with parameter $\tilde{q}=-q^{-1}$. 
%A
Moreover, by using $\U_X(\e)$, it is more natural to explain the connection with the representations of usual quantum groups from a viewpoint of super duality \cite{CLW}. Indeed, one can construct exact tensor functors naturally connecting $q$-oscillator representations of $\U_D(\e)$ (or $\U_C(\e)$) to finite-dimensional $\U_X(1^n)$-modules and infinite-dimensional $q$-oscillator $\U_Y(0^n)$-modules for $(X,Y)=(C,D), (D,C)$. These functors play an important role in the study of representations of quantum affine orthosymplectic algebras \cite{KLO24}.

In \cite{SW25}, there is a super analogue of a quantum symmetric pair, which is defined in terms of $U_A(\e)$ with an $\imath$quantum supergroup inside it. It would be interesting to find an analogue of a symmetric pair using $\U_A(\e)$  and characterize the centralizer of the associated $\imath$quantum supergroup in \eqref{eq:duality diagram}.

There exist crystals of an irreducible $q$-oscillator $\U_X(\e)$-module in $\W^{\ot \ell}$ or $(\W^2)^{\ot \ell}$ with $\e=(1^M,0^N)$ \cite{K15,K16} and a finite-dimensional irreducible ${\bf U}^\imath(\mf{so}_\ell)$-module with integral weights \cite{Wat23,Wat23-2}.
So it would be interesting to have a combinatorial picture of \eqref{eq:commuting action}, for example, an analogue of RSK algorithm in terms of $(\U_X(\e),{\bf U}^\imath(\mf{so}_\ell))$-crystals.
Also, it is natural to ask whether we can construct a canonical basis of an irreducible $q$-oscillator $\U_X(\e)$-module using an $\imath$canonical basis of ${\bf U}^\imath(\mf{so}_\ell)$-modules \cite{BW18,Wat23-2} through \eqref{eq:commuting action} in the sense of \cite{BW18-2,FKK,SW23}.
%Another direction related to this work is to consider the limit of $\mc{W}^{\otimes l}$ as $l\rightarrow \infty$ as recently studied in the case of type $A$ \cite{KL24}. On the side of ${\bf U}^\imath(\mf{k})$, one may expect to have representations of ${\bf U}^\imath(\mf{so}_\infty)$ and ${\bf U}^\imath(\mf{sp}_\infty)$ as $q$-analogues of those studied in \cite{DPS,PSt}.

\subsection{}
The paper is organized as follows. In Section \ref{sec:Qsuper}, we recall the quantized enveloping algebra $\U_X(\e)$ and its connection with $U_X(\e)$.
In Section \ref{sec:Fock space}, we introduce the $q$-oscillator representations $\W$ and $\W^{2}$ of $\U_X(\e)$, which are the building blocks of our duality.
In Section \ref{sec:QSP}, we prove the existence of $(\U_X(\e),{\bf U}^\imath(\mf{k}))$-action on $\W^{\ot \ell}$ and $(\W^2)^{\ot \ell}$ \eqref{eq:duality diagram}, and describe its semisimple decomposition. The well-definedness of the actions of $(\U_A(\e),{\bf U}(\mf{sl}_r))$ and $(\U_X(\e),{\bf U}^\imath(\mf{k}))$ are proved in Appendices \ref{app:proof type A} and \ref{app:proof type DC}, respectively.

\section{Quantum orthosymplectic superalgebras}\label{sec:Qsuper} 

\subsection{Generalized quantum groups of type $C$ and $D$}\label{subsec:pre}
We denote by $\mathbb{Z}_{\geq 0}$ the set of nonnegative integers. 
Let $\Bbbk=\mathbb{C}(q^{\frac{1}{2}})$ with $q^{\frac{1}{2}}$ a formal variable. %\comments{Do we need ``half-square"?} and $\Bbbk^{\times}=\Bbbk \setminus \left\{0\right\}$, where $q$ is an indeterminate.
We put 
{\allowdisplaybreaks
\begin{align*}
[m]&=\frac{q^m-q^{-m}}{q-q^{-1}}\quad (m\in \Z_{\geq 0}),\\
[m]!&=[m][m-1]\cdots [1]\quad (m\geq 1),\quad [0]!=1,\\
\begin{bmatrix} m \\ k \end{bmatrix}&= \frac{[m][m-1]\cdots [m-k+1]}{[k]!}\quad (0\leq k\leq m).
\end{align*}}
Fix $n\in \Z_{\ge 0}$ with $n\ge 4$.
We assume the following notations throughout the paper: 

\begin{itemize}
 
\item[$\bullet$] $\mathbb{I}= \{\,1<2<\cdots <n\,\}$,  

\item[$\bullet$] $\e=(\e_1,\dots,\e_n)$: a sequence with $\e_i\in \{0,1\}$ for $i \in \mathbb{I}$,

\item[$\bullet$] $P = \bigoplus_{i\in \I}\Z\de_i$: a $\Z$-lattice with a symmetric bilinear form $(\,\cdot\,|\,\cdot\,)$ satisfying 
\begin{equation*}
(\de_i|\de_j)=(-1)^{\e_i}\delta_{ij}, \quad 
%(\La|\La)=0, \quad
\end{equation*}

\item[$\bullet$] $q_i=\si q^{\si}$ $(i\in \mathbb{I})$, that is, $q_i=q$ if $\e_i=0$ and $q_i=\td{q}:=-q^{-1}$ if $\e_i=1$,

\item[$\bullet$] ${\bq}(\,\cdot\,,\,\cdot\,)$: a symmetric biadditive function from $P\times P$ to $\Bbbk^{\times}$ given by
\begin{equation*}
\begin{split}
&\bq(\mu,\nu) = \prod_{i\in\I}q_i^{\mu_i \nu_i} \quad (\mu,\nu\in P),
\end{split}
\end{equation*}
for $\mu=\sum_{i\in\I}\mu_i\de_i$ and $\nu=\sum_{i\in\I}\nu_i\de_i$,

\item[$\bullet$] $\Lambda = (-1)^{\e_1}\de_1 + \cdots + (-1)^{\e_n} \de_n$.
\end{itemize}
\medskip

Let us first recall the notion of a {\em generalized quantum group $\U_X(\e)$ of type $X=C, D$ associated to $\e$} \cite{KOS,Ma}. Put
\begin{itemize}
\item[$\bullet$] $I=\{0,1,\ldots,n-1\}$,

\item[$\bullet$] $I_{\rm even}=\{\,i\in I\,|\,(\alpha_i|\alpha_i)\neq 0\,\}$ and 
$I_{\rm odd}=\{\,i\in I\,|\,(\alpha_i|\alpha_i)=0\,\}$,

\item[$\bullet$] $\alpha_i=\de_{i}-\de_{i+1} \in P$ for $i\in I \setminus\{0\}$,

\item[$\bullet$] $\alpha_0=-\de_1-\de_2$ and $-2\de_1$ for $X=D$ and $C$, respectively.

%\item[$\bullet$] $I_{\rm even}=\{\,i\in I\,|\,(\alpha_i|\alpha_i)\neq 0\,\}$, $I_{\rm odd}=\{\,i\in I\,|\,(\alpha_i|\alpha_i)=0\,\}$. 
\end{itemize}
Note that there is no non-isotropic odd root in our setting.

We first define $\mathcal{U}_D(\epsilon)$ to be the associative $\Bbbk$-algebra with $1$ generated by 
$e_i,f_i,k_\mu$  $(i\in I, \mu\in P)$ satisfying
{\allowdisplaybreaks
\begin{equation}\label{eq:rel for GQ-D-1}
\begin{split}
&k_\mu=1 \quad(\mu=0), \quad k_{\mu +\mu'}=k_{\mu}k_{\mu'} \quad (\mu, \mu' \in P), \\ 
&k_\mu e_i k_{-\mu}=\bq(\mu,\alpha_i)e_i,\quad 
k_\mu f_i k_{-\mu}=\bq(\mu,\alpha_i)^{-1}f_i\quad (i\in I, \mu\in P), \\ 
&e_if_j - f_je_i =\delta_{ij}\frac{k_{i} - k_{i}^{-1}}{q-q^{-1}}\quad (i,j\in I), \\
&[e_i ,e_j]=[f_i ,f_j]=0 \quad\text{if $(\alpha_i|\alpha_j)=0$},\\
&e_i^2= f_i^2 =0 \quad \text{($i\in I_{\rm odd}$)}, \\
&e^2 _i e_j - (-1)^{\e_i}[2]e_i e_j e_i +e_j e^2 _i =(e \rightarrow f)=0
\quad \text{if $i,j\in I\setminus\{0\}$, $|i-j|=1$, $\e_i=\e_{j}$}, \\
&e_{i}e_{i-1}e_{i}e_{i+1}-e_{i}e_{i+1}e_{i}e_{i-1}+(-1)^{\e_i}[2]e_{i}e_{i-1}e_{i+1}e_{i} -e_{i-1}e_{i}e_{i+1}e_{i}\\
&\quad\quad\quad\quad\quad\quad +e_{i+1}e_{i}e_{i-1}e_{i}=(e \rightarrow f)=0
\quad\text{if $i\in I\setminus \{0,1,n-1\}$, $\e_i\neq \e_{i+1}$},\\
\end{split}
\end{equation}}
and
{\allowdisplaybreaks
\begin{equation}\label{eq:rel for GQ-D-2}
\begin{split}
&e_0^2 e_2 - (-1)^{\e_2}[2]e_0 e_2 e_0 +e_2 e_0^2 =(e \rightarrow f)=0 
\quad \text{if $\e_1=\e_2$},\\
&e_2^2 e_0 - (-1)^{\e_2}[2]e_2 e_0 e_2 +e_0 e_2^2 =(e \rightarrow f)=0 
\quad \text{if $\e_2=\e_3$}, \\
&e_0e_1e_2-e_1e_0e_2+(-1)^{\e_2}[2](e_1e_2e_0-e_0e_2e_1) \\
&\quad\quad\quad\quad\quad\quad\quad\quad\quad\quad
+e_2e_0e_1-e_2e_1e_0=(e\rightarrow f)=0 \quad \text{if $\e_1\neq \e_2$}, \\
&e_0e_2e_3e_2-e_3e_2e_0e_2+(-1)^{\e_3}[2]e_2e_3e_0e_2 \\
&\quad\quad\quad\quad\quad\quad\quad\quad\quad\quad
-e_2e_0e_2e_3+e_2e_3e_2e_0=(e \rightarrow f)=0
\quad \text{if $\e_2\neq \e_3$},
\end{split}
\end{equation}
where $k_i=k_{\alpha_i}$ for $i\in I$, and $(e\rightarrow f)$ denotes the relation obtained by replacing $e_i$ with $f_i$.}

Next, we define $\mathcal{U}_C(\epsilon)$ to be the associative $\Bbbk$-algebra with $1$ generated by $e_i,f_i,k_\mu$ $(i\in I, \mu\in P)$ satisfying the relations \eqref{eq:rel for GQ-D-1} except $e_0f_0-f_0e_0=\frac{k_0-k_0^{-1}}{q^2-q^{-2}}$ together with %\comments{Check the relations are correctly presented (24.02.20)}
{\allowdisplaybreaks
\begin{equation}\label{eq:rel for GQ-C-2}
\begin{split}
& e_0^2e_1-(q^2+q^{-2})e_0e_1e_0+e_1e_0^2=(e \rightarrow f)=0, \\
& e_1^3e_0-[3]e_1^2e_0e_1+[3]e_1e_0e_1^2-e_0e_1^3=(e \rightarrow f)=0 \quad\text{if } \e_1=\e_2, \\
&[[[[e_2, e_1]_{q_2}, e_0]_{q_1^2},[e_2, e_1]_{q_2}], e_1]=(e\rightarrow f)=0 \quad \text{if $\epsilon_1\neq \epsilon_2 \text{ and } \e_2\neq \epsilon_3$}, \\
&[[[[[[e_3, e_2]_{q_3}, e_1]_{q_2}, e_0]_{q_1^2}, e_1]_{q_1}, e_2]_{q_2}, e_1]=(e\rightarrow f)=0 \quad \text{if $\epsilon_1\neq \epsilon_2 \text{ and } \e_2= \epsilon_3$},
\end{split}
\end{equation}
where $[u, v]_t=uv-tvu$ and $[u,v]=[u,v]_1$}.

Note that $\U_X(\e)$ coincides with the usual quantum group $U_q(X_n)$ of type $X_n$ when $\e$ is homogeneous, that is, 
%\begin{equation}\label{eq:homogeneous}
% \U_X(0^n)=U_q(X_n),\quad \U_X(1^n)=U_{\td{q}}(X_n)
%\end{equation}
$\U_X(0^n)=U_q(X_n)$ and $\U_X(1^n)=U_{\td{q}}(X_n)$.
Let $\e^c=(1-\e_1,\dots,1-\e_n)$. Then there is an isomorphism $\td{\cdot}: \U_X(\e)\rightarrow \U_X(\e^c)$ such that $\td{q}=-q^{-1}$, and $\td{e_i}=e_i$, $\td{f_i}=f_i$ and $\td{k_\mu}=k_\mu$ for $i\in I$ and $\mu\in P$. There is a Hopf algebra structure on $\U_X(\e)$ with the comultiplication $\Delta$ and the antipode $S$ given by
\begin{equation}\label{eq:comult}
\begin{split}
\Delta(k_\mu)&=k_\mu\ot k_\mu, \\ 
\Delta(e_i)&= e_i\ot 1 + k_i \ot e_i, \\
\Delta(f_i)&= f_i\ot k_i^{-1} + 1 \ot f_i , \\  
S(k_\mu)=k_{-\mu},& \ \ S(e_i)=-k_i^{-1}e_i, \ \  S(f_i)=-f_ik_i,\\
\end{split}
\end{equation}
for $\mu\in P$ and  $i\in I$. 
Let $\U_A(\e)$ denote the subalgebra of $\U_X(\e)$ generated by $e_i,f_i,k_\mu$  $(i\in I\setminus\{0\}, \mu\in P)$. In particular, $\U_A(0^n)$ is the quantized enveloping algebra for $\gl_n$, which we denote by $U_q(\gl_n)$.

We say that a $\U_X(\e)$-module $V$ has a weight space decomposition if $V=\bigoplus_{\la\in {P}}V_\la$, where $V_\la = \{\,v\in V\,|\,k_{\mu} v= {\bq}(\la,\mu) v \ \ (\mu\in P) \,\}$.

%\begin{rem}
%\red{If $\e_{n-1}\neq \e_n$, there is an algebra isomorphism \cite[Sec 5.3]{Ma}
%$\tau_{n-1} : \U_C(\e_1, \ldots, \e_{n-1}, \e_n) \rightarrow \U_D(\e_1, \ldots, \e_n, \e_{n-1})$.}\comments{We have modified the defining relations for $\U_C(\e)$. So we should check whether we still have $\tau$. Also, following our convention, $\tau_{n-1}$ should be $\tau_{1}$. }
%\end{rem}

\subsection{Quantum superalgebras of type $C$ and $D$}\label{subsec:pre-2}
We briefly review how $\U_X(\e)$ is related to the quantized enveloping algebra for an orthosymplectic Lie superalgebra introduced in \cite{Ya94}.

Let $U_D(\epsilon)$ be an associative $\Bbbk$-algebra with 1 generated by $E_{i}$, $F_{i}$, $K_{\mu}$ $(i\in I, \mu\in P)$ satisfying %\comments{Check the relations are correctly presented (24.02.20)}

\allowdisplaybreaks{\small{}
\begin{equation}\label{eq:rel for YQ-D-1}
\begin{split}
& K_{\mu}=1 \quad(\mu=0),\quad K_{\mu+\mu'}=K_{\mu}K_{\mu'} \quad (\mu, \mu' \in P),\\
& K_{\mu}E_{i}K_{\mu}^{-1}=q^{(\mu|\alpha_{i})}E_{i},\quad K_{\mu}F_{i}K_{\mu}^{-1}=q^{-(\mu|\alpha_{i})}F_{i} \quad (i\in I, \mu\in P),\\
& E_{i}F_{j}-(-1)^{p(i)p(j)}F_{j}E_{i}=(-1)^{\epsilon_{i}}\delta_{ij}\frac{K_{i}-K_{i}^{-1}}{q-q^{-1}} \quad (i,j\in I\setminus\{0\}),\\
& E_{0}F_{0}-(-1)^{p(0)p(0)}F_{0}E_{0}=(-1)^{\epsilon_{1}}\frac{K_{0}-K_{0}^{-1}}{q-q^{-1}},\\
& E_{i}^{2}=F_{i}^{2}=0\quad\text{($i\in I_{\rm odd}$)},\\
&  E_{i}E_{j}-(-1)^{p(i)p(j)}E_{j}E_{i}=(E\rightarrow F)=0\quad\text{if $(\alpha_i | \alpha_j)=0$},\\
& E_{i}^{2}E_{j}-[2]E_{i}E_{j}E_{i}+E_{j}E_{i}^{2}=(E\rightarrow F)=0\quad \text{if $i,j\in I\setminus\{0\}$, $|i-j|=1$, and $\e_i=\e_{i+1}$},\\
& [E_{i},[[E_{i-1},E_{i}]_{(-1)^{p(i-1)}q},E_{i+1}]_{(-1)^{(p(i-1)+p(i))p(i+1)}q^{-1}}]_{(-1)^{p(i-1)+p(i)+p(i+1)}}\\
& \quad\quad\quad\quad\quad\quad\quad\quad\quad\quad =(E\rightarrow F)=0\quad\text{ if $\epsilon_{i}\neq\epsilon_{i+1}$ and $i\in I\setminus\{0,1,n-1\}$},\\
\end{split}
\end{equation}
and
\begin{equation}\label{eq:rel for YQ-D-2}
\begin{split}
& E_{2}^{2}E_{0}-[2]E_{2}E_{0}E_{2}+E_{0}E_{2}^{2}=(E\rightarrow F)=0\quad\text{ if }\epsilon_{2}=\epsilon_{3},\\
& E_{0}^{2}E_{2}-[2]E_{0}E_{2}E_{0}+E_{2}E_{0}^{2}=(E\rightarrow F)=0\quad\text{ if }\epsilon_{1}=\epsilon_{2},\\
& [E_{2},[[E_{3},E_{2}]_{(-1)^{p(3)}q},E_{0}]_{(-1)^{(p(3)+p(2))p(0)}q^{-1}}]_{(-1)^{p(3)+p(2)+p(0)}}= (E\rightarrow F)=0\quad\text{ if }\epsilon_{2}\neq\epsilon_{3},\\
& E_{0}E_{1}E_{2}-E_{1}E_{0}E_{2}+(-1)^{\epsilon_{2}+\epsilon_{3}}(E_{1}E_{2}E_{0}-E_{0}E_{2}E_{1}) +E_{2}E_{0}E_{1}-E_{2}E_{1}E_{0} \\
& \quad\quad\quad\quad\quad\quad\quad\quad\quad\quad\quad\quad\quad\quad\quad\quad
=(E\rightarrow F)=0\quad\text{ if \ensuremath{\epsilon_{1}\neq\epsilon_{2}}},
\end{split}
\end{equation}
} where $K_i=K_{\alpha_i}$ for $i\in I$, $p(i)=0$ (resp. $1$) for $i\in I_{\rm even}$ (resp. $i\in I_{\rm odd}$).

Similarly, let $U_C(\epsilon)$ be the associative $\Bbbk$-algebra with 1 generated by $E_{i}$, $F_{i}$, $K_{\mu}$ $(i\in I, \mu\in P)$ satisfying the relations \eqref{eq:rel for YQ-D-1} except $E_{0}F_{0}-F_{0}E_{0}=(-1)^{\epsilon_{1}}\frac{K_{0}-K_{0}^{-1}}{q^2-q^{-2}}$ together with %\comments{Check the relations are correctly presented (240220)}

\allowdisplaybreaks{\small{}
\begin{equation}\label{eq:rel for YQ-C-2}
\begin{split}
& E_1^3E_0-[3]E_1^2E_0E_1+[3]E_1E_0E_1^2-E_0E_1^3=(E\rightarrow F)=0 \quad\text{if } \e_1=\e_2, \\
& E_0^2E_1-(q^2+q^{-2})E_0E_1E_0+E_1E_0^2=(E\rightarrow F)=0, \\
& [[[[E_2, E_1]_{-q^{(-1)^{\e_2}}}, E_0]_{q^{2(-1)^{\e_1}}},[E_2, E_1]_{-q^{(-1)^{\e_2}}}], E_1]=(E\rightarrow F)=0 \quad \text{if } \e_1\neq\e_2 \text{ and } \e_2\neq\e_3,\\
& [[[[[[E_3, E_2]_{q^{(-1)^{\e_3}}}, E_1]_{(-1)^{p(3)}q^{(-1)^{\e_2}}}, E_0]_{q^{2(-1)^{\e_1}}}, E_1]_{(-1)^{(p(3)+1)}q^{(-1)^{\e_1}}}, E_2]_{q^{(-1)^{\e_2}}}, E_1]_{(-1)^{p(3)}}\\
& \quad\quad\quad\quad\quad\quad\quad\quad\quad\quad\quad\quad\quad\quad\quad\quad =(E\rightarrow F)=0 \quad \text{if $\e_1\neq\e_2$ and $\e_2=\e_3$}.
\end{split}
\end{equation}
}

Then there is an isomorphism $\tau$ of $\Bbbk$-algebras between suitable extensions of $\U_X(\e)$ and $U_X(\e)$ ($X=C,D$) (cf.~\cite{KLO24}).
Let $\Sigma$ be the $\Bbbk$-bialgebra generated by $\sigma_i$ for $i\in \mathbb I$ such that $\sigma_i \sigma_j = \sigma_j \sigma_i$ and $\sigma_i^2=1$ with comultiplication given by $\Delta(\sigma_i)=\sigma_i \ot \sigma_i$ for $i \in \mathbb{I}$. Define the action of $\Sigma$ on $\U_X(\e)$ by
$$\sigma_j k_\mu = k_\mu, \quad \sigma_j e_i = (-1)^{\e_j (\de_j | \alpha_i)} e_i, \quad \sigma_j f_i = (-1)^{\e_j (\de_j | \alpha_i)}f_i,$$
for $i \in I$, $j\in \mathbb{I}$ and $\mu \in P$.
Let $\U_X(\e)[\sigma]$ be the semidirect product of $\U_X(\e)$ and $\Sigma$ with respect to this action. Define $U_X(\e)[\sigma]$ in the same manner. 

  Let $\mathbb{I}=\mathbb{I}^{(1)} \sqcup \cdots \sqcup \mathbb{I}^{(l)}$ be a partition of $\mathbb{I}$ such that
\begin{itemize}
\item[(1)] $\e_i$ is constant on $\mathbb{I}^{(k)}$ for $1 \leq k \leq l$,
\item[(2)] $\e_i \neq \e_j$ if $i\in \mathbb{I}^{(k)}$ and $j\in \mathbb{I}^{(k+1)}$ for $1 \leq k \leq l-1$,
\item[(3)] $\mathbb{I}^{(k)}=\{i_{k-1}+1, i_{k-1}+2, \ldots, i_k\}$, where $i_k=\text{max }\mathbb{I}^{(k)}$ for $1 \leq k \leq l$ and $i_0=0$.
\end{itemize}
Let $\sigma_{\leq j}=\sigma_1 \sigma_2 \cdots \sigma_j$ for $j\in \mathbb{I}$ and {$\upvarsigma_i =\sigma_i \sigma_{i+1}$} for $i\in I\setminus\{0\}$. Define $\tau(E_i), \tau(F_i)$, and $\tau(K_i)$ for $i\in I$ as follows: 

\begin{itemize}
\item[(1)] if $i \in I_{\text{even}}\setminus\{0\}$ with $(\e_i, \e_{i+1})=(0, 0)$, then we put
$$\tau(E_i)=e_i, \quad \tau(F_i)=f_i, \quad \tau(K_i)=k_i,$$

\item[(2)] if $i \in I_{\text{odd}}\setminus\{0\}$ with $(\e_i, \e_{i+1})=(0,1)$, then we put
$$\tau(E_i)=e_i\sigma_{\leq i}, \quad \tau(F_i)=f_i \sigma_{\leq i}\upvarsigma_i, \quad \tau(K_i)=k_i \upvarsigma_i, $$

\item[(3)] if $i\in I_{\text{even}}\setminus\{0\}$  with $(\e_i, \e_{i+1})=(1, 1)$, then we have $\{i, i+1\} \subset \mathbb{I}^{(k)}$ for a unique $k$ and put
$$\tau(E_i)=e_i\upvarsigma_i^{i-i_{k-1}}, \quad \tau(F_i)=-f_i\upvarsigma_i^{i-i_{k-1}-1}, \quad \tau(K_i)=k_i\upvarsigma_i,$$

\item[(4)] if $i\in I_{\text{odd}}\setminus\{0\}$ with $(\e_i, \e_{i+1})=(1, 0)$, then we have $i=i_k$ for a unique $k$ and put
$$\tau(E_i)=e_i\sigma_{\leq i}\upvarsigma_i^{i_k-i_{k-1}}, \quad \tau(F_i)=(-1)^{i_k-i_{k-1}}f_i \sigma_{\leq i}\upvarsigma_i^{i_k-i_{k-1}-1}, \quad \tau(K_i)=k_i\upvarsigma_i,$$

\item[(5)] if $i=0$ and $X=D$, we define $\tau(E_0), \tau(F_0)$, and $\tau(K_0)$  in the same way as $i=1$,

\item[(6)] if $i=0$ and $X=C$, we put $\tau(F_0)=f_0$, $\tau(K_0)=k_0$, and $\tau(E_0)=(-1)^{\e_1}e_0$.%\comments{Check the well-definedness (240220)}
\end{itemize}

Let $\ms{X}$ (resp. $\ms{Y}$) be the subalgebra of {$U_X(\e)[\sigma]$ (resp. $\U_X(\e)[\sigma]$)} generated by $E_i,F_i,K_i$ and $\sigma_j$ (resp. $e_i,f_i,k_i$ and $\sigma_j$) for $i\in I$ and $j\in \I$.

\begin{prop}\label{prop:iso GQ and YQ}
The map $\tau$ extends to an isomorphism of $\Bbbk$-algebras from $\ms{X}$ to $\ms{Y}$ such that $\tau(\sigma_j)=\sigma_j$ for $j\in \I$.
\end{prop}
\pf The proof for the case when $X=D$ is given in \cite[Proposition 2.1]{KLO24}). So it suffices to consider the case when $X=C$ and check the relations where $e_0$ is involved.

Let us prove that the last two relations in \eqref{eq:rel for GQ-C-2} and \eqref{eq:rel for YQ-C-2} are preserved under $\tau$ up to $\Sigma$ modulo the other relations since the proof for the first two relations are easy. For simplicity, we use the following notations:
\begin{equation}\label{eq:monomial in U} 
\begin{split}
    (i_1 i_2 \cdots i_l)&:=E_{i_1}E_{i_2}\cdots E_{i_l} \in U_C(\e), \\
    \langle i_1 i_2 \cdots i_l \rangle&:=e_{i_1}e_{i_2} \cdots e_{i_l} \in \U_C(\e). 
\end{split}
\end{equation}

\noindent{\em Case 1.} Consider the second last relation. First, we have
{\allowdisplaybreaks\begin{align*}
    A_1&:=[[[[E_2, E_1]_{-q^{(-1)^{\e_2}}}, E_0]_{q^2},[E_2, E_1]_{-q^{(-1)^{\e_2}}}], E_1]\\
    &=q^{(-1)^{\e_2}}(210121)+q^{2(-1)^{\e_2}}(120121)-(012121)-(212101)-q^{2(-1)^{\e_2}}(121201)\\
    &\quad +q^{(-1)^{\e_1}}(210121)+(120121) -(121021)-q^{(-1)^{\e_2}}(121012)+q^{2(-1)^{\e_1}}(102121) \\
    &\quad +(101212) +(121210) -q^{2(-1)^{\e_1}}(121021)-q^{(-1)^{\e_1}}(121012).
\end{align*}}
On the other hand, we have
{\allowdisplaybreaks\begin{align*}
    B_1&:=[[[[e_2, e_1]_{q_2}, e_0]_{q_2^{-2}},[e_2, e_1]_{q_2}], e_1]\\
    &= -q_2 \langle 210121 \rangle +q_2^2 \langle 120121 \rangle - \langle 012121 \rangle -\langle 212101 \rangle  -q_2^2 \langle 121201 \rangle \\
    &\quad -q_2^{-1} \langle 210121 \rangle + \langle 120121 \rangle -\langle 121021 \rangle +q_2 \langle 121012 \rangle  +q_2^{-2} \langle 102121 \rangle \\
    &\quad + \langle 101212 \rangle +\langle 121210 \rangle -q_2^{-2} \langle 121021 \rangle  +q_2^{-1} \langle 121012 \rangle.
\end{align*}}
Here we use the relations $e_1^2=e_2^2=0$ and $E_1^2=E_2^2=0$.

\noindent{\em Case 1-1}. Assume that $(\e_1, \e_2, \e_3)=(0,1,0)$. Then
\begin{equation*}
    \tau(E_0)=e_0,\quad \tau(E_1)=e_1\sigma_1, \quad \tau(E_2)=e_2\sigma_1\sigma_3.
\end{equation*}
Since $\e_1=\e_3=0$, $\sigma_1$ and $\sigma_3$ commute with all the generators of $\U_C(\e)[\sigma]$ and $U_C(\e)[\sigma]$. Thus, $\tau( (i_1 \cdots i_6))= \langle i_1 \cdots i_6 \rangle \sigma_3$ for all $(i_1 \cdots i_6)\in U_C(\e)$. Thus, $\tau(A_1)=B_1 \sigma_3$. %Since $B_1=0$ in $\U_C(\e)[\sigma]$, we have $\tau(A_1)=0$.
\medskip

\noindent{\em Case 1-2}.
Assume that $(\e_1, \e_2, \e_3)=(1,0,1)$. Then
  \begin{equation*}
    \tau(E_0)=e_0, \quad \tau(E_1)=e_1\sigma_2, \quad \tau(E_2)=e_2\sigma_1\sigma_2.
  \end{equation*}
Since $\e_2=0$, $\sigma_2$ is central in $\U_C(\e)[\sigma]$. Since $\e_1=1$, we have $\sigma_1 e_2=e_2 \sigma_1$, $\sigma_1 e_0=e_0\sigma_1$ and $\sigma_1 e_1=-e_1 \sigma_1$. So we have $\tau((i_1 \ldots i_6))=\pm \langle i_1 \ldots i_6 \rangle$ where the sign is determined by the arrangement of indices $i_1, \ldots, i_6$. More precisely, note that the sign changes whenever $\sigma_1$ is exchanged with $e_1$. Since $\sigma_1^2=1$ and $\sigma_1$ does not appear in the images of $\tau(E_0)$ and $\tau(E_1)$, it suffices to count the number of 1's lying between two 2's in each term of $A_1$. If there are even number of 1's, then $\tau((i_1 \ldots i_6))=\langle i_1 \ldots i_6 \rangle$ and if there are odd number of 1's, then $\tau((i_1 \ldots i_6))=-\langle i_1 \ldots i_6 \rangle$. For example, we have $\tau((210121))=(-1)^2\langle 210121 \rangle \sigma_2=\langle 210121 \rangle \sigma_2$ and $\tau((120121))=-\langle 120121 \rangle \sigma_2$. Thus, we obtain
{\allowdisplaybreaks\begin{align*}
    \tau(A_1)&=\{q^{(-1)^{\e_2}}\langle 210121 \rangle-q^{2(-1)^{\e_2}}\langle 120121 \rangle +\langle 012121 \rangle +\langle 212101 \rangle +q^{2(-1)^{\e_2}}\langle 121201\rangle\\
    &\quad +q^{(-1)^{\e_1}}\langle 210121 \rangle-\langle 120121 \rangle  +\langle 121021 \rangle-q^{(-1)^{\e_2}}\langle 121012 \rangle -q^{2(-1)^{\e_1}}\langle 102121 \rangle \\
    &\quad -\langle101212\rangle  -\langle121210\rangle +q^{2(-1)^{\e_1}}\langle 121021 \rangle -q^{(-1)^{\e_1}}\langle 121012 \rangle \}\sigma_2 \\
    &=-B_1\sigma_2
\end{align*}}
\vskip 3mm

\noindent{\em Case 2.} Consider the last relation. First, we have
{\allowdisplaybreaks\begin{align*}
    A_2&:=[[[[[[E_3, E_2]_{q^{(-1)^{\e_3}}}, E_1]_{(-1)^{p(3)}q^{(-1)^{\e_2}}}, E_0]_{q^{2(-1)^{\e_1}}}, E_1]_{(-1)^{(p(3)+1)}q^{(-1)^{\e_1}}}, E_2]_{q^{(-1)^{\e_2}}}, E_1]_{(-1)^{p(3)}} \\
    &=(3210121)-(q^{(-1)^{\e_3}}+q^{(-1)^{\e_1}})(2310121)-(q^{(-1)^{\e_2}}+q^{3(-1)^{\e_1}})(3120121)\\
    &\quad +(-1)^{p(3)}(q^{2(-1)^{\e_2}}+q^{2(-1)^{\e_1}})(1230121) +q^{(-1)^{\e_1}}(3012121)-(-1)^{p(3)}(0123121) \\
    &\quad +q^{(-1)^{\e_1}}(3121021)-(-1)^{p(3)}(1231021) +q^{2(-1)^{\e_1}}(3101221)-q^{(-1)^{\e_1}}(1012321) \\
    &\quad  -(2312101)+(-1)^{p(3)}q^{(-1)^{\e_2}}(2123101) +(2101231)  -(3121012) \\
    &\quad +(-1)^{p(3)}(q^{(-1)^{\e_3}}+q^{(-1)^{\e_1}})(1231012) -q^{(-1)^{\e_1}}(3101212)+(1012312)\\
    &\quad +(-1)^{p(3)}q^{(-1)^{\e_2}}(1232101)-(-1)^{p(3)}q^{2(-1)^{\e_2}}(1223101) \\
    &\quad -(-1)^{p(3)}(q^{2(-1)^{\e_2}}+q^{2(-1)^{\e_1}})(1231201)  +(q^{3(-1)^{\e_2}}+q^{(-1)^{\e_1}})(1212301) \\
    &\quad  +(-1)^{p(3)}(1230121)-q^{(-1)^{\e_2}}(1201231) +(-1)^{p(3)}(1231210) \\
    &\quad -q^{(-1)^{\e_2}}(1212310) -(-1)^{p(3)}(1210123).
\end{align*}}
Here, we use the relations $E_1^2=0$, $[E_3, E_0]=[E_2, E_0]=0$ and $E_1E_3=(-1)^{p(3)}E_3E_1$ in $U_C(\e)[\sigma]$.  
On the other hand, we have
{\allowdisplaybreaks\begin{align*}
    B_2&:=[[[[[[e_3, e_2]_{q_3}, e_1]_{q_3}, e_0]_{q_3^{-2}}, e_1]_{-q_3^{-1}}, e_2]_{q_3}, e_1]\\
    &= \langle 3210121 \rangle -(q_3+q_3^{-1})\langle 2310121\rangle -(q_3+q_3^{-3})\langle 3120121\rangle +(q_3^2+1+q_3^{-2})\langle 1230121\rangle \\
    &\quad +q_3^{-1}\langle 3012121\rangle -\langle 0123121\rangle +q_3^{-1}\langle 3121201\rangle -(q_3^2+1+q_3^{-2})\langle 1231201\rangle \\
    &\quad +q_3^{-2}\langle 3101221\rangle -q_3^{-1}\langle 1012321\rangle -\langle 2312101\rangle +q_3\langle 2123101\rangle \\
    &\quad +\langle 2101231\rangle -\langle 3121012\rangle +(q_3+q_3^{-1})\langle 1231012\rangle -q_3^{-1}\langle 3101212\rangle \\
    &\quad  +\langle 1012312\rangle  +q_3\langle 1232101\rangle-q_3^2\langle 1223101\rangle +(q_3^3+q_3^{-1})\langle 1212301\rangle \\
    &\quad -q_3\langle 1201231\rangle +\langle 1231210\rangle -q_3\langle 1212310\rangle -\langle 1210123\rangle.
\end{align*}}
Similarly, we use here the relations $e_1^2=0$ and $[e_3,e_0]=[e_3,e_1]=[e_2,e_0]=0$.
There are four cases of $\e'=(\e_1, \e_2, \e_3, \e_4)$:
\begin{equation}\label{eq:case of e'}
 (1,0,0,0),\quad (1,0,0,1),\quad (0,1,1,0),\quad (0,1,1,1).
\end{equation}
By the direct computation as in {\em Case 1}, we have
\begin{equation*}
\tau(A_2)=B_2\sigma_2,\quad \tau(A_2)=-B_2\sigma_1\sigma_3,\quad \tau(A_2)=-B_2,\quad \tau(A_2)=B_2\sigma_1\sigma_2\sigma_3
\end{equation*}
for each case of $\e'=(\e_1, \e_2, \e_3, \e_4)$ in \eqref{eq:case of e'}, respectively.
% Since $B_2=0$ in $\U_C(\e)[\sigma]$, it follows that $\tau(A_2)=0$. 

By {\em Case 1} and {\em Case 2}, the claim is proved. Hence the proof completes.
\qed

\begin{rem}\label{rem:tau iso}
{\rm
The second last relation in \eqref{eq:rel for GQ-C-2} for $\U_C(\e)$ is slightly different from the one in \cite{Ma}. We adopt the one in \cite{Ya94} to verify the well-definedness of $\tau$. Also, we should remark that the comultiplication on $\ms{X}$ induced from $\Delta$ on $\ms{Y}$ under $\tau$ is not equal to the usual one on $\ms{X}$ in \cite{Ya94}. For example, assume that $(\e_1, \e_2)=(1,1)$. Then we have $\Delta(E_1)=E_1\ot 1 + K_1\ot E_1$. On the other hand, $((\tau^{-1}\ot \tau^{-1})\circ \Delta \circ \tau )(E_1)=E_1\ot \sigma_1\sigma_2 + K_1 \ot \sigma_1\sigma_2E_1$, whence $\Delta\neq(\tau^{-1}\ot \tau^{-1})\circ \Delta \circ \tau$ in this case. 
} 
\end{rem}

\subsection{Classical limit of $\U_X(\e)$-modules}\label{subsec:classical limit}
Let $\gl(V)$ denote the general linear Lie superalgebra for a complex superspace $V$.
Then $\mf{osp}(V)$ (resp. $\mf{spo}(V)$) is defined to be a subalgebra of $\gl(V)$ preserving a non-degenerate supersymmetric (resp. skew supersymmetric) bilinear form on $V$. We have $\mf{osp}(V)\cong \mf{spo}(\Pi V)$ as a Lie superalgebra, where $\Pi$ is the parity change functor (cf.~\cite[Remark 1.9]{CW12}).

Let $V=\C^{2n_0 | 2n_1}$ be a complex superspace such that $V_{\ov 0}\cong \C^{2n_0}$ and $V_{\ov 1}\cong \C^{2n_1}$, where $n_0$ and $n_1$ are the numbers of 0's and 1's in $\e$, respectively. 
We put 
\begin{equation}\label{eq:osp_X}
\mf{osp}_D(\e)=
\begin{cases} 
 \mf{spo}(V) & \text{if } {\e_1}=1, \\
 \mf{spo}(\Pi V) & \text{if } {\e_1}=0,
\end{cases}
\quad\quad 
\mf{osp}_C(\e)=\begin{cases} \mf{spo}(V) & \text{if } {\e_1}=0, \\
\mf{spo}(\Pi V) & \text{if } {\e_1}=1,
\end{cases}
\end{equation}
where the fundamental system for $\mf{spo}(V)$ (resp. $\mf{spo}(\Pi V)$) is assumed to be the one corresponding to the $\varepsilon\delta$-sequence determined by $\e$, that is, %\comments{`` with $\e_1=0$" is removed. (20250416)} 
a sequence whose $i$th coordinate is $\delta$ (resp. $\varepsilon$) if $\e_i=0$, and $\varepsilon$ (resp. $\delta$) if $\e_i=1$ (see \cite[Section 1.3]{CW12}, where the reading ordering of $\varepsilon\delta$-sequences is reverse to ours).

\begin{ex}{\rm
The following are examples of Dynkin diagram of $\mf{osp}_X(\e)$.

(1) For $\mf{osp}_D(\e)$, 

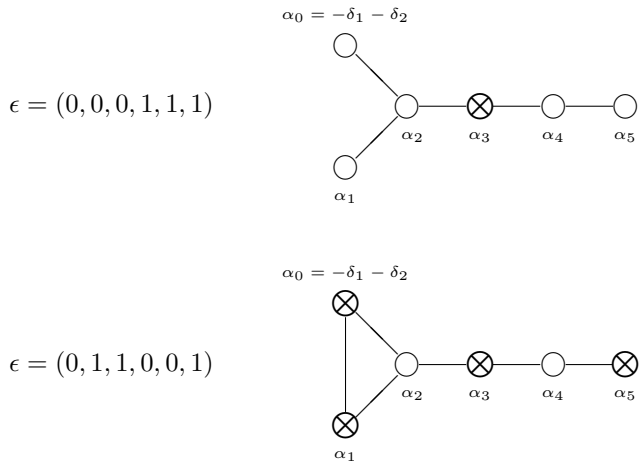
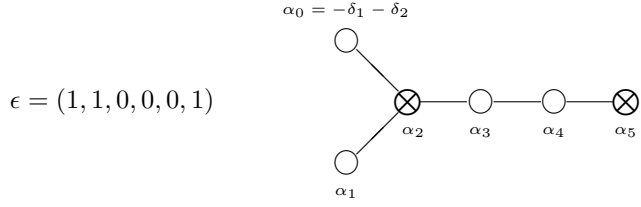
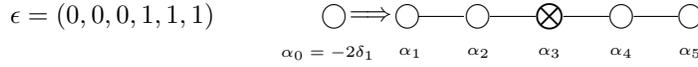
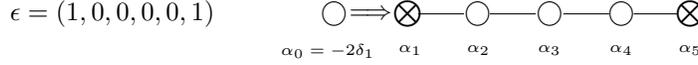
\begin{center} 
\hskip -1cm \setlength{\unitlength}{0.16in} \medskip
\begin{picture}(16,5.8)
\put(-5,1.8){$\e=(0, 0, 0, 1, 1, 1)$}
\put(6,0){\makebox(0,0)[c]{$\bigcirc$}}
\put(6,4){\makebox(0,0)[c]{$\bigcirc$}}
\put(8,2){\makebox(0,0)[c]{$\bigcirc$}}
\put(10.4,2){\makebox(0,0)[c]{$\bigotimes$}}
\put(12.8,2){\makebox(0,0)[c]{$\bigcirc$}}
\put(15.15,2){\makebox(0,0)[c]{$\bigcirc$}}
\put(6.35,0.3){\line(1,1){1.35}} 
\put(6.35,3.7){\line(1,-1){1.35}}
\put(8.4,2){\line(1,0){1.55}} 
\put(10.82,2){\line(1,0){1.55}}
\put(13.2,2){\line(1,0){1.55}} 
\put(6,5){\makebox(0,0)[c]{\tiny $\alpha_0=-\delta_1-\delta_2$}}
\put(6,-1){\makebox(0,0)[c]{\tiny $\alpha_1$}}
\put(8.2,1){\makebox(0,0)[c]{\tiny $\alpha_2$}}
\put(10.4,1){\makebox(0,0)[c]{\tiny $\alpha_3$}}
\put(12.8, 1){\makebox(0,0)[c]{\tiny $\alpha_4$}}
\put(15.15,1){\makebox(0,0)[c]{\tiny $\alpha_5$}}
\end{picture}\vskip 8mm
\end{center}

\begin{center}
\hskip -1cm \setlength{\unitlength}{0.16in} \medskip
\begin{picture}(16,5.8)
\put(-5,1.8){$\e=(0, 1, 1, 0, 0, 1)$}
\put(6,0){\makebox(0,0)[c]{$\bigotimes$}}
\put(6,4){\makebox(0,0)[c]{$\bigotimes$}}
\put(8,2){\makebox(0,0)[c]{$\bigcirc$}}
\put(10.4,2){\makebox(0,0)[c]{$\bigotimes$}}
\put(12.8,2){\makebox(0,0)[c]{$\bigcirc$}}
\put(15.15,2){\makebox(0,0)[c]{$\bigotimes$}}
\put(6.35,0.3){\line(1,1){1.35}} 
\put(6.35,3.7){\line(1,-1){1.35}}
\put(8.4,2){\line(1,0){1.55}} 
\put(10.82,2){\line(1,0){1.55}}
\put(6, 0.45){\line(0,1){3.1}}
\put(13.2,2){\line(1,0){1.55}} 
\put(6,5){\makebox(0,0)[c]{\tiny $\alpha_0=-\delta_1-\delta_2$}}
\put(6,-1){\makebox(0,0)[c]{\tiny $\alpha_1$}}
\put(8.2,1){\makebox(0,0)[c]{\tiny $\alpha_2$}}
\put(10.4,1){\makebox(0,0)[c]{\tiny $\alpha_3$}}
\put(12.8, 1){\makebox(0,0)[c]{\tiny $\alpha_4$}}
\put(15.15,1){\makebox(0,0)[c]{\tiny $\alpha_5$}}
\end{picture}\vskip 8mm
\end{center}

\begin{center}
\hskip -1cm \setlength{\unitlength}{0.16in} \medskip
\begin{picture}(17,5.8)
\put(-5,1.8){$\e=(1, 1, 0, 0, 0, 1)$}
\put(6,0){\makebox(0,0)[c]{$\bigcirc$}}
\put(6,4){\makebox(0,0)[c]{$\bigcirc$}}
\put(8,2){\makebox(0,0)[c]{$\bigotimes$}}
\put(10.4,2){\makebox(0,0)[c]{$\bigcirc$}}
\put(12.8,2){\makebox(0,0)[c]{$\bigcirc$}}
\put(15.15,2){\makebox(0,0)[c]{$\bigotimes$}}
\put(6.35,0.3){\line(1,1){1.35}} 
\put(6.35,3.7){\line(1,-1){1.35}}
\put(8.4,2){\line(1,0){1.55}} 
\put(10.82,2){\line(1,0){1.55}}
\put(13.2,2){\line(1,0){1.55}} 
\put(6,5){\makebox(0,0)[c]{\tiny $\alpha_0=-\delta_1-\delta_2$}}
\put(6,-1){\makebox(0,0)[c]{\tiny $\alpha_1$}}
\put(8.2,1){\makebox(0,0)[c]{\tiny $\alpha_2$}}
\put(10.4,1){\makebox(0,0)[c]{\tiny $\alpha_3$}}
\put(12.8, 1){\makebox(0,0)[c]{\tiny $\alpha_4$}}
\put(15.15,1){\makebox(0,0)[c]{\tiny $\alpha_5$}}
\end{picture}\vskip 8mm
\end{center}

(2) For $\mf{osp}_C(\e)$,

\begin{center}
\hskip -1cm \setlength{\unitlength}{0.16in}
\begin{picture}(17,4)
\put(-5,1.8){$\e=(0, 0, 0, 1, 1, 1)$}
\put(5.6,2){\makebox(0,0)[c]{$\bigcirc$}}
\put(8,2){\makebox(0,0)[c]{$\bigcirc$}}
\put(10.3,2){\makebox(0,0)[c]{$\bigcirc$}}
\put(12.65,2){\makebox(0,0)[c]{$\bigotimes$}}
\put(15,2){\makebox(0,0)[c]{$\bigcirc$}}
\put(17.25,2){\makebox(0,0)[c]{$\bigcirc$}}
\put(6.8,2){\makebox(0,0)[c]{$\Longrightarrow$}}
\put(8.35,2){\line(1,0){1.5}}
\put(10.7,2){\line(1,0){1.5}}
\put(13.1,2){\line(1,0){1.5}}
\put(15.4,2){\line(1,0){1.45}}
%\put(24.3,2){\line(1,0){1.28}}

\put(5.4,0.8){\makebox(0,0)[c]{\tiny $\alpha_{0}=-2\delta_1$}}
\put(8.1,0.8){\makebox(0,0)[c]{\tiny $\alpha_{1}$}}
\put(10.33,0.8){\makebox(0,0)[c]{\tiny $\alpha_{2}$}}
\put(12.65, 0.8){\makebox(0,0)[c]{\tiny $\alpha_{3}$}}
\put(15, 0.8){\makebox(0,0)[c]{\tiny $\alpha_{4}$}}
\put(17.25, 0.8){\makebox(0,0)[c]{\tiny $\alpha_{5}$}}
\end{picture}
\end{center}

\begin{center}
\hskip -1cm \setlength{\unitlength}{0.16in}
\begin{picture}(17,4)
\put(-5,1.8){$\e=(1, 0, 0, 0, 0, 1)$}
\put(5.6,2){\makebox(0,0)[c]{$\bigcirc$}}
\put(8,2){\makebox(0,0)[c]{$\bigotimes$}}
\put(10.3,2){\makebox(0,0)[c]{$\bigcirc$}}
\put(12.65,2){\makebox(0,0)[c]{$\bigcirc$}}
\put(15,2){\makebox(0,0)[c]{$\bigcirc$}}
\put(17.25,2){\makebox(0,0)[c]{$\bigotimes$}}
\put(6.8,2){\makebox(0,0)[c]{$\Longrightarrow$}}
\put(8.35,2){\line(1,0){1.5}}
\put(10.7,2){\line(1,0){1.5}}
\put(13.1,2){\line(1,0){1.5}}
\put(15.4,2){\line(1,0){1.45}}
%\put(24.3,2){\line(1,0){1.28}}

\put(5.4,0.8){\makebox(0,0)[c]{\tiny $\alpha_{0}=-2\delta_1$}}
\put(8.1,0.8){\makebox(0,0)[c]{\tiny $\alpha_{1}$}}
\put(10.33,0.8){\makebox(0,0)[c]{\tiny $\alpha_{2}$}}
\put(12.65, 0.8){\makebox(0,0)[c]{\tiny $\alpha_{3}$}}
\put(15, 0.8){\makebox(0,0)[c]{\tiny $\alpha_{4}$}}
\put(17.25, 0.8){\makebox(0,0)[c]{\tiny $\alpha_{5}$}}
\end{picture}
\end{center}

} 
\end{ex}

Let $V$ be a $\U_X(\e)$-module with weight space decomposition $V=\bigoplus_{\la\in {P}}V_\la$. We briefly recall how to associate $V$ with an $\mf{osp}_X(\e)$-module (see \cite[Section 2.3]{KLO24} for more details).
For $v\in V_\la$ with $\la=\sum_{i\in \mathbb{I}} \la_i \de_i$, define
$\sigma_j v =(-1)^{\e_j \la_j}v$. Then we may naturally regard $V$ to a $\ms{Y}$-module and then we denote by $V^\tau$ the $\ms{X}$-module obtained by $\tau$ in Proposition \ref{prop:iso GQ and YQ}.
Let ${\bf A}$ be the subring of $\Bbbk$ consisting of $f(q^{\frac{1}{2}})\in \Bbbk$, which are regular at $q^{\frac{1}{2}}=1$.
Suppose that $V^\tau$ has an ${\bf A}$-lattice $V^\tau_{\bf A}$ invariant under $E_i$, $F_i$, $K_i^{\pm 1}$ for $i\in I$. For example, if $V$ is a highest weight module with a highest weight vector $v$, then we may take $V_{\bf A}$ to be the ${\bf A}$-span of $F_{i_1}\dots F_{i_l}v$ for $l\ge 0$ and $i_k\in I$ ($k=1,\dots,l$).
Let $\ov{V^\tau}=V^\tau_{{\bf A}}\otimes_{{\bf A}}\mathbb{\C}$, where $\C$ is an ${\bf A}$-module such that $f(q^\hf)\cdot c = f(1)c$ for $f(q^\hf)\in {\bf A}$ and $c\in \C$.
Then $\ov{V^\tau}$ is invariant under the $\mathbb{C}$-linear endomorphisms ${\rm E}_i$, ${\rm F}_i$ and ${\rm H}_i$ induced from $E_i$, $F_i$ and $[E_i,F_i]$, respectively for $i\in I$, which satisfy the defining relations for the enveloping algebra of $\mf{osp}_X(\e)$ (cf.~\cite[Section 10.5]{Ya94}). Hence $\ov{V^\tau}$ becomes an $\mf{osp}_X(\e)$-module, which we call a classical limit of $V$.

\section{$q$-deformed supersymmetric space and $q$-oscillator representations}\label{sec:Fock space}
Let $\Z_{+}^n(\e)=\{\,{\bf m}=(m_1, \ldots, m_n) \mid m_i \in \Z_{\geq 0} \text{ if } \e_i=0 \text{ and } m_i\in\{0, 1\} \text{ if } \e_i=1 \,\}$.
Let
\begin{equation}\label{eq:W}
\W=\W_{\e}=\bigoplus_{\mathbf{m}\in \Z_+^n(\e)}\Bbbk \lvert \mathbf{m} \rangle 
\end{equation}
be the $\Bbbk$-space with basis $\{\,|{\bf m}\rangle\,|\,{\bf m}\in \Z_{+}^n(\e)\,\}$. Let $\{\,{\bf e}_i\,|\,i=1,\dots,n\,\}$ be the standard basis of $\Z^n$. 
Then $\W$ has a $\U_A(\e)$-module structure given by
\begin{align}\label{eq:type A action}
e_i \lvert {\bf m} \rangle &= [m_{i+1}] \lvert {\bf m}+{\bf e}_i-{\bf e}_{i+1} \rangle, \quad
f_i \lvert {\bf m} \rangle = [m_{i}] \lvert {\bf m}-{\bf e}_i+{\bf e}_{i+1} \rangle, \quad
{k_{\de_j}} \lvert {\bf m} \rangle = q_j^{m_j} \lvert {\bf m} \rangle,
\end{align}
for $i\in I\setminus\{0\}$ and $j\in \I$. Here and also in Section \ref{subsubsec:osc-DC}, we assume that the righthand side is zero if it does not belong to $\Z_{+}^n(\e)$.

\subsection{$q$-oscillator representation $\ms{W}$ of $\U_X(\e)$}\label{subsubsec:osc-DC}
Let us consider a $\U_X(\e)$-module structure on $\W$, and also consider a $\U_X(\e)$-module structure on $\W^{\ot 2}$ not induced from $\Delta$. 
These $\U_X(\e)$-modules provide orthosymplectic analogues of spin representations of type $D_n$ and fundamental representations of type $C_n$. 
\subsubsection{$q$-oscillator $\U_D(\e)$-modules}\label{subsubsec:osc-D} 
\begin{prop}{\rm  (\cite[Proposition 6.4]{Ma})}\label{prop:osc-D}
When $\e_1=1$, the $\Bbbk$-space $\W$ has a $\U_D(\e)$-module structure given by
\begin{align*}
e_0 \lvert {\bf m} \rangle &= [m_2] \lvert {\bf m}-{\bf e}_1-{\bf e}_2 \rangle, \\
f_0 \lvert {\bf m} \rangle &= \lvert {\bf m} + {\bf e}_1 + {\bf e}_2 \rangle, \\
k_0 \lvert {\bf m} \rangle &= q_1^{1-m_1}q_2^{-m_2} \lvert {\bf m} \rangle.
\end{align*}
\end{prop}
Note that $\W$ is infinite-dimensional unless $\e=(1^n)$, and $\W=\W_+\oplus \W_-$, where $\W_+$ and $\W_-$ are irreducible highest weight modules generated by $|(0^n)\rangle$ and $|{\bf e}_1\rangle$ of weight $\frac{1}{2}\La$ and $\frac{1}{2}\La+\de_1$, respectively, which are the spin representations of $\U_{\td{q}}(D_n)$ when $\e=(1^n)$.

\begin{prop}{\rm (\cite[Proposition 2.8]{KLO24})}\label{prop:osc-D-2} 
When $\e_1=0$, the $\Bbbk$-space $\W^{\ot 2}$ has a $\U_D(\e)$-module structure given by
{\allowdisplaybreaks \small
\begin{align*}
&e_0 (\lvert {\bf m} \rangle \ot \lvert {\bf m}' \rangle) \\
&=  [m_2][m_1'] \lvert {\bf m}-{\bf e}_2 \rangle \ot \lvert {\bf m}'-{\bf e}_1 \rangle -q_1^{-m_1'}q_2^{-m_2}q^{-1} [m_1][m_2'] \lvert {\bf m} -{\bf e}_1 \rangle \ot \lvert {\bf m}'-{\bf e}_2 \rangle, \\
&f_0 (\lvert {\bf m} \rangle \ot \lvert {\bf m}' \rangle)   
= -q_1^{m_1}q_2^{m_2'}q \lvert {\bf m}+{\bf e}_2 \rangle \ot \lvert {\bf m}'+{\bf e}_1 \rangle + \lvert {\bf m}+{\bf e}_1 \rangle \ot \lvert {\bf m}'+{\bf e}_2 \rangle, \\
&k_0 (\lvert {\bf m} \rangle \ot \lvert {\bf m}' \rangle) = q_1^{-m_1-m_1'}q_2^{-m_2-m_2'}q^{-2} \lvert {\bf m} \rangle \ot \lvert {\bf m}' \rangle, \\
&e_i (\lvert {\bf m} \rangle \ot \lvert {\bf m}' \rangle)
=[m_{i+1}] \lvert {\bf m}+{\bf e}_i -{\bf e}_{i+1} \rangle \ot \lvert {\bf m}' \rangle + q_i^{m_i}q_{i+1}^{-m_{i+1}}[m_{i+1}'] \lvert {\bf m} \rangle \ot \lvert {\bf m}' +{\bf e}_i-{\bf e}_{i+1} \rangle, \\
&f_i (\lvert {\bf m} \rangle \ot  \lvert {\bf m}' \rangle)
= q_i^{-m_i'}q_{i+1}^{m_{i+1}'}[m_i] \lvert {\bf m}-{\bf e}_i+{\bf e}_{i+1} \rangle \ot \lvert {\bf m}' \rangle + [m_i'] \lvert {\bf m} \rangle \ot \lvert {\bf m}'-{\bf e}_i +{\bf e}_{i+1} \rangle, \\
& k_i (\lvert {\bf m} \rangle \ot \lvert {\bf m}' \rangle) = q_i^{m_i+m_i'}q_{i+1}^{-m_{i+1}-m_{i+1}'} \lvert {\bf m} \rangle \ot \lvert {\bf m}' \rangle,
\end{align*}}
for $i\in I\setminus\{0\}$. We denote this $\U_D(\e)$-module by $\W^2$.
\end{prop}
\medskip

Note that $\W^{2}$ is infinite-dimensional and semisimple (see also Proposition \ref{prop:irr comp in Fock space}). Its explicit decomposition together with highest weight vectors are given in \cite[Lemma 7.1]{KLO24}.

\subsubsection{$q$-oscillator $\U_C(\e)$-modules}\label{subsubsec:osc-C}
Let us give an analogue of Propositions \ref{prop:osc-D} and \ref{prop:osc-D-2} for $\U_C(\e)$.
\begin{prop}{\rm (\cite[Propositiom 6.3]{Ma})}\label{prop:osc-C}
When $\e_1=0$, the $\Bbbk$-space $\W$ has a $\U_C(\e)$-module structure given by
{\allowdisplaybreaks 
\begin{align*}
e_0 \lvert {\bf m} \rangle &= \dfrac{[m_1][m_1-1]}{[2]^2} \lvert {\bf m}-2{\bf e}_1 \rangle, \\
f_0 \lvert {\bf m} \rangle &= - \lvert {\bf m}+2 {\bf e}_1 \rangle, \\
k_0 \lvert {\bf m} \rangle &= q_1^{-2m_1-1} \lvert {\bf m} \rangle.
\end{align*}}
\end{prop}
Note that $\W$ is infinite-dimensional and $\W=\W_+\oplus \W_-$, where $\W_+$ and $\W_-$ are irreducible highest weight modules generated by $|(0^n)\rangle$ and $|{\bf e}_1\rangle$ of weight {$\frac{1}{2}\La$} and {${\frac{1}{2}\La}+\de_1$, respectively. It can be viewed as a super analogue of $q$-oscillator $U_q(C_n)$-module in \cite{Ha}. 

Finally, we introduce a $\U_C(\e)$-module structure on $\W^{\ot 2}$, which is a super analogue of $2$-spin realizations of a $U_{\td{q}}(C_n)$-module isomorphic to a direct sum of fundamental representations \cite{K14}.

\begin{prop}\label{prop:osc-C-2}
When $\e_1=1$, the $\Bbbk$-space $\W^{\ot 2}$ has a $\U_C(\e)$-module structure given by
{\allowdisplaybreaks \small
\begin{align*}
e_0 (\lvert {\bf m} \rangle \ot \lvert {\bf m}' \rangle) &= [m_1][m_1'] \lvert {\bf m}-{\bf e}_1 \rangle \ot \lvert {\bf m}'-{\bf e}_1 \rangle, \\
f_0 (\lvert {\bf m} \rangle \ot \lvert {\bf m}' \rangle) &= - \lvert {\bf m}+{\bf e}_1 \rangle \ot \lvert {\bf m}'+{\bf e}_1 \rangle, \\
k_0 (\lvert {\bf m} \rangle \ot \lvert {\bf m}' \rangle) &= q^{-2}q_1^{-2m_1-2m_1'} \lvert {\bf m} \rangle \ot \lvert {\bf m}' \rangle, \\
e_i (\lvert {\bf m} \rangle \ot \lvert {\bf m}' \rangle) &= [m_{i+1}]\lvert {\bf m}+{\bf e}_i -{\bf e}_{i+1} \rangle \ot \lvert {\bf m}' \rangle + q_i^{m_i}q_{i+1}^{-m_{i+1}} [m'_{i+1}]\lvert {\bf m}\rangle \ot \lvert {\bf m}'+{\bf e}_i -{\bf e}_{i+1} \rangle, \\
f_i (\lvert {\bf m} \rangle \ot \lvert {\bf m}' \rangle) &= q_i^{-m_i'}q_{i+1}^{m_{i+1}'}[m_i] \lvert {\bf m}-{\bf e}_i+{\bf e}_{i+1} \rangle \ot \lvert {\bf m}' \rangle + [m'_i]\lvert {\bf m} \rangle \ot \lvert {\bf m}' -{\bf e}_i +{\bf e}_{i+1} \rangle, \\
k_i (\lvert {\bf m} \rangle \ot \lvert {\bf m}' \rangle) &= q_i^{m_i+m_i'}q_{i+1}^{-m_{i+1}-m_{i+1}'} \lvert {\bf m} \rangle \ot \lvert {\bf m}' \rangle,
\end{align*}}
\! for $i\in I\setminus\{0\}$. We also denote this $\U_C(\e)$-module by $\W^2$.
\end{prop}
\pf It can be checked in a straightforward manner by using the defining relations of $\U_C(\e)$.
We show only the last two relations in \eqref{eq:rel for GQ-C-2}, which are the most involved. We use the notation \eqref{eq:monomial in U}.

\noindent {\em Case 1}. Let $(\e_1, \e_2, \e_3)=(1, 0, 1)$. Consider $B_1$ in the proof of Proposition \ref{prop:iso GQ and YQ} with $q_2=q$.

%We have
%\begin{equation}\label{eq: B_1}
%    \begin{split}
%        B_1&:=[[[[e_2, e_1]_{q_2}, e_0]_{q_1^{2}},[e_2, e_1]_{q_2}], e_1] \\
%        &= -q \langle 210121 \rangle +q^2 \langle 120121 \rangle - \langle 012121 \rangle -\langle 212101 \rangle  -q^2 \langle 121201 \rangle \\
%        &\quad -q^{-1} \langle 210121 \rangle + \langle 120121 \rangle -\langle 121021 \rangle +q \langle 121012 \rangle  +q^{-2} \langle 102121 \rangle \\
%        &\quad + \langle 101212 \rangle +\langle 121210 \rangle -q^{-2} \langle 121021 \rangle  +q^{-1} \langle 121012 \rangle,
%    \end{split}
%\end{equation}\comments{$q_2 \rightarrow q$}
%where we use the relations $\langle 11\rangle =\langle 22\rangle =0$ since $1, 2\in I_{\rm odd}$. 

Let us write $B_1$ as an operator on $\lvert {\bf m}\rangle \ot \lvert {\bf m}' \rangle$. For example, we have
\begin{equation}\label{eq: example of B_1}
    \begin{split}
    \langle 210121 \rangle %&=\langle 21012 \rangle (e_1 \ot 1 + k_1 \ot e_1)\\
    %&=\langle 2101 \rangle (e_2e_1 \ot 1 + k_2e_1\ot e_2 + e_2k_1 \ot e_1 + k_2k_1 \ot e_2e_1)\\
    %&=\langle 210 \rangle (k_1e_2e_1 \ot e_1 + k_1k_2e_1 \ot e_1e_2 + e_1e_2k_1\ot e_1 + e_1k_2k_1 \ot e_2e_1) \\
    %&=\langle 21 \rangle (e_0k_1e_2e_1 \ot e_0e_1 + e_0k_1k_2e_1 \ot e_0e_1e_2 \\
    %&\hspace{1cm} + e_0e_1e_2k_1\ot e_0e_1 + e_0e_1k_2k_1 \ot e_0e_2e_1) \\
    %&=\langle2\rangle (e_1e_0k_1e_2e_1 \ot e_0e_1 + k_1e_0k_1e_2e_1 \ot e_1e_0e_1\\
    %&\qquad + e_1e_0k_1k_2e_1 \ot e_0e_1e_2 + k_1e_0k_1k_2e_1 \ot e_1e_0e_1e_2 \\
    %&\qquad + e_1e_0e_1e_2k_1\ot e_0e_1 + k_1e_0e_1e_2k_1\ot e_1e_0e_1\\
    %&\qquad + e_1e_0e_1k_2k_1 \ot e_0e_2e_1 + k_1e_0e_1k_2k_1 \ot e_1e_0e_2e_1)\\
    &=k_2e_1e_0'k_1e_2e_1 \ot e_2e_0'e_1 + k_2k_1e_0'k_1e_2e_1 \ot e_2e_1e_0'e_1\\
    &\quad + e_2e_1e_0'k_1k_2e_1 \ot e_0'e_1e_2 + e_2k_1e_0'k_1k_2e_1 \ot e_1e_0'e_1e_2 \\
    &\quad + k_2e_1e_0'e_1e_2k_1\ot e_2e_0'e_1 + k_2k_1e_0'e_1e_2k_1\ot e_2e_1e_0'e_1\\
    &\quad + e_2e_1e_0'e_1k_2k_1 \ot e_0'e_2e_1 + e_2k_1e_0'e_1k_2k_1 \ot e_1e_0'e_2e_1,
    \end{split}
\end{equation}
where $e_0' \lvert {\bf m} \rangle:=[m_1]\lvert {\bf m}-{\bf e}_1 \rangle$.
Here we use the fact that, for $x=x_{1}x_{2}\cdots x_{l}$ with $x_{t}\in \{e_0', e_1, e_2, k_1, k_2\}$  ($t=1,\dots,l$) and $\lvert {\bf n}\rangle \in \W$, we have
\begin{itemize}
    \item $x \lvert {\bf n}\rangle=0$ if $|\{t\, |\, x_t=e_1\}|\geq 2$ and $x_{s}\neq e_0'$ for $l_1 < s < l_2$ with $x_{{l_1}}=x_{{l_2}}=e_1$,
    \item $x \lvert {\bf n}\rangle=0$ if $|\{t \mid x_t=e_2\}|=2$.
\end{itemize}
The other terms of the form $\langle i_1\dots i_6 \rangle$ in $B_1$ can be written similarly.

Note that $B_1(\lvert {\bf m} \rangle \ot \lvert {\bf m}' \rangle$) is a linear combination of ${\bf v}_1$ and ${\bf v}_2$, where ${\bf v}_1=\lvert {\bf m}+{\bf e}_1-{\bf e}_2-{\bf e}_3 \rangle \ot \lvert {\bf m}'-{\bf e}_3 \rangle$, and  ${\bf v}_2=\lvert {\bf m}-{\bf e}_3 \rangle \ot \lvert {\bf m}'+{\bf e}_1 -{\bf e}_2-{\bf e}_3 \rangle$. The vector ${\bf v}_1$ is obtained by applying the following terms in $B_1$ to $\lvert {\bf m}\rangle \ot \lvert {\bf m}' \rangle$:
\begin{equation*}
    \begin{split}
    B_{1,1}&:=[2]\{(qe_1e_2-e_2e_1)e_0'e_1-q^{-1}e_1e_0'(qe_1e_2-e_2e_1)\}k_2k_1 \ot e_2e_0'e_1 \\
    &\quad +q^{-1}[2]\{(qe_1e_2-e_2e_1)e_0'e_1-q^{-1}e_1e_0'(qe_1e_2-e_2e_1)\}k_2k_1 \ot e_1e_2e_0' \\
    &\quad -q^{-1}[2]\{(qe_1e_2-e_2e_1)e_0'e_1 -q^{-1}e_1e_0'(qe_1e_2-e_2e_1)\}k_2k_1 \ot e_0'e_1e_2 \\
    &\quad -[2]\{(qe_1e_2-e_2e_1)e_0'e_1-q^{-1}e_1e_0'(qe_1e_2-e_2e_1)\}k_2k_1 \ot e_2e_1e_0'.
    \end{split}
\end{equation*}
We can check directly that
\begin{gather*}
\{(qe_1e_2-e_2e_1)e_0'e_1-q^{-1}e_1e_0'(qe_1e_2-e_2e_1)\} \lvert {\bf m} \rangle=0,
\end{gather*}
which implies that $B_{1,1} (\lvert {\bf m}\rangle \ot \lvert {\bf m}' \rangle)=0$. The vector ${\bf v}_2$ is obtained by applying the following terms to $\lvert {\bf m}\rangle \ot \lvert {\bf m}' \rangle$:
\begin{equation*}
    \begin{split}
    B_{1,2}&:=\{-(q^{-3}+2q^{-1}+q)e_2e_0'e_1k_1^2k_2 +(q^2+2+q^{-2})e_0'e_1e_2k_1^2k_2 \\
    &\quad +(q^3+2q+q^{-1})e_2e_1e_0'k_1^2k_2-(q^4+2q^2+1)e_1e_2e_0'k_1^2k_2 \}\ot e_1e_0'e_2e_1 \\
    &\quad +\{(1+q^{-2})e_0'e_2e_1k_1^2k_2-(q+q^{-1})e_0'e_1e_2k_1^2k_2 \\
    &\quad -(1+q^2)e_2e_1e_0'k_1^2k_2+(q+q^3)e_1e_2e_0'k_1^2k_2\}\ot e_2e_1e_0'e_1 \\
    &\quad +\{(1+q^{-2})e_2e_0'e_1k_1^2e_2 -(q+q^{-1})e_0'e_1e_2k_1^2k_2 \\
    &\quad -(1+q^2)e_2e_1e_0'k_1^2k_2+(q+q^3)e_1e_2e_0'k_1^2k_2\}\ot e_1e_0'e_1e_2 \\
    &=-(q+q^{-1})E\ot e_1e_0'e_2e_1 +E\ot e_2e_1e_0'e_1 +E\ot e_1e_0'e_1e_2,
    \end{split}
\end{equation*}
where
\begin{equation*}
    E=[2](q^{-1}e_2e_0'e_1k_1^2k_2-e_0'e_1e_2k_1^2k_2-qe_2e_1e_0'k_1^2k_2+q^2e_1e_2e_0'k_1^2k_2).
\end{equation*}
We obtain $E\lvert {\bf m}\rangle=q^{-2m_1+4m_2-m_3}[2]\de_{m_3,1}(\de_{m_1,1}q^2-\de_{m_1,0})\lvert {\bf m}-{\bf e}_3\rangle$, and hence
\begin{align*}
    B_{1,2}(\lvert {\bf m}\rangle \ot \lvert {\bf m}'\rangle) &= \{-(q+q^{-1})E\ot e_1e_0'e_2e_1 + E \ot e_2e_1e_0'e_1 + E \ot e_1e_0'e_1e_2\} (\lvert{\bf m}\rangle \ot \lvert {\bf m}'\rangle)\\
    &=q^{-2m_1+4m_2-m_3}[2][m_2'](\de_{m_1,1}q^2-\de_{m_1,0})\de_{m_1',0}\de_{m_3,1}\de_{m_3',1}\\
    &\quad (-[2][m_2']+[m_2'-1]+[m_2'+1])(\lvert {\bf m}-{\bf e}_3\rangle \ot \lvert {\bf m}'+{\bf e}_1-{\bf e}_2\rangle )=0.
\end{align*}
Therefore, $B_1(\lvert {\bf m}\rangle \ot \lvert {\bf m}'\rangle)=B_{1,1}(\lvert {\bf m}\rangle \ot \lvert {\bf m}'\rangle)+B_{1,2}(\lvert {\bf m}\rangle \ot \lvert {\bf m}'\rangle)=0$.\smallskip

\noindent {\em Case 2.} Let $(\e_1, \e_2, \e_3)=(1,0,0)$. Consider $B_2$ in the proof of Proposition \ref{prop:iso GQ and YQ} with $q_3=q$.
Similarly, we write $B_2$ as an operator on $\lvert {\bf m}\rangle \ot \lvert {\bf m}' \rangle$.
For example, we have
\begin{equation}\label{eq: example of B_2}
    \begin{split}
    \langle 3210121 \rangle &=e_3e_2e_1e_0'Ek_1 \ot e_0'e_1 + q^{-2}e_2e_1e_0'Ek_1k_3 \ot e_3e_0'e_1 \\
    &\quad +e_3e_1e_0'Ek_1k_2 \ot e_2e_0'e_1 +e_3e_2e_1e_0'e_1k_1k_2 \ot e_0'E' \\
    &\quad +q^{-1}e_1e_0'Ek_1k_2k_3\ot e_3e_2e_0'e_1 +q^{-1}e_2e_1e_0'e_1k_1k_2k_3 \ot e_3e_0'E' \\
    &\quad +e_3e_2e_0'Ek_1^2 \ot e_1e_0'e_1 +q^{-2}e_2e_0'Ek_1^2k_3 \ot e_3e_1e_0'e_1 \\
    &\quad +qe_3e_0'Ek_1^2k_2 \ot e_2e_1e_0'e_1 +qe_3e_2e_0'e_1k_1^2k_2\ot e_1e_0'E' \\
    &\quad +e_0'Ek_1^2k_2k_3\ot e_3e_2e_1e_0'e_1 +e_2e_0'e_1k_1^2k_2k_3\ot e_3e_1e_0'E' \\
    &\quad +q^{-2}e_3e_1e_0'e_1k_1k_2^2 \ot e_2e_0'E' +q^{-2}e_1e_0'e_1k_1k_2^2k_3 \ot e_3e_2e_0'E' \\
    &\quad +e_3e_0'e_1k_1^2k_2^2 \ot e_2e_1e_0'E' +e_0'e_1k_1^2k_2^2k_3 \ot e_3e_2e_1e_0'E',
    \end{split}
\end{equation}
where $E=e_1e_2-q^{-1}e_2e_1$ and $E'=e_2e_1-q^{-1}e_1e_2$.

Note that $B_2(\lvert {\bf m} \ot \lvert {\bf m}'\rangle)$ is a linear combination of $\lvert {\bf m}''\rangle \ot \lvert {\bf m}'''\rangle$, where $\lvert {\bf m}'''\rangle$ is one of the following:
\begin{equation}\label{eq: 12 cases of B_2}
    \begin{aligned}
    &  \lvert {\bf m}'-{\bf e}_2\rangle ,& &\lvert {\bf m}'-{\bf e}_3 \rangle ,& & \lvert {\bf m}'-{\bf e}_4 \rangle ,& \\
    & \lvert {\bf m}'+{\bf e}_1-2{\bf e}_2 \rangle ,& & \lvert {\bf m}'+{\bf e}_1-2{\bf e}_3 \rangle ,& & \lvert {\bf m}'+{\bf e}_2-2{\bf e}_3 \rangle ,& \\
    & \lvert {\bf m}'+{\bf e}_1-{\bf e}_2-{\bf e}_3 \rangle ,& &  \lvert {\bf m}'+{\bf e}_1-{\bf e}_2-{\bf e}_4 \rangle ,& &  \lvert {\bf m}'+{\bf e}_1-{\bf e}_3-{\bf e}_4 \rangle ,& \\
    &  \lvert {\bf m}'+{\bf e}_2-{\bf e}_3-{\bf e}_4 \rangle ,& &  \lvert {\bf m}'-{\bf e}_2+{\bf e}_3-{\bf e}_4\rangle ,& &  \lvert {\bf m}'+{\bf e}_1-2{\bf e}_2+{\bf e}_3-{\bf e}_4\rangle &.
    \end{aligned}
\end{equation}
Then, for example, $\lvert {\bf m}'' \rangle \ot \lvert {\bf m}'-{\bf e}_2\rangle$ is obtained by applying the following terms in $B_2$:
{\small \allowdisplaybreaks{
\begin{align*}
    &e_3e_2e_1e_0'Ek_1\ot e_0'e_1 -[2]e_2e_3e_1e_0'Ek_1 \ot e_0'e_1 -(q+q^{-3}) e_3e_1e_2e_0'Ek_1 \ot e_0'e_1 \\
    &+[3]e_1e_2e_3e_0'Ek_1 \ot e_0'e_1 +q^{-1}e_3Ee_2e_0'e_1k_1 \ot e_1e_0' -[3]Ee_3e_2e_0'e_1k_1 \ot e_1e_0' \\
    &+q^{-2}e_3e_1e_0'(e_1e_2e_2-q^{-2}e_2e_2e_1)k_1 \ot e_0'e_1 -q^{-1}e_1e_0'(e_1e_2e_3e_2-q^{-2}e_2e_3e_2e_1)k_1 \ot e_0'e_1 \\
    &-qe_2e_3Ee_0'e_1k_1 \ot e_1e_0' +q^2e_2Ee_3e_0'e_1k_1 \ot e_1e_0' +e_2e_1e_0'Ee_3k_1 \ot e_0'e_1 -e_3Ee_0'e_1e_2k_1 \ot e_1e_0' \\
    &+[2]Ee_3e_0'e_1e_2k_1 \ot e_1e_0' -{q^{-2}}e_3e_1e_0'Ee_2k_1 \ot e_0'e_1 +q^{-1}e_1e_0'Ee_3e_2k_1 \ot e_0'e_1 \\
    &+{q^2}(e_1e_2e_3e_2-q^{-2}e_2e_3e_2e_1)e_0'e_1k_1 \ot e_1e_0' -{q^3}(e_1e_2e_2-q^{-2}e_2e_2e_1)e_3e_0'e_1k_1 \ot e_1e_0' \\
    &+(q^3+q^{-1})Ee_2e_3e_0'e_1k_1 \ot e_1e_0' -qe_1e_2e_0'Ee_3k_1 \ot e_0'e_1 -Ee_0'e_1e_2e_3k_1 \ot e_1e_0'
\end{align*}
which yields
\begin{align*}    
&q_1^{m_1}\de_{m_1,0}\de_{m_1',0} [m_4][m_3] \{[m_2]([m_3-1]-[2][m_3]+[m_3+1]) \\
    &{\quad +[m_2+1](-q[m_3-1]+(q^2+1)[m_3]-q[m_3+1])\}\lvert {\bf m}''\rangle \ot \lvert {\bf m}'-{\bf e}_2\rangle} \\
    &{\quad + q_1^{m_1}\de_{m_1,0}\de_{m_1',1} [m_4][m_3] \{[m_2](q^{-1}[m_3-1]-(1+q^{-2})[m_3]+q^{-1}[m_3+1])} \\
    &{\quad +[m_2+1](-[m_3-1]+[2][m_3]-[m_3+1])\}\lvert {\bf m}''\rangle \ot \lvert {\bf m}'-{\bf e}_2\rangle} \\
    &=0 \ot \lvert {\bf m}'-{\bf e}_2\rangle + 0 \ot \lvert {\bf m}'-{\bf e}_2\rangle=0.
\end{align*}}}
Thus the coefficient of $\lvert {\bf m}'' \rangle \ot \lvert {\bf m}'-{\bf e}_2\rangle$ is $0$.
Similarly, we can check the coefficient of $\lvert {\bf m}'' \rangle \ot \lvert {\bf m}'''\rangle$ is $0$ for the other $\lvert {\bf m}'''\rangle$ in \eqref{eq: 12 cases of B_2}. Therefore, $B_2(\lvert {\bf m}\rangle \ot \lvert {\bf m}' \rangle)=0$.
\qed}
\medskip

From now on, we denote by $\ms{W}=\ms{W}_\e$ either $\W$ or $\W^{2}$, the $\U_X(\e)$-module defined in Propositions \ref{prop:osc-D}--\ref{prop:osc-C-2}. Note that the weight of $|{\bf m}\rangle$ for ${\bf m}=(m_i)\in \Z^n_+(\e)$ is  
\begin{equation}\label{eq:wt of m}
{\hf\La +\sum_{i\in \I}m_i\de_i.}
\end{equation}

\subsection{Classical limit of $\ms{W}^{\ot \ell}$ and Howe duality}\label{subsec:classical Howe duality}

Let $\ms{W}_{\bf A}$ be the ${\bf A}$-span of $|{\bf m}\rangle$ or $|{\bf m}\rangle \ot |{\bf m}'\rangle$ for ${\bf m}, {\bf m}'\in \Z_{+}^n(\e)$. % except that when $\ms{W}=\W$ is a $\U_C(\e)$-module in Section \ref{subsubsec:osc-C} we assume that $\ms{W}_{\bf A}$ is the ${\bf A}$-span of $F_{i_1}\dots F_{i_l}|(0^n)\rangle$ for $l\ge 0$ and $i_k\in I$ ($k=1,\dots,l$). 
It is not difficult to see that $\ms{W}_{\bf A}$ admits a well-defined classical limit, which is an $\mf{osp}_X(\e)$-module, and hence so does $(\ms{W}_{\bf A})^{\ot \ell}$ ($\ell\ge 1$).
%\comments{$\ov{\left(\ms{W}_{\bf A}^{\ot \ell}\right)^\tau}=\left(\ov{\ms{W}_{\bf A}^\tau}\right)^{\ot \ell}$ is removed since it may not be true. See Remark \ref{rem:tau iso} (20250416)}
Let
\begin{equation}\label{eq:classical limit of W}
 W=\ov{\ms{W}_{\bf A}^\tau}.
\end{equation}
We may view it as a $\C$-space
\begin{equation}\label{eq:classical limit of W-2}
 W^{\ot \ell} =  
\begin{cases}
 S(\C^\ell \ot \C^{n_0|n_1}) = S(\C^{n_0|n_1})^{\ot \ell}& \text{if $\ms{W}=\W$},\\
 S(\C^{2\ell} \ot \C^{n_0|n_1}) = S(\C^{n_0|n_1})^{\ot 2\ell}& \text{if $\ms{W}=\W^{2}$},\\
\end{cases}
\end{equation}
where $S(V)$ denotes the supersymmetric algebra generated by a superspace, and $n_0$ and $n_1$ are the numbers of 0's and 1's in $\e$, respectively.

Let $G_\ell$ denote either a complex algebraic group $O_\ell$ ($\ell\ge 2$) or $Sp_{2\ell}$ ($\ell\ge 1$), and let   
\begin{equation}\label{eq:partitions for classical groups}
\begin{split}
&\mc{P}(O_\ell) =\{\,\la=(\la_1,\dots,\la_\ell)\in \Z_{\ge 0}^\ell\,|\, \la_1\ge \dots\ge \la_\ell,\ \la'_1+\la'_2\le \ell\,\},\\
&\mc{P}(Sp_{2\ell}) =\{\,\la=(\la_1,\dots,\la_\ell)\in \Z_{\ge 0}^\ell\,|\, \la_1\ge \dots\ge \la_\ell \,\}, 
\end{split}
\end{equation}
which parametrize the finite-dimensional irreducible $G_\ell$-modules, say $V_{G_\ell}(\la)$. Here $\la'=(\la'_1,\la'_2,\dots)$ is the conjugate of $\la$ when we regard it as a partition.

For $G_\ell=Sp_{2\ell}$ and $\la\in \mc{P}(Sp_{2\ell})$, $V_{Sp_{2\ell}}(\la)$ is also an irreducible $\mf{sp}_{2\ell}$-module, which we denote by $V_{\mf{sp}_{2\ell}}(\la)$. 
For $G_\ell=O_{\ell}$ and $\la\in \mc{P}(O_{\ell})$, $V_{O_{\ell}}(\la)$ is an irreducible $\mf{so}_{\ell}$-module, say $V_{\mf{so}_{\ell}}(\la)$, when either $\ell$ is odd or $\ell$ is even with $\la_\ell=0$. Otherwise, $V_{O_{\ell}}(\la)$ is a direct sum of two irreducible $\mf{so}_{\ell}$-modules, which we denote by $V_{\mf{so}_{\ell}}(\la)$ and $V_{\mf{so}_{\ell}}(\la^{-})$ (see \cite{BT} and \cite[Chapter 5.3]{CW12} for more details).

Let $\ms{P}_\e$ be the set of partitions $\la = (\la_i)_{i\ge 1}$ such that  $\la_{n_0+1}\le n_1$, which parametrizes the irreducible polynomial representations of $\U_A(\e)$ (see \cite[Section 3.2]{KY}). %\comments{Should be corrected (2025.07.02)}
Put $\mc{P}(G_\ell)_\e=\mc{P}(G_\ell)\cap \ms{P}_\e$.
For $\lambda \in \mc{P}(G_\ell)_\e$, let
\begin{equation}\label{eq:highest weight La}
\Lambda_{\lambda, \e}= r\ell \Lambda + \sum_{i\in\I} m_i \de_i,
\end{equation}
where $r=\frac{1}{2}$ (resp. $r=1$) for $G_{{\ell}}=O_{{\ell}}$ (resp. $Sp_{2{\ell}}$), and $m_i=(-1)^{\e_i}(\omega_\la|\de_i)$. Here $\omega_\la$ denotes the highest weight of an irreducible polynomial representations of $\U_A(\e)$ corresponding to $\la\in \cP_\e$ (see \cite{CW12} or \cite[Section 3.2]{KY} for its explicit combinatorial description).
For $\lambda \in \mc{P}(G_\ell)_\e$, let $V^\la=V_\e^\lambda$ be the irreducible highest weight $\mf{osp}_X(\e)$-module with highest weight $\Lambda_{\lambda,\e}$.

\begin{thm}\label{thm:classical Howe duality}{\rm (\cite[Theorems 5.1 and 5.2]{CZ})}
There exists an $(\mf{osp}_X(\e),G_\ell)$-action on $W^{\ot \ell}$ which gives the following multiplicity-free decomposition as an $(\mf{osp}_X(\e),G_\ell)$-module:
\begin{equation}\label{eq:osp Howe duality}
W^{\ot \ell} = 
\bigoplus_{\la\in \mc{P}(G_\ell)_{\e}} V^\la \ot V_{G_\ell}(\la),
\end{equation}
where $G_\ell=O_\ell$ when $W^{\ot \ell}=S(\C^{n_0|n_1})^{\ot \ell}$ and $G_\ell=Sp_{2\ell}$ when $W^{\ot \ell}=S(\C^{n_0|n_1})^{\ot 2\ell}$.
\end{thm}

\subsection{Semisimplicity and irreducible components of $\ms{W}^{\ot\ell}$}
From now on, we assume the following convention:
\begin{equation}\label{eq:G convention}
 G=
\begin{cases}
 O & \text{if $\ms{W}=\W$},\\
 Sp & \text{if $\ms{W}=\W^{2}$}. 
\end{cases}
\end{equation}
For $\lambda \in \mc{P}(G_\ell)_\e$, let $\V^\la=\V_\e^\lambda$ be the irreducible highest weight $\U_X(\e)$-module with highest weight $\Lambda_{\lambda, \e}$.

\begin{prop}\label{prop:irr comp in Fock space}
$\ms{W}^{\ot \ell}$ is a semisimple $\U_X(\e)$-module, and an irreducible $\U_X(\e)$-module in $\ms{W}^{\ot \ell}$ is isomorphic to $\V^\la$ for some $\lambda \in \mc{P}(G_\ell)_\e$. Moreover,
\begin{equation*}
 \ms{W}^{\ot \ell} = 
\bigoplus_{\la\in \mc{P}(G_\ell)_{\e}} (\V^\la)^{\oplus \dim V_{G_\ell}(\la)}.
\end{equation*}
\end{prop}
\pf %\comments{Revised. (20250416)}
Let $V$ be a submodule of $\ms{W}^{\ot\ell}$.
We see from \eqref{eq:wt of m} and the definition of $\U_X(\e)$-actions on $\ms{W}$ that
the set of weights of $\ms{W}^{\ot\ell}$ is finitely dominated , that is, it is contained in $\bigcup_{k=1}^s(\Lambda_k - \sum_{i\in I}\Z_{\ge 0}\alpha_i)$ for some $s\ge 1$ and $\Lambda_k\in P$ ($k=1,\dots,s$), since the set of weights of $\ms{W}$ is contained in a finite (at most two) union of $r\Lambda +c\de_1 - \sum_{i\in I}\Z_{\ge 0}\alpha_i$ with $r$ in \eqref{eq:highest weight La} and $c=0,1$. This implies that $V$ has a weight vector $v$ of maximal weight $\omega$, that is, $e_i v=0$ for all $i\in I$. 
Let $V_0$ be the submodule of $V$ generated by $v$. 
Then $\ov{V_0^\tau}$ is also a highest weight module in $W^{\ot \ell}$.

Since $\ov{\left(\ms{W}_{\bf A}^{\ot \ell}\right)^\tau}$ is semisimple (see Remark \ref{rem:semisimplicity of osc module}), $\ov{V_0^\tau}$ is isomorphic to $V^\la$ for some $\la\in \mc{P}(G_\ell)_\e$, which implies that $V_0$ is irreducible with highest weight $\omega=\La_{\la,\e}$.
Hence every highest weight submodule in $\ms{W}^{\ot \ell}$ is irreducible and isomorphic to $\mc{V}^\la$ for some $\la\in \mc{P}(G_\ell)_\e$.  

Let us prove that $\W^{\ot \ell}$ is semisimple.
Let $\ms{W}^\vee = \bigoplus_{\mu} {\rm Hom}_\Bbbk(\ms{W}_\mu,\Bbbk)$ be the restricted dual of $\ms{W}$, which is a $\U_X(\e)$-module with respect to the antipode $S$ in \eqref{eq:comult}. Then by the same argument as in $\ms{W}^{\ot \ell}$, we see that every lowest weight submodule in $(\ms{W}^\vee)^{\ot \ell}$ is irreducible with lowest weight $-\La_{\la,\e}$ for some $\la\in \mc{P}(G_\ell)_\e$.
Now applying the same arguments as in \cite[Lemma 3.5.3]{HK}, we may conclude that $\ms{W}^{\ot \ell}$ is semisimple.%\comments{This paragraph has been added. (20250506)}

Note that $\ov{\left(\ms{W}_{\bf A}^{\ot \ell}\right)^\tau}=\left(\ov{\ms{W}_{\bf A}^\tau}\right)^{\ot \ell}=W^{\ot \ell}$ as a $\C$-space. It is not difficult to see from Section \ref{subsubsec:osc-DC} that the character of $\ov{\left(\ms{W}_{\bf A}^{\ot \ell}\right)^\tau}$ as an $\mf{osp}_X(\e)$-module is equal to that of $W^{\ot \ell}$ in Theorem \ref{thm:classical Howe duality} (see \cite[Section 6]{CZ}). Hence, the multiplicity of each $\mc{V}^\la$ follows from Theorem \ref{thm:classical Howe duality} since the characters of $V^\la$ are linearly independent.
\qed
\medskip

In the proof of Proposition \ref{prop:irr comp in Fock space}, we also have proved the following.

\begin{cor}\label{cor:classical limit of Vlambda}
We have $\ov{(\mc{V}^\la)^\tau}\cong V^\la$ for $\la\in \mc{P}(G_\ell)_\e$. 
\end{cor}

\begin{rem}\label{rem:semisimplicity of osc module}
{\rm From a viewpoint of super duality \cite{CLW}, the $\mf{osp}_X(\e)$-module $W^{\ot \ell}$ belongs to a semisimple category, which corresponds to a semisimple category having integrable highest weight modules of types $C_\infty$ and $D_\infty$ as irreducible objects, when it is extended to modules over orthosymplectic Lie superalgebra of infinite-rank. %More details on this semisimple category of $\mf{osp}_X(\e)$-modules can be found in \cite[Section 3]{K15} when $\e$ is standard, that is, $\e=(0^{n_0},1^{n_1})$, while the proof for general $\e$ is similar.
} 
\end{rem}

\subsection{Polarization on $\ms{W}^{\ot\ell}$}

Define a non-degenerate symmetric bilinear form on {$\W$} by
\begin{equation}\label{eq:polarization}
(\lvert {\bf m}\rangle, \lvert {\bf m}'\rangle )= \delta_{{\bf m}, {\bf m}'}q^{\sum_{i\in\I}\frac{m_i(m_i-1)}{2}}\prod_{i\in\I} [m_i]!,
\end{equation}
for ${\bf m}=(m_i)$ and ${\bf m}'=(m'_i)$.
On {$\W^{\ot 2}$}, we define a bilinear form by
$(\lvert {\bf m}\rangle\ot \lvert {\bf m}'\rangle, \lvert {\bf m}''\rangle\ot \lvert {\bf m}'''\rangle)=(\lvert {\bf m}\rangle, \lvert {\bf m}'' \rangle)(\lvert {\bf m}' \rangle, \lvert {\bf m}''' \rangle)$. Similarly, we define a bilinear form on {$\W^{\ot l}$} for $l\ge 1$.

Let $\eta_D$ be an anti-involution on $\U_D(\e)$ given by
\begin{equation}\label{eq:eta_D}
\begin{split}
\eta_D(k_\mu)&=k_\mu, \\
\eta_D(e_i)&=\begin{cases}
q_i^{-1}f_ik_i &\quad \text{if $i\in I\setminus\{0\}$}, \\
(-1)^{\e_1+1}q_1^{(-1)^{\e_1}}f_0k_0 &\quad \text{if $i=0$},
\end{cases}, \\
\eta_D(f_i)&=\begin{cases}
q_ik_i^{-1}e_i &\quad \text{if $i\in I\setminus\{0\}$}, \\
(-1)^{\e_1+1}q_1^{(-1)^{\e_1+1}}k_0^{-1}e_0 &\quad \text{if $i=0$}.
\end{cases}
\end{split}
\end{equation}
Also, let $\eta_C$ and $\eta'_C$ be anti-involutions on $\U_C(\e)$ given by
\begin{equation}\label{eq:eta_C}
\begin{split}
\eta_C(k_\mu)&=k_\mu, \\
\eta_C(e_i)&=
\begin{cases}
q_i^{-1}f_ik_i &\quad \text{if $i\in I\setminus\{0\}$}, \\
-\frac{1}{[2]^2}f_0k_0 &\quad \text{if $i=0$}, 
\end{cases}, \\
\eta_C(f_i)&=\begin{cases}
q_ik_i^{-1}e_i &\quad \text{if $i\in I\setminus\{0\}$}, \\
-[2]^2k_0^{-1}e_0 &\quad \text{if $i=0$},
\end{cases}
\end{split}
\end{equation}
\begin{equation}\label{eq:eta'_C}
    \begin{split}
        \eta_C'(k_\mu)&=k_\mu, \\
        \eta_C'(e_i)&=\begin{cases}
            q_i^{-1}f_ik_i &\quad \text{if $i\in I\setminus\{0\}$}, \\
            -q^2f_0k_0 &\quad \text{if $i=0$},
        \end{cases}, \\
        \eta_C'(f_i)&=\begin{cases}
            q_ik_i^{-1}e_i &\quad \text{if $i\in I\setminus\{0\}$}, \\
            -q^{-2}k_0^{-1}e_0 &\quad \text{if $i=0$}.
        \end{cases}
    \end{split}
\end{equation}

Then we have the following lemmas, which can be checked in a straightforward manner. The proofs are quite similar to that of \cite[(2.9)]{BKK} for Lemma \ref{lem:compatibility of eta}, and those of \cite[Lemma 5.2]{KLO24} and \cite[Lemma 3.6]{KwO21} for Lemma \ref{lem:polarization on level 1}. So we leave the details to the reader.

\begin{lem}\label{lem:compatibility of eta}
We have $(\eta \ot \eta) \circ \Delta = \Delta \circ \eta$, where $\eta$ denotes one of \eqref{eq:eta_D}--\eqref{eq:eta'_C}.
\end{lem}

\begin{lem}\label{lem:polarization on level 1}
We have
\begin{itemize}
\item[(1)] $(xv, w)=(v, \eta_D(x)w)$ on {$\W$ and $\W^{2}$} for $x\in \U_D(\e)$, $v \in {\W^{2}}$, $w \in {\W^{2}}$,
\item[(2)] $(xv, w)=(v, \eta_C(x)w)$ on {$\W$}  for $x\in \U_C(\e)$ and $v, w \in {\W}$, 
\item[(3)] $(xv, w)=(v, \eta'_C(x)w)$ on {$\W^{2}$} for $x\in \U_C(\e)$ and $v, w \in {\W^{2}}$. 
\end{itemize}
\end{lem}

\begin{prop}
Let $\eta$ be one of \eqref{eq:eta_D}--\eqref{eq:eta'_C}, and $\ell\ge 1$.
We have 
\begin{equation*}
(xv, w)=(v, \eta(x)w),
\end{equation*}
for $x\in \U_X(\e)$ and $v, w \in {\ms{W}^{\ot \ell}}$, where $(\, ,\,)$ is a bilinear form on {$\ms{W}^{\ot \ell}$} induced from \eqref{eq:polarization}.
\end{prop}
\pf It follows from Lemmas \ref{lem:compatibility of eta} and \ref{lem:polarization on level 1}.
\qed

\begin{cor}
$\ms{W}^{\ot \ell}$ is a semisimple $\U_X(\e)$-module. 
\end{cor}
\pf It follows from the same arguments in the proof of \cite[Theorem 2.12]{BKK}.
\qed

\section{Quantum symmetric pairs and Howe duality} \label{sec:QSP}

\subsection{Quantum symmetric pairs $({\bf U}(\mf{sl}_\ell),{\bf U}^\imath(\mf{so}_\ell))$ and $({\bf U}(\mf{sl}_{2\ell}),{\bf U}^\imath(\mf{sp}_{2\ell}))$}\label{subsubsection:QSP}

For $r\ge 2$, let ${\bf U}(\mf{sl}_r)$ be the subalgebra of $\U_A((0^r))\cong U_q(\mf{gl}_r)$ generated by $e_i, f_i, k_i^{\pm 1}$ for $i=1,\dots,r-1$. 
To avoid confusion with the notations for $\U_X(\e)$ in Section \ref{sec:Qsuper}, we denote the generators of ${\bf U}(\mf{sl}_r)$ by $\mathsf{e}_i,\mathsf{f}_i,\mathsf{k}_i^{\pm 1}$  for $i\in \mathsf{I}:=\{\,1,\dots,r-1\,\}$.

Let $W$ be the Weyl group of $\mf{sl}_r$ generated by $s_i$ ($i\in \msf{I}$). 
For $i\in \msf{I}$, let $T_i$ be the automorphism on ${\bf U}(\mf{sl}_r)$ given by
\begin{gather*}
T_{i}(\mathsf{k}_j)=\mathsf{k}_{s_i(\alpha_j)},\quad
T_{i}(\mathsf{e}_i)=-\mathsf{f}_i\mathsf{k}_i, \quad 
T_{i} (\mathsf{f}_i)=-\mathsf{k}_i^{-1}\mathsf{e}_i, \\
T_{i}(\mathsf{e}_j)=\sum_{r+s=-(\alpha_i|\alpha_j)}(-1)^r{q}^{-r}\mathsf{e}_i^{(s)}\mathsf{e}_j\mathsf{e}_i^{(r)} \quad  
T_{i}(\mathsf{f}_j)=\sum_{r+s=-(\alpha_i|\alpha_j)}(-1)^r{q}^r\mathsf{f}_i^{(r)}\mathsf{f}_j\mathsf{f}_i^{(s)} \quad (j\neq i),
\end{gather*}%\comments{$q_i$ has been replaced by $q$}
for $j\in \msf{I}$, which is $T_{i, 1}''$ in \cite[37.1.3]{Lu93}.
For a reduced expression $w=s_{i_1} \cdots s_{i_k}\in W$, let $T_w=T_{i_1} \cdots T_{i_k}$.

Let us consider the quantum symmetric pairs associated to the symmetric pair $(\mf{g},\mf{k})=(\mf{sl}_\ell,\mf{so}_\ell)$ ($\ell\ge 2$), and $(\mf{sl}_{2\ell},\mf{sp}_{2\ell})$ ($\ell\ge 1$) which correspond to the Satake diagrams of type AI and AII, respectively. We follow the convention in \cite{BW18}\medskip

{\em Case 1.} $(\mf{g},\mf{k})=(\mf{sl}_\ell,\mf{so}_\ell)$.
Let ${\bf U}^\imath(\mf{so}_\ell)$ be the $\Bbbk$-subalgebra of ${\bf U}(\mf{sl}_\ell)$ generated by
\begin{equation}\label{eq:B_i for so}
\msf{B}_i=\msf{f}_i+\varsigma_i \msf{e}_i \msf{k}_i^{-1}+\kappa_i \msf{k}_i^{-1} \quad (i\in \msf{I}), 
\end{equation}
where $\varsigma_i\in \Bbbk$ is such that $\varsigma_i^2=q^{-2}$ for $i\in \msf{I}$, and $\kappa_i\in \Bbbk$ is such that $\kappa_i=0$ for $i\in \msf{I}$ unless $\msf{I}=\{1\}$ and $\kappa_1=\ov{\kappa_1}$ in this case. For later use, we put $\msf{I}_\circ=\msf{I}$ and $\msf{I}_\bullet=\emptyset$.

{\em Case 2.} $(\mf{g},\mf{k})=(\mf{sl}_{2\ell},\mf{sp}_{2\ell})$. 
Let ${\bf U}^\imath(\mf{sp}_{2\ell})$ be the $\Bbbk$-subalgebra of ${\bf U}(\mf{sl}_{2\ell})$ generated by 
\begin{gather}\label{eq:B_i for sp}
\msf{B}_i=\msf{f}_i+\varsigma_i T_{w_\bullet}(\msf{e}_i)\msf{k}_i^{-1} \ (i\in \msf{I}_\circ),\quad \msf{e}_j, \msf{f}_j, \msf{k}_j \ (j\in \msf{I}_\bullet),
\end{gather}
where $w_\bullet=s_1 s_3 \cdots s_{2\ell-1}$, $\msf{I}_\circ=\{\,2, 4, \ldots, 2\ell-2\,\}$, $\msf{I}_\bullet =\{1, 3, \ldots, 2\ell-1\}$, and $\varsigma_i\in \Bbbk$ is such that $\varsigma_i^2=q^{2}$ for $i\in \msf{I}_\circ$.
\medskip

Recall from \cite[Chapter 12.3]{GW} that the Lie algebra $\mf{k}$ in the above cases is given by $\mf{k} =\{\,X\,|\,X\in \mf{sl}_{r}, X^TJ_1+J_1X={\mathbb O}\,\}$, where ${\mathbb O}$ denotes the zero matrix, $J_1=I_\ell$, the identity matrix for $\mf{k}=\mf{so}_\ell$, and $J_1={\rm diag}(\Omega,\dots,\Omega)$, the $2\ell \times 2\ell$-block diagonal matrix with  
$\Omega=\begin{bmatrix} 0 & 1 \\ -1 & 0 \end{bmatrix}$ for $\mf{k}=\mf{sp}_{2\ell}$.
On the other hand, the action of $G_\ell$ in Theorem \ref{thm:classical Howe duality} (see also \cite[Chapter 5]{CW12}) is defined by lifting the action of the Lie algebra $\mf{x}=\mf{so}_\ell, \mf{sp}_{2\ell}$ given by 
$\mf{x} =\{\,X\,|\,X\in \mf{sl}_{r}, X^TJ_2+J_2X={\mathbb O}\,\}$, where $J_2$ is
\begin{equation}
\texttt{J}_\ell, \quad 
\begin{bmatrix}
0 & \texttt{J}_\ell \\
-\texttt{J}_\ell & 0 
\end{bmatrix},
\end{equation} 
for $\mf{so}_\ell$ and $\mf{sp}_{2\ell}$, respectively.
Here $\texttt{J}_\ell=(a_{ij})_{1\le i,j\le \ell}$ denotes the matrix given by $a_{ij}=\de_{i+j-\ell-1,0}$.
Let $\xi$ be an isomorphism of Lie algebras given by
\begin{equation}\label{eq:Xi}
  \xymatrixcolsep{3pc}\xymatrixrowsep{0pc}\xymatrix{ 
 \mf{x}\ \ar@{->}[r]^{\xi} & \  \mf{k} \\
  X  \ \ar@{|->}[r] & PXP^{-1} \ 
}, 
\end{equation}
where $P$ is the complex matrix such that $PJ_1P^{-1}=J_2$. So we may also regard $V_{G_\ell}(\la)$ as a $\mf{k}$-module by $\xi^{-1}$.%\comments{This paragraph added (20250413)}

\subsection{$(\U_X(\e),{\bf U}^\imath(\mf{k}))$-action on $\ms{W}^{\ot\ell}$}
Let $(\mf{g},\mf{k})$ be the symmetric pair in Section \ref{subsubsection:QSP}. 
Let us define a ${\bf U}(\mf{g})$-action on $\ms{W}^{\ot\ell}$, where
\begin{equation*}
 \ms{W}=
\begin{cases}
 \W & \text{for $\g=\mf{sl}_\ell$},\\
 \W^{2} & \text{for $\g=\mf{sl}_{2\ell}$},
\end{cases}
\end{equation*}
so that $\ms{W}^{\ot\ell}$ becomes a ${\bf U}^\imath(\mf{k})$-module, where 
\begin{equation}\label{eq:k convention}
 \mf{k}=
\begin{cases}
 \mf{so}_\ell & \text{if $\ms{W}=\W$},\\
 \mf{sp}_{2\ell} & \text{if $\ms{W}=\W^{2}$},
\end{cases}
\end{equation}
(cf.~\eqref{eq:G convention}).

From now on, we assume that
\begin{equation}\label{eq:twist by psi}
\text{the $\U_D(\e)$-module $\ms{W}$ in Proposition \ref{prop:osc-D} is twisted by $\psi$,} 
\end{equation}
where $\psi$ is an involution of $\U_D(\e)$ as a $\Bbbk$-algebra such that $\psi(f_0)=-f_0$, $\psi(k_0)=-k_0$, and $\psi(x)=x$ for the other generators $x$ of $\U_D(\e)$. 
This assumption is necessary when we consider the classical limit of  $\ms{W}^{\ot \ell}$ as a ${\bf U}^\imath(\mf{k})$-module, more precisely for $(\varsigma_i)_{i\in \msf{I}_\circ}$ to be specializable, that is, $\varsigma_i \in {\bf A}$ and $\varsigma_i(1)=-1$ for all $i\in \mathsf{I}_\circ$. Note that the weight spaces for the classical limit of $\ms{W}^{\ot \ell}$ as an $\mf{osp}_X(\e)$-module is twisted by $\psi$.

Suppose that $\ms{W}=\W$.
For $1\le t\le n$, let $\W'_t =\W_{\e^{(t)}}$, where $\e^{(t)}=(\e_t,\dots\e_t)=(\e_t^\ell)$. 
For $\lvert {\bf m}_1 \rangle \ot \lvert {\bf m}_2 \rangle \ot \cdots \ot \lvert {\bf m}_{\ell} \rangle \in \W^{\ot \ell}$, we write ${\bf m}_s=(m_{s1},\dots,m_{s n})$ for $1\le s\le \ell$.
Then we have a $\Bbbk$-linear isomorphism
\begin{equation}\label{eq:transpose iso}
 \xymatrixcolsep{3pc}\xymatrixrowsep{0pc}\xymatrix{ 
 \ms{W}^{\ot \ell} = \W^{\ot \ell} \ \ar@{->}[r] & \ \W'_1\ot \dots\ot \W'_n\\
  \lvert {\bf m}_1 \rangle \ot  \cdots \ot \lvert {\bf m}_\ell \rangle  \ \ar@{|->}[r] & \lvert {\bf m}^1 \rangle \ot \cdots \ot \lvert {\bf m}^n \rangle \ 
 }, 
\end{equation}
where ${\bf m}^t=(m_{1t},\dots,m_{\ell t})$ for $1\le t\le n$. One may consider an $(\ell\times n)$-matrix ${\bf M}=(m_{st})$ whose $s$-th row is ${\bf m}_s$ and $t$-th column is ${\bf m}^t$. So the action of $\U_X(\e)$ is given on ${\bf M}$ row by row from top to bottom via comultiplication.

For each $1\le t\le n$, we define a ${\bf U}(\mf{sl}_\ell)$-module structure on $\W'_t$ following \eqref{eq:type A action} with $\e=(0^\ell)$ (not depending on $\e^{(t)}$). So $\W'_1\ot \dots \ot \W'_n$ is a ${\bf U}(\mf{sl}_\ell)$-module with respect to \eqref{eq:comult}. %\comments{* Even though $\e^{t}=(1^\ell)$, the ${\bf U}(\mf{sl}_\ell)$-action is given as in $\W_{(0^\ell)}$. Should be confirmed by Bae. (250330)}
Thus we may define a ${\bf U}(\mf{sl}_\ell)$-action on $\ms{W}^{\ot\ell}$ by \eqref{eq:transpose iso}, which induces a ${\bf U}^\imath(\mf{k})$-action on $\ms{W}^{\ot\ell}$.
When $\ms{W}=\W^{\ot 2}$, we define a ${\bf U}(\mf{sl}_{2\ell})$-action on $\ms{W}^{\ot\ell}$ in the same way except $\e^{(t)}=(\e_t^\ell)$ replaced by $\e^{(t)}=(\e_t^{2\ell})$.

\begin{prop}\label{prop:type A commuting actions}
The action of ${\bf U}(\mf{g})$ on $\ms{W}^{\ot \ell}$ commutes with that of $\U_A(\e)$.
\end{prop}
\pf See Appendix \ref{app:proof type A}.
\qed

\begin{rem}{\rm
Let ${\bf U}(\gl_{r})=\U_A((0^r))$. Let $r=\ell$ and $2\ell$ when $\ms{W}=\W$ and $\W^{2}$, respectively. As a $(\U_A(\e),{\bf U}(\gl_r))$-module, $\ms{W}^{\ot \ell}$ has a multiplicity-free decomposition into tensor product of irreducible polynomial representations of $\U_A(\e)$ and ${\bf U}(\gl_r)$ parametrized by partitions $\la=(\la_1,\dots,\la_r)$ in $\ms{P}_\e$.
} 
\end{rem}

\begin{prop}\label{prop:joint actions on Fock space}
We have the following:
\begin{itemize}
 \item[(1)]  the ${\bf U}^\imath(\mf{so}_\ell)$-action on $\ms{W}^{\ot \ell}$ commutes with that of $\U_X(\e)$ if and only if $\varsigma_i=-q^{-1}$ and $\kappa_i=0$ for $i\in \msf{I}$,
 
 \item[(2)]  the ${\bf U}^\imath(\mf{sp}_{2\ell})$-action on $\ms{W}^{\ot \ell}$ commutes with that of $\U_X(\e)$ if and only if $\varsigma_i=-q$ for $i\in \msf{I}_\circ$.

\end{itemize}

\end{prop}
\pf See Appendix \ref{app:proof type DC}.
\qed
\medskip

%\red{
%\begin{rem}
%Actually, \ref{iqg class lim} holds for any subspace of $U$ instead of $U^\imath$.
%\end{rem}\comments{?}
%}
%\section{Howe duality on $\ms{W}^{\ot \ell}$}

%\subsection{Semisimplicity of $\ms{W}^{\ot \ell}$ as a ${\bf U}^\imath(\mf{k})$-module}

\subsection{$(\U_X(\e),{\bf U}^\imath(\mf{k}))$-duality on $\ms{W}^{\ot\ell}$}

For $\la\in \mc{P}(G_\ell)_\e$, let 
\begin{equation}\label{eq:mult space}
\mc{V}^\imath_{\mf{k}}(\la)={\rm Hom}_{\U_X(\e)}\left(\mc{V}^\la,\ms{W}^{\ot \ell}\right)
\end{equation}
be the multiplicity space of $\mc{V}^\la$ in $\ms{W}^{\ot \ell}$, which is a ${\bf U}^\imath(\mf{k})$-module by Proposition \ref{prop:joint actions on Fock space}. %
We have an isomorphism of $(\U_X(\e),{\bf U}^\imath(\mf{k}))$-modules:
\begin{equation}\label{eq:decomp iso}
  \xymatrixcolsep{3pc}\xymatrixrowsep{0pc}\xymatrix{ 
\bigoplus_{\la\in \mc{P}(G_\ell)_{\e}} \mc{V}^\la \ot \mc{V}^\imath_{\mf{k}}(\la) \ \ar@{->}[r] & \  \ms{W}^{\ot \ell} \\
  v\ot {\rho}  \ \ar@{|->}[r] & {\rho}(v) \ 
 },
\end{equation}
for $v\in \mc{V}^\la$ and $\rho\in \mc{V}^\imath_{\mf{k}}(\la)$. 
Let $v_\la$ be a highest weight vector of $\mc{V}^\la$, and let $B^\la=\{\,v^1,\dots,v^d\,\}$ with $d=\dim V_{G_\ell}(\la)$ be a $\Bbbk$-linearly independent subset in $\ms{W}^{\ot\ell}$ consisting of highest weight vectors of weight $\La_{\la,\e}$ \eqref{eq:highest weight La} (cf.~Proposition \ref{prop:irr comp in Fock space}). 
By \eqref{eq:decomp iso}, we identify the $\Bbbk$-span of $B^\la$ with $\mc{V}^\imath_{\mf{k}}(\la)$ as a ${\bf U}^\imath(\mf{k})$-module, where $v^k$ corresponds to a $\U_X(\e)$-linear map $\rho_k$ sending $v_\la$ to $v^k$ ($k=1,\dots,d$).

\begin{lem}\label{lem:semisimplicity of mult space}
For $\la\in \mc{P}(G_\ell)_\e$, $\mc{V}^\imath_{\mf{k}}(\la)$ is a semisimple ${\bf U}^\imath(\mf{k})$-module.
\end{lem}
\pf For ${\bf m}=(m_1,\dots,m_n)\in \Z^n_+(\e)$, let $|{\bf m}|=m_1+\dots+m_n$. For  $d\in \Z_{\ge 0}$, let $(\ms{W}^{\ot \ell})_d$ be the $\Bbbk$-span of $|{\bf m}_1\rangle \ot\dots \ot |{\bf m}_r\rangle$ such that $|{\bf m}_1|+\dots+ |{\bf m}_r|=d$, where $r=\ell$ or $2\ell$.

Then $\mc{V}^\imath_{\mf{k}}(\la)$, the $\Bbbk$-span of $B^\la$ is a subspace of $(\ms{W}^{\ot \ell})_d$, where $d=\la_1+\dots+\la_\ell$ (see \eqref{eq:wt of m} and \eqref{eq:highest weight La}). 
Since $(\ms{W}^{\ot \ell})_d$ is a finite-dimensional ${\bf U}(\mf{g})$-submodule of $\ms{W}^{\ot \ell}$, we have by \cite[Proposition 3.1.4]{Wat21} that $\mc{V}^\imath_{\mf{k}}(\la)$ is a finite-dimensional {\em classical weight module} in the sense of \cite[Section 3.1]{Wat21}.
Hence it follows from \cite[Corollaries 4.1.7, 4.2.4, 4.3.4]{Wat21} that  $\mc{V}^\imath_{\mf{k}}(\la)$ is semisimple.
\qed
\medskip

Let us consider the classical limit of \eqref{eq:decomp iso}.
Let ${\bf U}_{\bf A}$ be the ${\bf A}$-subalgebra of ${\bf U}={\bf U}(\mf{g})$ generated by $\mathsf{e}_i, \mathsf{f}_i, \mathsf{k}_i^{\pm 1}$, and $(\mathsf{k}_i ; 0)_q:= \frac{\mathsf{k}_i-1}{q-1}$ for $i\in \mathsf{I}$. 
Let ${\bf U}_1={\bf U}_{\bf A}\ot_{\bf A}\C $, where $\C$ is an ${\bf A}$-module given by evaluating $f(q^{\hf})\in {\bf A}$ at 1. 
Then there is an isomorphism of $\C$-algebras $\pi: {\bf U}_1 \rightarrow U(\mf{g})$ given by $\pi({\mathsf{e}_i\ot 1})= \ov{\mathsf{e}}_i$, $\pi({\mathsf{f}_i\ot 1}) = \ov{\mathsf{f}}_i$, and $\pi({(\mathsf{k}_i ; 0)_q\ot 1})= \ov{\mathsf{h}}_i$ for $i\in \msf{I}$, where $\ov{\mathsf{e}}_i, \ov{\mathsf{f}}_i, \ov{\mathsf{h}}_i$ denote the generators of $\mf{g}$. %, where $\overline{\mathsf{x}}=\mathsf{x}\ot 1\in {\bf U}_1$ for $\mathsf{x}=\mathsf{e}_i, \mathsf{f}_i, \mathsf{k}_i$, and $e_i, f_i, h_i$ denote the generators of $\mf{sl}_r$.
Let ${\bf U}^\imath_{\bf A}={\bf U}^\imath(\mf{k}) \cap {\bf U}_{\bf A}$ and ${\bf U}^\imath_1= {\bf U}^\imath_{\bf A} \ot_{\bf A}\C \subset {\bf U}_1$. If $(\varsigma_i)_{i\in \msf{I}_\circ}$ are specializable, that is, $\varsigma_i \in {\bf A}$ and $\varsigma_i(1)=-1$ for all $i\in \mathsf{I}_\circ$, then the map $\pi$ restricts to an isomorphism ${\bf U}^\imath_1 \rightarrow {U}(\mf{k})$ {\cite[Theorem 10.8]{Kol14}}.
Let $M$ be a ${\bf U}^\imath(\mf{k})$-module with an ${\bf A}$-lattice $M_{\bf A}$ stable under ${\bf U}^\imath_{\bf A}$. Then $\ov{M_{\bf A}}=M_{\bf A}\ot_{\bf A} \C$ becomes a ${\bf U}^\imath_1$-module, which we also call a classical limit of $M$.

Now, let $\ms{W}_{\bf A}^{\ot \ell}$ be the ${\bf A}$-submodule of $\ms{W}^{\ot \ell}$ in Section \ref{subsec:classical Howe duality}. Then it is clearly stable under the action of ${\bf U}^\imath_{\bf A}$. By Proposition \ref{prop:joint actions on Fock space}, its classical limit $W^{\ot \ell}$ \eqref{eq:classical limit of W-2} becomes a $(U(\mf{osp}_X(\e)),U(\mf{k}))$-module (cf.~\eqref{eq:osp_X}).%\comments{It may not coincide with the actions in Theorem \ref{thm:classical Howe duality} (250308)}  

Let $\la\in \mc{P}(G_\ell)_\e$ be given. Let 
\begin{equation}\label{eq:mult space A}
\mc{V}^\imath_{\mf{k}}(\la)_{\bf A}=
{\rm Hom}_{U_X(\e)_{\bf A}}
\left(
(\mc{V}^\la)^\tau_{\bf A},
\left(\ms{W}^{\ot \ell}_{\bf A}\right)^\tau
\right),
\end{equation}
where $U_X(\e)_{\bf A}$ is the ${\bf A}$-subalgebra of $U_X(\e)$ generated by $E_i$, $F_i$, and $[E_i,F_i]$ for $i\in I$, and $(\mc{V}^\la)^\tau_{\bf A}=U_X(\e)_{\bf A}v_\la$. 
Then it is stable under the action of ${\bf U}^\imath_{\bf A}$. Let 
\begin{equation*}
 {V}_{\mf{k}}^\imath(\la)=\ov{\mc{V}^\imath_{\mf{k}}(\la)_{\bf A}}.
\end{equation*}
If we choose a set $B^\la=\{\,v^1,\dots,v^d\,\}\subset \ms{W}_{\bf A}^{\ot \ell}$  such that $v^k$ is a highest weight vector of weight $\La_{\la,\e}$ and $\{\,v^k\ot 1\,|\,k=1,\dots,d\,\}\subset \ms{W}_{\bf A}^{\ot \ell}\ot_{\bf A} \C$ is linearly independent over $\C$, then the corresponding set $\{\,{\rho_k}\,|\,k=1,\dots,d\,\}\subset \mc{V}^\imath_{\mf{k}}(\la)$ forms an ${\bf A}$-basis of $\mc{V}^\imath_{\mf{k}}(\la)_{\bf A}$ (cf.~\eqref{eq:decomp iso}). 

Since each weight vector of $\left(\ms{W}^{\ot \ell}\right)^\tau$ as a $U_X(\e)$-module is contained only in a finitely many $(\mc{V}^\la)^\tau$'s,
we may conclude that the following map is an isomorphism of $({U_X(\e)}_{\bf A},{\bf U}^\imath_{\bf A})$-modules:
\begin{equation}
  \xymatrixcolsep{3pc}\xymatrixrowsep{0pc}\xymatrix{ 
\bigoplus_{\la\in \mc{P}(G_\ell)_{\e}} (\mc{V}^\la)^\tau_{\bf A} \ot \mc{V}^\imath_{\mf{k}}(\la)_{\bf A} \ \ar@{->}[r] & \  \left(\ms{W}^{\ot \ell}_{\bf A}\right)^\tau\\
  v\ot {\rho}  \ \ar@{|->}[r] & {\rho}(v) \ 
 }.
\end{equation}
Taking $\,\cdot\, \ot_{\bf A}\C$, we have by Corollary \ref{cor:classical limit of Vlambda}
\begin{equation}\label{eq:classical limit of decomp of Fock space}
{W}^{\ot \ell}\cong\bigoplus_{\la\in \mc{P}(G_\ell)_{\e}} {V}^\la \ot {V}^\imath_{\mf{k}}(\la),
\end{equation}
as a $({U(\mf{osp}_X(\e))},U(\mf{k}))$-module.

Let $\la\in \mc{P}(G_\ell)_\e$ be given. We have by Lemma \ref{lem:semisimplicity of mult space} 
$$\mc{V}^\imath_{\mf{k}}(\la)=V_1\oplus\dots\oplus V_s,$$ where $V_i$ is an irreducible ${\bf U}^\imath(\mf{k})$-module. To describe the above decomposition of $\mc{V}^\imath_{\mf{k}}(\la)$, let us briefly recall some of the results in \cite{Wat21}.

Let $\mathbb{K}$ be an algebraic closure of $\Bbbk$ in $\bigcup_{n=1}^\infty\C(\!(q^{\frac{1}{n}})\!)$.
Let $\mc{U}^\imath$ be the $\mathbb{K}$-algebra defined in \cite[Definition 3.1.13]{Wat21}, which includes ${\bf U}^\imath(\mf{k}) \ot_{\Bbbk}\mathbb{K}$ as a subalgebra.
Let $\mathbb{K}_1$ be its subring of $f\in \mathbb{K}$ regular at $1$. 
Let $\mc{U}^\imath_{\mathbb{K}_1}$ be its $\mathbb{K}_1$-form defined in \cite[Section 3.3]{Wat21} and ${\bf U}^\imath_{\mathbb{K}_1}={\bf U}^\imath_{\bf A}\ot_{\bf A}\mathbb{K}_1$. 
By \cite[Propositions 2.4.6, 3.3.1]{Wat21}, we have
\begin{equation}\label{eq:classical limits are equal}
 \ov{\mc{U}^\imath_{\mathbb{K}_1}} := \mc{U}^\imath_{\mathbb{K}_1}\ot_{\mathbb{K}_1}\C=
 \ov{{\bf U}^\imath_{\mathbb{K}_1}} = 
 \ov{{\bf U}^\imath_{\bf A}}\cong U(\mf{k}),
\end{equation} 
where we assume that $\C$ is a $\mathbb{K}_1$-module given by $f\cdot 1 = a_0$ where $a_0$ is the evaluation of $f$ at $q^{\frac{1}{n}}=1$ ($n\ge 1$).

The ${\bf U}^\imath(\mf{k})$-modules $\mc{V}^\imath_{\mf{k}}(\la)$ and $V_i$ ($i=1,\dots,s$) can be extended to $\mc{U}^\imath$-modules.
Let $\mc{V}^\imath_{\mf{k}}(\la)_{\mathbb{K}_1}=\mc{V}^\imath_{\mf{k}}(\la)_{\bf A}\ot_{\bf A}\mathbb{K}_1$. 
Note that 
\begin{equation}\label{eq:classical limit from K_1}
\ov{\mc{V}^\imath_{\mf{k}}(\la)_{\mathbb{K}_1}}:=
\mc{V}^\imath_{\mf{k}}(\la)_{\mathbb{K}_1}\ot_{\mathbb{K}_1}\C = \ov{\mc{V}^\imath_{\mf{k}}(\la)_{\bf A}}, 
\end{equation} 
as a $\C$-space. 

\begin{lem}\label{lem:classical limits of multiplicity space are equal}
We have $\ov{\mc{V}^\imath_{\mf{k}}(\la)_{\mathbb{K}_1}}=\ov{\mc{V}^\imath_{\mf{k}}(\la)_{\bf A}}$ as $U(\mf{k})$-modules. 
\end{lem}
\pf It follows from \eqref{eq:classical limit from K_1} and \eqref{eq:classical limits are equal}.
\qed
\medskip

By \cite[Sections 3.3, 4.1--4.3]{Wat21},
there exist finite subsets $\mc{X}_+, \mc{X}_-, \mc{X}_0\subset \mc{U}^\imath_{\mathbb{K}_1}$ such that 
\begin{itemize}
 \item[(1)] $\mc{U}^\imath_{\mathbb{K}_1}=\mc{U}^\imath_{-,\mathbb{K}_1}\mc{U}^\imath_{0,\mathbb{K}_1}\mc{U}^\imath_{+,\mathbb{K}_1}$ where $\mc{U}^\imath_{\pm,\mathbb{K}_1}$ is the $\mathbb{K}_1$-subalgebra generated by $\mc{X}_\pm$ and $\mc{U}^\imath_{0,\mathbb{K}_1}$ is a subalgebra including $\mc{X}_0$,
 
 \item[(2)] $\mf{k}=\C\ov{\mc{X}_+}\oplus \C\ov{\mc{X}_0}\oplus \C\ov{\mc{X}_-}\subset \ov{\mc{U}^\imath_{\mathbb{K}_1}}$ is a triangular decomposition of $\mf{k}$, where $\C\ov{\mc{X}_0}$, the $\C$-span of the image of ${\mc{X}_0}$, is the Cartan subalgebra and $\C\ov{\mc{X}_+}$ (resp. $\C\ov{\mc{X}_-}$) is the subalgebra spanned by the positive (resp. negative) root vectors.
\end{itemize}
Moreover, it is shown that every irreducible finite-dimensional classical weight module $M$ is generated by $v$ such that $\mc{X}_+v=0$, and $v$ is a simultaneous eigenvector with respect to $\mc{X}_0$ (a highest weight vector), while $\U^\imath_{0,\mathbb{K}_1}$ acts on $v$ semisimply.

\begin{lem}\label{lem:decomp of mult space}
Let $(V_i)_{\mathbb{K}_1}=V_i\cap \mc{V}^\imath_{\mf{k}}(\la)_{\mathbb{K}_1}$ for $i=1,\dots,s$. 
Then we have
\begin{itemize}
 \item[(1)] $\mc{V}^\imath_{\mf{k}}(\la)_{\mathbb{K}_1}=(V_1)_{\mathbb{K}_1}\oplus\dots\oplus (V_s)_{\mathbb{K}_1}$,
 
 \item[(2)] $\ov{(V_i)_{\mathbb{K}_1}}$ is an irreducible $U(\mf{k})$-module for $i=1,\dots,s$.
\end{itemize}
\end{lem}
\pf For each $i$, let $v_i$ be a highest weight vector of $V_i$. 
We may assume that $v_i\in (V_i)_{\mathbb{K}_1}$ and $\ov{v_i}\in \ov{(V_i)_{\mathbb{K}_1}}$ is non-zero.
From the triangular decomposition of $\U^\imath_{\mathbb{K}_1}$, we see that 
$(V_i)_{\mathbb{K}_1}\supset \mc{U}^\imath_{-,\mathbb{K}_1} \, v_i$ is finitely generated over $\mathbb{K}_1$ and hence free over $\mathbb{K}_1$. This implies that $\dim_{\mathbb{K}}V_i={\rm rank}_{\mathbb{K}_1}(\mc{U}^\imath_{-,\mathbb{K}_1} \, v_i)$ (cf.~\cite[Proposition 3.3.9]{Wat21}).
Since $\mc{V}^\imath_{\mf{k}}(\la)_{\mathbb{K}_1}$ is also $\mathbb{K}_1$-free with $\dim_{\mathbb{K}}\mc{V}^\imath_{\mf{k}}(\la)={\rm rank}_{\mathbb{K}_1}\mc{V}^\imath_{\mf{k}}(\la)_{\mathbb{K}_1}$, we apply Nakayama lemma to have $\mc{V}^\imath_{\mf{k}}(\la)_{\mathbb{K}_1}=\bigoplus_{i=1}^s\mc{U}^\imath_{-,\mathbb{K}_1} \, v_i=\bigoplus_{i=1}^s (V_i)_{\mathbb{K}_1}$ with $\mc{U}^\imath_{-,\mathbb{K}_1} \, v_i=(V_i)_{\mathbb{K}_1}$. This proves (1). The proof of (2) is easy.
\qed

\begin{rem}{\rm
Indeed,  due to \cite[Propositions 4.1.6, Corollaries 4.2.2 and 4.3.2]{Wat21}, one may replace $\mathbb{K}_1$ in Lemma \ref{lem:decomp of mult space} with ${\bf A}$.
} 
\end{rem}
%\comments{Revised (20250511)}

Now we have the following, which can be viewed as a $q$-analogue of Theorem \ref{thm:classical Howe duality}.

\begin{thm}\label{thm:main result}\mbox{}
\begin{itemize}
 \item[(1)] 
As a $(\U_X(\e),{\bf U}^\imath(\mf{k}))$-module, we have
\begin{equation*}
\ms{W}^{\ot \ell} = 
\bigoplus_{\la\in \mc{P}(G_\ell)_{\e}} \V^\la \ot \mc{V}^\imath_{\mf{k}}(\la).
\end{equation*} 
 
 \item[(2)] The classical limit \eqref{eq:classical limit of decomp of Fock space} of $\ms{W}^{\ot \ell}$ is isomorphic to \eqref{eq:osp Howe duality} as a $({U(\mf{osp}_X(\e))},U(\mf{x}))$-module when twisted by $\xi$ \eqref{eq:Xi}. In particular, the classical limit of $\mc{V}^\imath_{\mf{k}}(\la)$ is isomorphic to $V_{G_\ell}(\la)$ as a $U(\mf{x})$-module.
 
 \item[(3)] Let $\la\in \mc{P}(G_\ell)_{\e}$ be given. For $G_\ell=Sp_{2\ell}$, $\mc{V}^\imath_{\mf{k}}(\la)$ is irreducible. For $G_\ell=O_\ell$, $\mc{V}^\imath_{\mf{k}}(\la)$ is irreducible if $\ell$ is odd or $\ell$ is even with $\la_\ell=0$, and it is a sum of two irreducible ${\bf U}^\imath(\mf{k})$-modules otherwise.
\end{itemize}
%We have as a $(\U_X(\e),{\bf U}^\imath(\mf{k}))$-module
%\begin{equation*}
%\ms{W}^{\ot \ell} = 
%\bigoplus_{\la\in \mc{P}(G_\ell)_{\e}} \V^\la \ot \mc{V}^\imath_{\mf{k}}(\la),
%\end{equation*} 
%where the classical limit ${V}^\imath_{\mf{k}}(\la)$ of $\mc{V}^\imath_{\mf{k}}(\la)$ is isomorphic to $V_{G_\ell}(\la)$ as a $U(\mf{k})$-module.
% given as follows:
\end{thm}
\pf (1) It follows from \eqref{eq:decomp iso}.

(2) For simplicity, let $\mc{A}=U(\mf{osp}_X(\e))$, $\mc{B}_1=U(\mf{x})$, and $\mc{B}_2=U(\mf{k})$. Let $\mathbb{W}^\ell_1$ and $\mathbb{W}^\ell_2$ denote the $(\mc{A},\mc{B}_1)$ and $(\mc{A},\mc{B}_2)$-modules in \eqref{eq:osp Howe duality} and \eqref{eq:classical limit of decomp of Fock space}, respectively. We may also regard $\mathbb{W}^\ell_2$ as an $(\mc{A},\mc{B}_1)$-module by twisting the action of $\mc{B}_2$ by $\xi$  \eqref{eq:Xi}.

It is enough to show that they have the same characters since they are semisimple $(\mc{A},\mc{B}_1)$-modules with irreducible components $V^\la\ot U$ for $\la\in \mc{P}(G_\ell)_\e$ and a finite-dimensional irreducible $\mc{B}_1$-module $U$, and the dimension of each weight space with respect to the action of $(\mc{A},\mc{B}_1)$ is finite-dimensional.

Let $\alpha=\{\,v_1,\dots,v_r\,\}$ be the standard basis of $\C^r$, the natural representation of $\mf{sl}_r$ ($r=\ell, 2\ell$).
For $\mf{k}=\mf{x}=\mf{so}_\ell$, let $\beta=\{\,w_{\pm s}\,|\,s=1,\dots,n\,\}$ for $\ell=2n$ and $\beta=\{\,w_{\pm s}\,|\,s=1,\dots,n\,\}\cup\{\,w_0\,\}$ for $\ell=2n+1$, where
\begin{equation*}
w_{\pm s}=\frac{1}{\sqrt{2}}(v_s \pm \sqrt{-1}v_{\ell-s+1}),\quad 
w_{0}= v_{n+1}.
\end{equation*}

For $\mf{k}=\mf{x}=\mf{sp}_{2\ell}$, let $\beta=\{\,w_{\pm s}\,|\,s=1,\dots,\ell\,\}$, where $w_s=v_{\pi(s)}$ and $w_{-s}=v_{\pi(2\ell-s+1)}$ for $s=1,\dots,\ell$, and $\pi$ is the permutation on $2\ell$ letters given by $\pi(i)=2i-1$ ($i=1,\dots,\ell$) and $\pi(i)=2(2\ell-i+1)$ for $i=\ell+1,\dots,2\ell$. %\comments{Corrected (20250518)} 
Then the matrix $P$ in \eqref{eq:Xi} is the transition matrix $[{\rm id}]_\beta^\alpha$ in both cases.

Let $\mathbb{U}_i$ denote the representation of $\mc{B}_i$ ($i=1,2$) given by restriction of the $\mf{sl}_r$-module $\C^r$.
Then $\mathbb{U}_1$ is the first fundamental representation of $\mc{B}_1$.
If we denote by $\mathbb{U}_2^\xi$ the $\mc{B}_1$-module twisted by $\xi$, then we can check without difficulty that
\begin{equation}\label{eq:natural twisted}
\mathbb{U}_2^\xi \cong \mathbb{U}_1,
\end{equation}
where $\beta$ is a basis of $\mathbb{U}_2^\xi$ consisting of weight vectors.
Note that $\mathbb{W}_i^\ell\cong \mathbb{V}_i^1 \ot \dots\ot \mathbb{V}_i^n$ as a $\mc{B}_i$-module ($i=1,2$), where each $\mathbb{V}_i^t$ is either the symmetric algebra $S(\mathbb{U}_i)$ or the exterior algebra $\Lambda(\mathbb{U}_i)$ generated by $\mathbb{U}_i$ (depending on the parity of $\e_t$).
Hence $\mathbb{W}^\ell_2$ when twisted by $\xi$ is isomorphic to $\mathbb{W}_1^\ell$ as a $\mc{B}_1$-module by \eqref{eq:natural twisted}.

Since an element in the monomial basis $\beta^\ell$ of $\mathbb{W}^\ell_2$ generated by $\beta$ is also a weight vector with respect to the action of $\mc{A}$, we can check directly that the character of $\mathbb{W}^\ell_2$ as an $(\mc{A},\mc{B}_1)$-module with respect to $\beta^\ell$ is equal to that of $\mathbb{W}^\ell_1$ given as a product form in \cite[(6.11), (6.12)]{CZ} for $G_\ell=O_\ell$ with $\varepsilon=1$ and \cite[(6.17)]{CZ} for $G_\ell=Sp_{2\ell}$.
So we have $\mathbb{W}^\ell_1 \cong \mathbb{W}^\ell_2$ as $(\mc{A},\mc{B}_1)$-modules, which in particular implies that ${V}^\imath_{\mf{k}}(\la)=\ov{\mc{V}^\imath_{\mf{k}}(\la)_{\bf A}}$ is isomorphic to $V_{G_\ell}(\la)$ as a $\mc{B}_1$-module, equivalently as a $\mc{B}_2$-module by \eqref{eq:Xi}.

(3) It follows from Lemma \ref{lem:decomp of mult space} and (2).
\qed 

\appendix

\section{Proof of Proposition \ref{prop:type A commuting actions}}\label{app:proof type A}

Since $\W^2$ is isomorphic to $\W^{\ot 2}$ as a $\U_A(\e)$-module, which is induced from \eqref{eq:type A action} by comultiplication, we may assume that $\ms{W}=\W=\W_\e$. According to \eqref{eq:transpose iso}, we may identify both of $\lvert {\bf m}_1 \rangle \ot  \cdots \ot \lvert {\bf m}_\ell \rangle\in \W^{\ot \ell}$ and $\lvert {\bf m}^1 \rangle \ot \cdots \ot \lvert {\bf m}^n \rangle\in \W'_1\ot \dots\ot \W'_n$ with ${\bf M}=(m_{st})$ for $1\le s\le \ell$ and $1\le t\le n$.

Let $i\in I\setminus\{0\}$ and $j\in \msf{I}$ be given. 
It is enough to show that $x_i \msf{x}_j {\bf M} = \msf{x}_j x_i {\bf M}$ for $x=e,f$ and $\msf{x}=\msf{e}, \msf{f}$. Let us consider the case when $x=e$ and $\msf{x}=\msf{e}$ since the other cases can be checked similarly. 

As operators on $\W^{\ot \ell}$ and $\W'_1\ot \dots\ot \W'_n$, 
we have $e_i=\sum_{s=1}^\ell e_{i,s}$ and $\msf{e}_j=\sum_{t=1}^n \msf{e}_{j,t}$, respectively, where  $e_{i, s}:=k_i^{\ot (s-1)} \ot e_i \ot 1^{\ot (\ell-s)}$ for $1\le s\le \ell$ and $\msf{e}_{j, t}$ is defined similarly for $1\le t\le n$. 
Note that $\msf{e}_{j,t} e_{i,s} {\bf M}=e_{i,s} \msf{e}_{j,t} {\bf M}$ if { either $s \neq j, j+1$ or $t \neq i, i+1$}, since $e_{i,s}$ acts on the $i$-th and $(i+1)$-th {columns} and $\msf{e}_{j,t}$ acts on the $j$-th and $(j+1)$-th {rows} of ${\bf M}$. 
Thus it is enough to show that
\begin{equation*}
(e_{i,j}+e_{i,j+1})(\msf{e}_{j, i}+\msf{e}_{j, i+1}) {\bf M} = (\msf{e}_{j, i}+\msf{e}_{j, i+1})(e_{i,j}+e_{i,j+1}) {\bf M}, 
\end{equation*}
that is, $e_i \msf{e}_j {\bf M}' = \msf{e}_j e_i  {\bf M}' $, where ${\bf M}'$ is the $2 \times 2$ submatrix
\begin{equation*}
 {\bf M}'=
\begin{bmatrix}
x & y \\ z & w
\end{bmatrix}
=
\begin{bmatrix}
m_{j i} & m_{j i+1} \\ 
m_{j+1 i} & m_{j+1 i+1}
\end{bmatrix},
\end{equation*}
and the action of $e_i$ (resp. $\msf{e}_j$) is given on each row (resp. column).

We consider the case when $(\e_i,\e_{i+1})=(0,1)$, and leave the other cases to the reader. We have
{\small
\begin{align*}\label{eq:aux-1}
& e_i \msf{e}_j {\bf M}' \\
&= e_i \Bigg\{ [z] 
\begin{bmatrix}
x+1 & y \\ z-1 & w
\end{bmatrix} 
+ q^{x-z} \de_{y,0}\de_{w,1} 
\begin{bmatrix}
x & y+1 \\ z & w-1
\end{bmatrix} \Bigg\} \\
&= [z] \de_{y,1} 
\begin{bmatrix}
x+2 & y-1 \\ z-1 & w
\end{bmatrix} 
+ ([z]q_i^{x+1}q_{i+1}^{-y} \de_{w,1} + {q^{x-z}\de_{y,0}\de_{w,1}})
\begin{bmatrix}
x+1 & y \\
z & w-1
\end{bmatrix},
\end{align*}}%\comments{Bae: removed $\de_{z\geq 1}$ (20250801)}
and
{\small
\begin{align*}
&\msf{e}_je_i {\bf M}' \\
&= 
\msf{e}_j \Bigg\{ \de_{y,1} 
\begin{bmatrix}
x+1 & y-1 \\ z & w
\end{bmatrix} 
+ q_i^x q_{i+1}^{-y} \de_{w,1} 
\begin{bmatrix}
x & y \\ z+1 & w-1
\end{bmatrix} \Bigg\} \\
&= 
\de_{y,1} \Bigg\{ [z] 
\begin{bmatrix}
x+2 & y-1 \\ z-1 & w
\end{bmatrix} 
+ q^{x+1-z}\de_{w,1} 
\begin{bmatrix}
x+1 & y \\ z & w-1
\end{bmatrix} \Bigg\} \\
&\quad + q_i^xq_{i+1}^{-y}\de_{w,1}[z+1] 
\begin{bmatrix}
x+1 & y \\ z & w-1
\end{bmatrix}.
\end{align*}}%\comments{Bae: removed $\de_{z\geq 1}$ (20250801)}
It is clear that the coefficients of  
{\small $\begin{bmatrix}
x+2 & y-1 \\ z-1 & w
\end{bmatrix}$}
in $e_i\msf{e}_j {\bf M}$ and $\msf{e}_je_i {\bf M}$ are the same, while those of 
{\small $\begin{bmatrix}
x+1 & y \\ z & w-1
\end{bmatrix}$} are
\begin{equation*}
\begin{split}
A:&=[z]q_i^{x+1}q_{i+1}^{-y} \de_{w,1} + {q^{x-z}\de_{y,0}\de_{w,1}},\\
B:&=q^{x+1-z} \de_{y,1}\de_{w,1} + q_i^x q_{i+1}^{-y}[z+1] \de_{w,1},
\end{split} 
\end{equation*}%\comments{Bae: removed $\de_{z\geq 1}$ (20250801)}
respectively. It is straightforward to check that $A=B$.
\qed

\section{Proof of Proposition \ref{prop:joint actions on Fock space}}\label{app:proof type DC}
We use the notation ${\bf M}=(m_{st})$ as in Appendix \ref{app:proof type A}.
By Proposition \ref{prop:type A commuting actions}, it suffices to check 
\begin{align}
x_0\msf{B}_i {\bf M} = \msf{B}_ix_0 {\bf M},\label{aux:comm-1}\\
x_0\msf{x}_j {\bf M} = \msf{x}_jx_0 {\bf M},\label{aux:comm-2}
\end{align}
for $x_0\in \U_X(\e)$ ($x=k, e, f$), $\msf{B}_i\in {\bf U}^\imath(\mf{k})$ ($i\in \msf{I}_\circ$) in \eqref{eq:B_i for so} and \eqref{eq:B_i for sp}, and $\msf{x}_j$ ($\msf{x}=\msf{k}, \msf{e},\msf{f}$, $j\in \msf{I}_\bullet$). 
Let us prove \eqref{aux:comm-1} in the case of $x=e$ since the proof for the case of $x=f$ is similar. The proofs of \eqref{aux:comm-1} with $x=k$ and \eqref{aux:comm-2} are straightforward.
Note that as an operator on $\ms{W}^{\ot \ell}$, we have $e_0=\sum_{s=1}^\ell e_{0,s}$, where  $e_{0, s}=k_0^{\ot (s-1)} \ot e_0 \ot 1^{\ot (\ell-s)}$ for $1\le s\le \ell$. Let $\lang k\rang=\de_{\e_2,0}+\de_{\e_2,1}\de_{k,0}$ for $k\in \Z_{\geq 0}$. %\comments{Bae: added $\lang k\rang$ which will be used to handle both cases of $\e_2=0$ and $\e_2=1$ together (20250724)}

\subsection{Proof of Proposition \ref{prop:joint actions on Fock space}(1)}
It is clear that $\msf{B}_i e_{0,s} {\bf M}=e_{0,s} \msf{B}_i {\bf M}$ for $s \neq i, i+1$. So, it remains to show that
\begin{equation}\label{eq:commutativity for so}
 \msf{B}_i(e_{0, i}+e_{0, i+1}) {\bf M} = (e_{0, i}+e_{0, i+1})\msf{B}_i {\bf M},
\end{equation}

{\em Case 1.} Suppose that $X=D$. %\comments{Revised (250325)}
As an operator on $\ms{W}^{\ot \ell}$, let us write \eqref{eq:B_i for so} as
\begin{equation}\label{eq:B_i as operator on tensor for so}
\msf{B}_i=\sum_{1\le t\le n}\msf{f}_{i,t}+\varsigma_i \sum_{1\le t\le n}\msf{e}_{i,t} (\msf{k}_i^{-1})^{\ot n},
\end{equation}
where $\msf{f}_{i,t}= 1^{\ot (n-t)} \ot \msf{f}_i \ot (\msf{k}_0^{-1})^{\ot (t-1)}$ and $\msf{e}_{i,t}=(\msf{k}_0)^{\ot (t-1)} \ot \msf{e}_i \ot 1^{\ot (n-t)}$.
First, one can check without difficulty that for $t\ge 3$
\begin{equation*}
\begin{split}
 \msf{f}_{i,t}(e_{0, i}+e_{0, i+1}) {\bf M} &= (e_{0, i}+e_{0, i+1})  \msf{f}_{i,t}{\bf M},\\
  \msf{e}_{i,t}(\msf{k}_i^{-1})^{\ot n}(e_{0, i}+e_{0, i+1}) {\bf M} &= (e_{0, i}+e_{0, i+1})  \msf{e}_{i,t}(\msf{k}_i^{-1})^{\ot n}{\bf M}.
\end{split}
\end{equation*}%\comments{Will be checked again by Bae}
Hence in order to prove \eqref{eq:commutativity for so}, it is enough to prove $e_0\msf{B}_i {\bf M}' = \msf{B}_ie_0 {\bf M}'$, where ${\bf M}'$ is the $2 \times 2$ submatrix 
\begin{equation*}
 {\bf M}'=
\begin{bmatrix}
x & y \\ z & w
\end{bmatrix}
=
\begin{bmatrix}
m_{j 1} & m_{j 2} \\ 
m_{j+1 1} & m_{j+1 2}
\end{bmatrix},
\end{equation*}
and the action of $e_0$ (resp. ${\bf B}_i$) is given on each row (resp. column).
%We consider the case when $(\e_1, \e_2)=(1, 0)$, since the proof for the case of $(\e_1, \e_2)=(1, 1)$ are almost the same.
We have
{\small{\allowdisplaybreaks
\begin{align*}
&e_0 \msf{B}_i 
\begin{bmatrix}
x &y \\ z&w
\end{bmatrix} \\ 
&= e_0 \Bigg[ \de_{z,0}[x] q^{-y+w} 
\begin{bmatrix}
x-1 & y \\ z+1 & w
\end{bmatrix} +[y]\lang w\rang 
\begin{bmatrix}
x & y-1 \\ z & w+1
\end{bmatrix} \\
&\quad + \varsigma_i q^{-x-y+z+w} \Bigg\{ \de_{x,0}[z] \begin{bmatrix}
x+1 & y \\ z-1 & w
\end{bmatrix} + q^{x-z}[w]\lang y\rang \begin{bmatrix}
x & y+1 \\ z & w-1
\end{bmatrix} \Bigg\} + \kappa_i q^{z+w-x-y}
\begin{bmatrix}
x & y \\
z & w
\end{bmatrix}  \Bigg] \\
&=- \de_{z,0}[x] q^{-y+w} \td{q}^{-x+2}q_2^{-y}[w] \begin{bmatrix}
x-1 & y \\ z & w-1
\end{bmatrix} \\
& \quad + [y]\lang w\rang \Bigg\{ [y-1] \de_{x,1} \begin{bmatrix}
x-1 & y-2 \\ z & w+1
\end{bmatrix} - \td{q}^{-x+1}q_2^{-y+1} \de_{z,1}[w+1] \begin{bmatrix}
x & y-1 \\
z-1 & w
\end{bmatrix} \Bigg\} \\
& \quad +\varsigma_i q^{-x-y+z+w}\de_{x,0}[y][z] \begin{bmatrix}
x & y-1 \\
z-1 & w
\end{bmatrix} \\
& \quad + \varsigma_i q^{-y+w}[w]\lang y \rang \Bigg\{ \de_{x,1}[y+1]\begin{bmatrix}
x-1 & y \\
z& w-1
\end{bmatrix} - \td{q}^{-x+1}q_2^{-y-1} \de_{z,1}[w-1] \begin{bmatrix}
x & y+1 \\
z-1 & w-2
\end{bmatrix} \Bigg\} \\
& \quad + \kappa_i q^{z+w-x-y} \Bigg\{ \de_{x,1}[y] 
\begin{bmatrix}
x-1 & y-1 \\
z & w
\end{bmatrix}
 - \td{q}^{1-x}q^{-y} \de_{z,1}[w] 
\begin{bmatrix}
x & y \\
z-1 & w-1
\end{bmatrix} \Bigg\},
\end{align*}}} %\comments{Bae added constants $\lang k \rang$ (20250724)}

and
{\small{\allowdisplaybreaks
\begin{align*}
&\msf{B}_i e_0 
\begin{bmatrix}
x & y \\
z & w
\end{bmatrix}  \\
& = \msf{B}_i \Bigg\{ \de_{x,1}[y] \begin{bmatrix}
x-1 & y-1 \\
z & w
\end{bmatrix} - \td{q}^{1-x}q^{-y}\de_{z,1}[w] \begin{bmatrix}
x & y \\
z-1 & w-1
\end{bmatrix} \Bigg\} \\
&= \de_{x,1} [y] \Bigg[ [y-1] \begin{bmatrix}
x-1 & y-2 \\
z & w+1
\end{bmatrix} \\
&\quad + \varsigma_i q^{z-x+1}q^{w-y+1} \Bigg\{ \de_{z,1} \begin{bmatrix}
x & y-1 \\
z-1 & w
\end{bmatrix} + q^{x-z-1}[w] \begin{bmatrix}
x-1 & y \\
z & w-1
\end{bmatrix} \Bigg\} \Bigg] \\
&\quad + \kappa_i q^{z-x+1}q^{w-y+1} 
\begin{bmatrix}
x-1 & y-1 \\
z & w
\end{bmatrix} \\
& \quad - \td{q}^{1-x}q^{-y} \de_{z,1}[w] \Bigg\{ \de_{x,1}q^{-y+w-1} \begin{bmatrix}
x-1 & y \\
z & w-1
\end{bmatrix} + [y] \begin{bmatrix}
x & y-1 \\
z-1 & w
\end{bmatrix} \\
&\quad + \varsigma_i q^{z-1-x}q^{w-1-y}q^{x-z+1}[w-1] \begin{bmatrix}
x & y+1 \\
z-1 & w-2
\end{bmatrix} + \kappa_i q^{-x+z-1}q^{-y+w-1} 
\begin{bmatrix}
x & y \\
z-1 & w-1
\end{bmatrix}\Bigg\} .
\end{align*}}}
Comparing the coefficients appearing on both sides, one can check that $e_0 \msf{B}_i {\bf M'}= \msf{B}_i e_0 {\bf M'}$ if and only if $\varsigma_i=-q^{-1}$ and $\kappa_i=0$.
\medskip

{\em Case 2.} Suppose that $X=C$.
By the same reason as in {\em Case 1}, it is enough to show that $e_0\msf{B}_i {\bf M}' = \msf{B}_ie_0 {\bf M}'$,  where ${\bf M}'$ is the $2 \times 1$ submatrix
\begin{equation*}
{\bf M}' =
\begin{bmatrix}
x \\ y
\end{bmatrix}
=
\begin{bmatrix}
m_{i 1} \\ m_{i+1 1}
\end{bmatrix}. 
\end{equation*}
We have
\begin{align*}
e_0 \msf{B}_i \begin{bmatrix}
x \\
y
\end{bmatrix} &= e_0 \Bigg\{ [x] \begin{bmatrix}
x-1 \\
y+1
\end{bmatrix} + \varsigma_i q^{y-x} [y] \begin{bmatrix}
x+1 \\
y-1
\end{bmatrix}  + \kappa_i q^{y-x} \begin{bmatrix}
x \\
y
\end{bmatrix} \Bigg\} \\
&= [x] \Bigg\{ \frac{[x-1][x-2]}{[2]^2} \begin{bmatrix}
x-3 \\ y+1
\end{bmatrix} + q^{-2(x-1)-1} \frac{[y+1][y]}{[2]^2} \begin{bmatrix}
x-1 \\ y-1
\end{bmatrix} \Bigg\} \\
&\quad +\varsigma_i q^{y-x} [y] \Bigg\{ \frac{[x+1][x]}{[2]^2} \begin{bmatrix}
x-1 \\ y-1
\end{bmatrix} + q^{-2(x+1)-1} \frac{[y-1][y-2]}{[2]^2}\begin{bmatrix}
x+1 \\ y-3
\end{bmatrix} \Bigg\} \\
&\quad + \kappa_i q^{y-x} \Bigg\{ \frac{[x][x-1]}{[2]^2} \begin{bmatrix}
x-2 \\ y
\end{bmatrix} + q^{-2x-1} \frac{[y][y-1]}{[2]^2} \begin{bmatrix}
x \\ y-2
\end{bmatrix} \Bigg\},
\end{align*}
and
\begin{align*}
\msf{B}_ie_0 \begin{bmatrix}
x \\ y
\end{bmatrix} &= \msf{B}_i \Bigg\{ \frac{[x][x-1]}{[2]^2} \begin{bmatrix}
x-2 \\ y
\end{bmatrix} + q^{-2x-1} \frac{[y][y-1]}{[2]^2} \begin{bmatrix}
x \\ y-2
\end{bmatrix} \Bigg\} \\
&= \frac{[x][x-1]}{[2]^2} \Bigg\{ [x-2] \begin{bmatrix}
x-3 \\ y+1
\end{bmatrix} + \varsigma_i q^{y-x+2}[y] \begin{bmatrix}
x-1 \\ y-1
\end{bmatrix} + \kappa_i q^{y-x+2} \begin{bmatrix}
x-2 \\ y
\end{bmatrix} \Bigg\} \\
&\quad + q^{-2x-1} \frac{[y][y-1]}{[2]^2} \Bigg\{ [x] \begin{bmatrix}
x-1 \\ y-1
\end{bmatrix} + \varsigma_iq^{y-2-x}[y-2] \begin{bmatrix}
x+1 \\ y-3
\end{bmatrix} + \kappa_i q^{y-2-x} \begin{bmatrix}
x \\ y-2
\end{bmatrix} \Bigg\}.
\end{align*}
As in {\em Case 1}, one can check that $e_0 \msf{B}_i {\bf M'}= \msf{B}_i e_0 {\bf M'}$ if and only if $\varsigma_i=-q^{-1}$ and $\kappa_i=0$.

\subsection{Proof of Proposition \ref{prop:joint actions on Fock space}(2)}
\mbox{}%
Since $\msf{B}_i e_{0,s} {\bf M}=e_{0,s} \msf{B}_i {\bf M}$ for $i\in \msf{I}_\circ=\{\,2,4,\dots,2\ell-2\,\}$ and $s \neq \tfrac{i}{2}, \tfrac{i}{2}+1$, it remains to show that
\begin{equation}\label{eq:commutativity for sp}
 \msf{B}_i(e_{0, i/2}+e_{0, i/2+1}) {\bf M} = (e_{0, i/2}+e_{0, i/2+1})\msf{B}_i {\bf M}.
\end{equation}
Also we have for $i\in \msf{I}_\circ$
\begin{equation}\label{eq:T_bullet(e_i)}
\begin{split}
T_{w_\bullet}(\msf{e}_i)&= T_{1} T_{3}\cdots T_{2l-1}(\msf{e}_i) = T_{1} T_{3} \cdots T_{i-1} T_{i+1}(\msf{e}_i) \\
&= T_{1} T_{3} \cdots T_{i-1}(\msf{e}_{i+1}\msf{e}_i -q^{-1}\msf{e}_i\msf{e}_{i+1}) \\
&= T_{1} T_{3} \cdots T_{i-3}\{\msf{e}_{i+1}(\msf{e}_{i-1}\msf{e}_i-q^{-1}\msf{e}_i\msf{e}_{i-1})-q^{-1}(\msf{e}_{i-1}\msf{e}_i-q^{-1}\msf{e}_i\msf{e}_{i-1})\msf{e}_{i+1}\} \\
&= \msf{e}_{i+1}(\msf{e}_{i-1}\msf{e}_i-q^{-1}\msf{e}_i\msf{e}_{i-1})-q^{-1}(\msf{e}_{i-1}\msf{e}_i-q^{-1}\msf{e}_i\msf{e}_{i-1})\msf{e}_{i+1}.
\end{split}
\end{equation}
%Note that we have $\msf{e}_{i-1}\msf{e}_{i+1} = \msf{e}_{i+1}\msf{e}_{i-1}$ and $\msf{k}_{i}\msf{e}_{j}=q^{-1}\msf{e}_{j}\msf{k}_{i}$ for $|i-j|=1$, which will be used in further computations. %\comments{Bae: added (20250724)}

{\em Case 1.} Suppose that $X=D$.%\comments{Revised (250325)}
%
%We consider the case when $(\epsilon_1,\epsilon_2)=(0,1)$ since the proof for $(\epsilon_1,\epsilon_2)=(0,0)$ is similar. 

For a permutation ${\bf i}=(i_1,i_2,i_3)$ of $(i-1,i,i+1)$ and ${\bf t}=(t_1,t_2,t_3)$ with $1\le t_1,t_2,t_3\le n$, let
$\msf{e}_{{\bf i},{\bf t}}=\msf{e}_{i_1,t_1}\msf{e}_{i_2,t_2}\msf{e}_{i_3,t_3}$
As in \eqref{eq:B_i as operator on tensor for so}, we may write \eqref{eq:B_i for sp} as
\begin{equation}\label{eq:B_i as operator on tensor for sp}
\msf{B}_i=\sum_{1\le t\le n}\msf{f}_{i,t}+\varsigma_i \sum_{{\bf i},{\bf t}}c_{{\bf i},{\bf t}}\msf{e}_{{\bf i},{\bf t}}(\msf{k}_i^{-1})^{\ot n},
\end{equation}
for some $c_{{\bf i},{\bf t}}\in \Bbbk$ (see \eqref{eq:T_bullet(e_i)}).
Here $\msf{f}_{i,t}$ and $\msf{e}_{i,t}$ are the same as in {\em Case 1}.\smallskip

\noindent{\em Step 1.}
First, we can check easily that
\begin{equation}\label{eq:aux-1 Prop4.3(1)Case 1}
\msf{f}_{i,t}(e_{0, i/2}+e_{0, i/2+1}) {\bf M} = (e_{0, i/2}+e_{0, i/2+1}) \msf{f}_{i,t}{\bf M} \quad (t\ge 3).
\end{equation}%\comments{Will be checked by Bae}
\smallskip

\noindent{\em Step 2.}
Next, we claim that
\begin{equation}\label{eq:aux-2 Prop4.3(1)Case 1}
\msf{B}^{\le 2}_i(e_{0, i/2}+e_{0, i/2+1}) {\bf M} = (e_{0, i/2}+e_{0, i/2+1})\msf{B}^{\le 2}_i {\bf M},
\end{equation}
where $\msf{B}^{\le 2}_i$ is the partial sum of $\msf{B}_i$ \eqref{eq:B_i as operator on tensor for sp} with index $t\le 2$ for $\msf{f}_{i,t}$ and $\msf{e}_{i_k,t}$ ($k=1,2,3$). This is equivalent to proving \eqref{eq:aux-2 Prop4.3(1)Case 1} when restricted to the $4 \times 2$ submatrix ${\bf M}'$ of  ${\bf M}$ given by 
\begin{equation*}
{\bf M}'=
\begin{bmatrix}
x & y \\
z & w \\
a & b \\
c & d
\end{bmatrix} 
=
\begin{bmatrix}
m_{i-1,1} & m_{i-1,2} \\
m_{i,1} & m_{i,2} \\
m_{i+1,1} & m_{i+1,2} \\
m_{i+2,1} & m_{i+2,2} \\
\end{bmatrix}.
\end{equation*}

We compare the coefficients of the matrices ${\bf N}$ appearing in $e_0\msf{B}_i {\bf M}'$ and $\msf{B}_ie_0 {\bf M}'$. 
Recall that $\msf{B}_i=\msf{f}_i+{\varsigma_i}T_{w_\bullet}(\msf{e}_i)\msf{k}_i^{-1}$. Let us write  
\begin{equation*}
\begin{split}
&e_0\msf{B}_i{\bf M}'= \sum A_{\bf N}{\bf N}=\sum (A^-_{\bf N}+A^+_{\bf N}){\bf N},\\
&\msf{B}_ie_0{\bf M}'= \sum B_{\bf N}{\bf N}=\sum (B^-_{\bf N}+B^+_{\bf N}){\bf N},
\end{split}
\end{equation*}
where
\begin{equation*}
\begin{split}
&e_0\msf{f}_i{\bf M}'= \sum A^-_{\bf N}{\bf N},\ \
e_0 T_{w_\bullet}(\msf{e}_i)\msf{k}_i^{-1}{\bf M}'=\sum A^+_{\bf N}{\bf N},\\
&\msf{f}_ie_0{\bf M}'= \sum B^-_{{\bf N}}{\bf N},\ \
 T_{w_\bullet}(\msf{e}_i)\msf{k}_i^{-1}e_0{\bf M}'=\sum B^+_{\bf N}{\bf N}.
\end{split}
\end{equation*}
We also write $e_0=e_0^{(1)}+e_0^{(2)}$ as an operator on $\W^{\ot 2}$, 
where 
\begin{equation*}
\begin{split}
& e_0^{(1)} \lvert \mathbf{m} \rangle \otimes \lvert \mathbf{m}' \rangle = [m_2][m_1'] \lvert \mathbf{m}-\mathbf{e}_2 \rangle \otimes \lvert \mathbf{m}'-\mathbf{e}_1 \rangle,\\ 
&e_0^{(2)} \lvert \mathbf{m} \rangle \otimes \lvert \mathbf{m}' \rangle=-q_1^{m_1'}q_2^{-m_2}q^{-1}[m_1][m_2']\lvert \mathbf{m}-\mathbf{e}_{1} \rangle \otimes \lvert \mathbf{m}'-\mathbf{e}_2 \rangle,
\end{split}
\end{equation*}
for $\lvert \mathbf{m} \rangle \otimes \lvert\mathbf{m}'\rangle$.
When the action of $T_{w_\bullet}(\mathsf{e}_i)$ is restricted to ${\bf M'}$, it can be simplified as follows:
{\allowdisplaybreaks \begin{equation}\label{eq:simplified T}
    \begin{split}
    T_{w_\bullet}(\mathsf{e}_i) 
	&=(\msf{e}_{i+1}\msf{e}_{i-1}\msf{e}_{i} -q^{-1}\msf{e}_{i+1}\msf{e}_{i}\msf{e}_{i-1} -q^{-1}\msf{e}_{i-1}\msf{e}_{i}\msf{e}_{i+1} +q^{-2}\msf{e}_{i}\msf{e}_{i-1}\msf{e}_{i+1}) \ot 1  \\
    &\quad +(1-q^{-2})^2 \msf{e}_{i+1}\msf{e}_{i-1}\msf{k}_{i} \ot \msf{e}_{i} \\
    &\quad +(1-q^{-2}) \msf{e}_{i+1}\msf{k}_{i}\msf{k}_{i-1} \ot (\msf{e}_{i-1}\msf{e}_{i}-q^{-1}\msf{e}_{i}\msf{e}_{i-1}) \\
    &\quad +(1-q^{-2})\msf{e}_{i-1}\msf{k}_{i+1}\msf{k}_{i} \ot (\msf{e}_{i+1}\msf{e}_{i}-q^{-1}\msf{e}_{i}\msf{e}_{i+1}) \\
    &\quad +\msf{k}_{i+1}\msf{k}_{i}\msf{k}_{i-1} \ot (\msf{e}_{i+1}\msf{e}_{i-1}\msf{e}_{i}-q^{-1}\msf{e}_{i+1}\msf{e}_{i}\msf{e}_{i-1}-q^{-1}\msf{e}_{i-1}\msf{e}_{i}\msf{e}_{i+1}+q^{-2}\msf{e}_{i}\msf{e}_{i-1}\msf{e}_{i+1}).
    \end{split}
\end{equation}}%\comments{Bae: equation added (20250724)}
%\comments{Equation number should be the single one. (Use the expression "the first line in the second equation") (20250801) \\ Done by Bae}
\smallskip

\noindent {\em Step 2-1.} Consider first $A^-_{\bf N}$ and $B^-_{\bf N}$.
We have
{\small
\allowdisplaybreaks{\begin{align*}
    &e_0\mathsf{f}_i {\bf M}'
    =e_0 \Bigg\{q^{b-w}[z] \begin{bmatrix}
        x & y \\ z-1 & w \\ a+1 & b \\ c & d
    \end{bmatrix} + [w]\lang b \rang \begin{bmatrix}
        x & y \\ z & w-1 \\ a & b+1 \\ c & d
    \end{bmatrix} \Bigg\} \\
    &= q^{b-w}[z] \Bigg\{[y][z-1] \begin{bmatrix}
        x & y-1 \\ z-2 & w \\ a+1 & b \\ c & d
    \end{bmatrix} -q^{-z+1}q_2^{-y}q^{-1}[x][w] \begin{bmatrix}
        x-1 & y \\ z-1 & w-1 \\ a+1 & b \\ c & d
    \end{bmatrix} \\
    &\quad +q^{-x-z-1}q_2^{-y-w}[b][c] {\begin{bmatrix}
        x & y \\ z-1 & w \\ a+1 & b-1 \\ c-1 & d
    \end{bmatrix}} -q^{-x-z-c-2}q_2^{-y-w-b}[a+1][d] {\begin{bmatrix}
        x & y \\ z-1 & w \\ a & b \\ c & d-1
    \end{bmatrix}} \Bigg\} \\
    &\quad + [w]\lang b \rang \Bigg\{[y][z]\begin{bmatrix}
        x & y-1 \\ z-1 & w-1 \\ a & b+1 \\ c & d
    \end{bmatrix} - q^{-z-1}q_2^{-y}[x][w-1]\begin{bmatrix}
        x-1 & y \\ z & w-2 \\ a & b+1 \\ c & d
    \end{bmatrix} \\
    &\quad +q^{-x-z-2}q_2^{-y-w+1}[b+1][c]{\begin{bmatrix}
        x & y \\ z & w-1 \\ a & b \\ c-1 & d
    \end{bmatrix}} -q^{-x-z-c-3}q_2^{-y-w-b}[a][d]{\begin{bmatrix}
        x & y \\ z & w-1 \\ a-1 & b+1 \\ c & d-1
    \end{bmatrix}} \Bigg\}.
\end{align*}}
and 
\allowdisplaybreaks{
\begin{align*}
    &\mathsf{f}_ie_0{\bf M}'
    =\mathsf{f}_i \Bigg\{[y][z]\begin{bmatrix}
        x & y-1 \\ z-1 & w \\ a & b \\ c & d
    \end{bmatrix}-q^{-z-1}q_2^{-y}[x][w]\begin{bmatrix}
        x-2 & y \\ z & w-1 \\ a & b \\ c & d
    \end{bmatrix} \\
    &\quad +q^{-x-z-2}q_2^{-y-w}[b][c] \begin{bmatrix}
        x & y \\ z & w \\ a & b-1 \\ c-1 & d
    \end{bmatrix} -q^{-x-z-c-3}q_2^{-y-w-b}[a][d]\begin{bmatrix}
        x & y \\ z & w \\ a-1 & b \\ c & d-1
    \end{bmatrix} \Bigg\} \\
    &=[y][z]\Bigg\{q^{b-w}[z-1]\begin{bmatrix}
        x & y-1 \\ z-2 & w \\ a+1 & b \\ c & d
    \end{bmatrix} + [w]\lang b \rang \begin{bmatrix}
        x & y-1 \\ z-1 & w-1 \\ a & b+1 \\ c & d
    \end{bmatrix} \Bigg\} \\
    &\quad -q^{-z-1}q_2^{-y}[x][w] \Bigg\{ q^{b-w+1}[z] \begin{bmatrix}
        x-1 & y \\ z-1 & w-1 \\ a+1 & b \\ c & d
    \end{bmatrix} + [w-1]\lang b \rang \begin{bmatrix}
        x-2 & y \\ z & w-2 \\ a & b+1 \\ c & d
    \end{bmatrix} \Bigg\} \\
    &\quad +q^{-x-z-2}q_2^{-y-w}[b][c] \Bigg\{q^{b-w-1}[z]{\begin{bmatrix}
        x & y \\ z-1 & w \\ a+1 & b-1 \\ c-1 & d
    \end{bmatrix}} + [w]{\begin{bmatrix}
        x & y \\ z & w-1 \\ a & b \\ c-1 & d
    \end{bmatrix}} \Bigg\} \\
    &\quad -q^{-x-z-c-3}q_2^{-y-w-b}[a][d]\Bigg\{q^{b-w}[z]{\begin{bmatrix}
        x & y \\ z-1 & w \\ a & b \\ c & d-1
    \end{bmatrix}} + [w]\lang b \rang {\begin{bmatrix}
        x & y \\ z & w-1 \\ a-1 & b+1 \\ c & d-1
    \end{bmatrix}} \Bigg\}.
\end{align*}}}
We see that $A^-_{\bf N}=B^-_{\bf N}$ except when ${\bf N}$ is one of the following form
{\small
\begin{equation}\label{eq:nontrivial terms}
    \begin{bmatrix}
	x & y \\
	z-1 & w \\
	a+1 & b-1 \\
	c-1 & d
	\end{bmatrix}, \quad
	\begin{bmatrix}
	x & y \\
	z-1 & w \\
	a & b \\
	c & d-1
	\end{bmatrix}, \quad
	\begin{bmatrix}
	x & y \\
	z & w-1 \\
	a & b \\
	c-1 & d
	\end{bmatrix}, \quad
	\begin{bmatrix}
	x & y \\
	z & w-1 \\
	a-1 & b+1 \\
	c & d-1
	\end{bmatrix}.
    \end{equation}}
Also we can check directly that $A^+_{\bf N}=B^+_{\bf N}=0$ for ${\bf N}$'s appearing in $\msf{f}_ie_0{\bf M}'$ and $e_0\msf{f}_i{\bf M}'$ except for the ones in \eqref{eq:nontrivial terms}. 
\smallskip
    
\noindent {\em Step 2-2.} Let us check that $A_{\bf N}=B_{\bf N}$ for ${\bf N}$ in \eqref{eq:nontrivial terms}.
We consider, for example, the case when
{\small
\begin{equation}\label{eq:M''}
{\bf N}=
\begin{bmatrix}
x & y \\
z-1 & w \\
a+1 & b-1 \\
c-1 & d
\end{bmatrix}
\end{equation}}
since the proof of the other cases are similar.

Note that $e_0\msf{B}_i$ is a sum of linear operators of the form $T_1\ot T_2$, where $T_1$ and $T_2$ are products of $k_0, e_0^{(t)}, \msf{e}_s, \msf{k}_s^{\pm 1}$ and $\msf{f}_i$ ($t=1,2$, $s=i-1,i,i+1$). 
In fact, \eqref{eq:M''} appears when applying to ${\bf M}'$ the following summands in $e_0\msf{B}_i$; $(k_0 \otimes e_0^{(1)})(\msf{f}_i \otimes \msf{k}_i^{-1})$, the summands corresponding to the second line of \eqref{eq:simplified T} followed by $e_0^{(1)}\ot 1$, and the summands corresponding to the third line of \eqref{eq:simplified T} followed by $e_0^{(2)}\ot 1$. More precisely, it appears only in
{\allowdisplaybreaks
\begin{align*}
&\Big[(k_0 \otimes e_0^{(1)})(\msf{f}_i \otimes \msf{k}_i^{-1})+\varsigma_i \Big\{ (1-q^{-2})^2(e_0^{(1)} \otimes 1) (\msf{e}_{i+1}\msf{e}_{i-1}\msf{k}_{i} \ot \msf{e}_{i}) \\
&\qquad + (1-q^{-2})(e_0^{(2)} \otimes 1) \{\msf{e}_{i+1}\msf{k}_{i}\msf{k}_{i-1} \ot (\msf{e}_{i-1}\msf{e}_{i}-q^{-1}\msf{e}_{i}\msf{e}_{i-1})\} \Big\} (\msf{k}_i^{-1}\otimes \msf{k}_i^{-1})\Big]{\bf M'}.
\end{align*}}
Then we have
{\small{\allowdisplaybreaks
\begin{align*}
A_{\bf N} &= q^{-x-z-w+b-1}q_2^{-y-w}[z][b][c] \\
&\quad +\varsigma_i q^{-z-w+a+b}\{(1-q^{-2})q^{x-a}[y+1][z][b][c]\lang y \rang ([w+1]\lang w\rang -q^{-1}[w]\lang w-1\rang) \\
&\quad \qquad \qquad \qquad -(1-2q^{-2}+q^{-4})q^{-a}q_2^{-y}[x+1][z][w+1][b][c]\lang w\rang\}. 
\end{align*}}}
Similarly, \eqref{eq:M''} appears only in
{\allowdisplaybreaks
\begin{align*}
&\Big[(\msf{f}_i \otimes \msf{k}_i^{-1})(k_0 \otimes e_0^{(1)})+ \varsigma_i\Big\{ (1-q^{-2})^2(\msf{e}_{i+1}\msf{e}_{i-1}\msf{k}_i\ot \msf{e}_i) (\msf{k}_i^{-1} \otimes \msf{k}_i^{-1})(e_0^{(2)} \otimes 1) \\
&\quad + (1-q^{-2})\{\msf{e}_{i+1}\msf{k}_i\msf{k}_{i-1}\ot(\msf{e}_{i-1}\msf{e}_i-q^{-1}\msf{e}_i\msf{e}_{i-1})\}(\msf{k}_i^{-1}\otimes \msf{k}_i^{-1})(e_0^{(1)} \otimes 1)\Big\}\Big]{\bf M}',
\end{align*}}
so we have
{\small{\allowdisplaybreaks
\begin{align*}
B_{\bf N}&= q^{-x-z-w+b-3}q_2^{-y-w}[z][b][c] \\
&\quad +\varsigma_i(1-q^{-2})q^{x-z-w+b+1}[y][z][b][c]\lang y-1\rang ([w+1]\lang w\rang -q^{-1}[w]\lang w-1\rang) \\
&\quad - \varsigma_i (1-2q^{-2}+q^{-4})q^{-z-w+b}q_2^{-y}[x][z][w][b][c]\lang w-1\rang.
\end{align*}}} %\comments{Bae: "blue" area revised  (20250724)}
Now, it is straightforward to check that $A_{\bf N}=B_{\bf N}$ if and only if $\varsigma_i=-q$.

The coefficient $A_{\bf N}$ for the second matrix ${\bf N}$ in \eqref{eq:nontrivial terms} appears only when applying to ${\bf M}'$ the following summands in $e_0\msf{B}_i$; $(k_0 \ot e_0^{(2)})(\msf{f}_{i}\ot \msf{k}_{i}^{-1})$, the summands corresponding to the fifth line in \eqref{eq:simplified T} followed by $(e_0^{(1)}\ot ~1)$, and the summands corresponding to the fourth line in \eqref{eq:simplified T} followed by $(e_0^{(2)}\ot 1)$. Similarly, $B_{\bf N}$ for the same ${\bf N}$ appears only when applying to ${\bf M}'$ the summands in $\msf{B}_ie_0$ which correspond to the same terms used to compute $A_{\bf N}$, but in the opposite order of composition.

Next, the coefficient $A_{\bf N}$ for the third matrix in \eqref{eq:nontrivial terms} appears only when applying to ${\bf M}'$ the following summands in $e_0\msf{B}_i$; $(k_0 \ot e_0^{(1)})(1\ot\msf{f}_{i})$, and the summands corresponding to the first line in \eqref{eq:simplified T} followed by $e_0^{(2)}\ot 1$. Similarly, $B_{\bf N}$ for the same ${\bf N}$ appears only when applying to ${\bf M}'$ the summands in $\msf{B}_ie_0$ which correspond to the same terms used to compute $A_{\bf N}$, but in the opposite order. We also obtain $\varsigma_i=-q$ by comparing $A_{\bf N}$ and $B_{\bf N}$ for these two cases. 

The last matrix in \eqref{eq:nontrivial terms} appears only when applying to ${\bf M}'$ the terms $k_0\ot e_0^{(2)}$ in $e_0$ and $1\ot \msf{f}_{i}$ in $\msf{B}_i$, which commute each other.
\smallskip

\noindent{\em Step 2-3.} Now it remains to consider the case of ${\bf N}$ such that $A^-_{\bf N}=B^-_{\bf N}=0$, equivalently, ${\bf N}'$s appearing in $e_0 T_{w_\bullet}(\msf{e}_i)\msf{k}_i^{-1}{\bf M}'$ and $T_{w_\bullet}(\msf{e}_i)\msf{k}_i^{-1}e_0{\bf M}'$ except for the ones in \eqref{eq:nontrivial terms}, which is one of the following forms:

{\small{\allowdisplaybreaks
\begin{align*} 
    &\begin{bmatrix}
        x+1 & y-1 \\ z-1 & w \\ a & b \\ c-1 & d
    \end{bmatrix} \ \
    \begin{bmatrix}
        x+1 & y \\ z & w \\ a & b-1 \\ c-2 & d
    \end{bmatrix} \ \
    \begin{bmatrix}
        x+1 & y \\ z & w \\ a-1 & b \\ c-1 & d-1
    \end{bmatrix} \ \ 
    \begin{bmatrix}
        x+1 & y-1 \\ z-2 & w+1 \\ a+1 & b-1 \\ c-1 & d
    \end{bmatrix} \\
    &\begin{bmatrix}
        x+1 & y \\ z-1 & w+1 \\ a+1 & b-2 \\ c-2 & d
    \end{bmatrix} \ \
    \begin{bmatrix}
        x+1 & y \\ z-1 & w+1 \\ a & b-1 \\ c-1 & d-1
    \end{bmatrix} \ \ 
    \begin{bmatrix}
        x-1 & y+1 \\ z & w-1 \\ a+1 & b-1 \\ c-1 & d
    \end{bmatrix} \ \
    \begin{bmatrix}
        x & y+1 \\ z & w \\ a+1 & b-2 \\ c-2 & d
    \end{bmatrix} \\
    &
    \begin{bmatrix}
        x & y+1 \\ z & w \\ a & b-1 \\ c-1 & d-1
    \end{bmatrix} \ \
    \begin{bmatrix}
        x+1 & y-1 \\ z-2 & w+1 \\ a & b \\ c & d-1
    \end{bmatrix} \ \
    \begin{bmatrix}
        x+1 & y \\ z-1 & w+1 \\ a-1 & b \\ c & d-2
    \end{bmatrix} \ \ 
    \begin{bmatrix}
        x-1 & y+1 \\ z & w-1 \\ a & b \\ c & d-1
    \end{bmatrix} \ \
    \begin{bmatrix}
        x & y+1 \\ z & w \\ a-1 & b \\ c & d-2
    \end{bmatrix}.
\end{align*}}} %\comments{Bae: revised (20250724) due to \eqref{eq:simplified T} the cases of matrices reduced}

%\comments{List them if possible}
Let us consider the following term:
\begin{equation*} \label{eq:N1 in Step 2-3}
{\bf N}=
    \begin{bmatrix}
    x+1 & y-1 \\
    z-1 & w \\
    a & b \\
    c-1 & d
    \end{bmatrix}.
\end{equation*}
%
%\red{...the difference arising from the order of actions of $e_0$ and $T_{w_\bullet}(\mathsf{e}_i)$ is canceled out by the order of actions of $e_0$ and $\mathsf{k}_i^{-1}$...}

It appears when applying to ${\bf M}'$ the summands in $e_0(\varsigma_i T_{w_\bullet}(\msf{e}_i)\msf{k}_i^{-1})$ corresponding to the first line of \eqref{eq:simplified T}, that is, ${\bf N}$ appears only in
{\small{\allowdisplaybreaks
\begin{align*}
\varsigma_i (e_0^{(1)}\ot 1)\{( \msf{e}_{i+1}\msf{e}_{i-1}\msf{e}_i -q^{-1} \msf{e}_{i+1}\msf{e}_i\msf{e}_{i-1} -q^{-1}\msf{e}_{i-1}\msf{e}_i\msf{e}_{i+1} +q^{-2}\msf{e}_i\msf{e}_{i-1}\msf{e}_{i+1}) \ot 1\}(\msf{k}_i^{-1}\ot \msf{k}_i^{-1}){\bf M}'.
\end{align*}}}
%\comments{Bae: "blue" added (20250724)}

Then we have 
\begin{align*}
A^+_{\bf N}&=q^{a-z+b-w}\de_{y,1}[z][c]([a][z+1]-q^{-1}[a][z]-q^{-1}[a+1][z+1]+q^{-2}[a+1][z]) \\
&=q^{a-z+b-w}\de_{y,1}[z][c]([a]-q^{-1}[a+1])([z+1]-q^{-1}[z]).
\end{align*}
Similarly, we have
\begin{align*}
B^+_{\bf N}&=q^{a-z+b-w+1}\de_{y,1}[z][c]([a][z]-q^{-1}[a][z-1]-q^{-1}[a+1][z]+q^{-2}[a+1][z-1]) \\
&=q^{a-z+b-w {+1}}\de_{y,1}[z][c]([a]-q^{-1}[a+1])([{z}]-q^{-1}[{z-1}]).
\end{align*}
Hence $A^+_{\bf N}=B^+_{\bf N}$. The other cases can be checked in a similar way.

Therefore, by {\em Steps 2-1--2-3}, we have $A_{\bf N}=B_{\bf N}$ for all ${\bf N}$, which proves \eqref{eq:aux-2 Prop4.3(1)Case 1}. %that is, $e_0\msf{B}_i{\bf M}'=\msf{B}_ie_0{\bf M}'$.
\smallskip

\noindent{\em Step 3.} Finally, we claim that
\begin{equation}\label{eq:aux-3 Prop4.3(1)Case 1}
\msf{B}^{>2}_i(e_{0, i/2}+e_{0, i/2+1}) {\bf M} = (e_{0, i/2}+e_{0, i/2+1})\msf{B}^{>2}_i {\bf M},
\end{equation}
where 
\begin{equation*}
\msf{B}^{>2}_i=\varsigma_i \sum_{\substack{{\bf i},{\bf t}\\ \text{$t_k\ge 3$ for some $k$}}}c_{{\bf i},{\bf t}}\msf{e}_{{\bf i},{\bf t}}(\msf{k}_i^{-1})^{\ot n}.
\end{equation*}
This can be checked by tedious computations similar to but less involved than {\em Step 2}. We leave the details to the reader.%\comments{Will be checked by Bae}

Therefore \eqref{eq:commutativity for sp} follows from taking the sum of \eqref{eq:aux-1 Prop4.3(1)Case 1}, \eqref{eq:aux-2 Prop4.3(1)Case 1}, and \eqref{eq:aux-3 Prop4.3(1)Case 1}.
\smallskip

{\em Case 2.} Suppose that $X=C$. The proof is done by similar steps as in {\em Case 1}.%\comments{Revised (250325)}

\noindent{\em Step 1.} We can check \eqref{eq:aux-1 Prop4.3(1)Case 1} with $t\ge 3$ replaced by $t\ge 2$.
\smallskip

\noindent{\em Step 2.} We prove \eqref{eq:aux-2 Prop4.3(1)Case 1} with $\msf{B}^{\le 2}_i$ replaced by $\msf{B}^{\le 1}_i$, which is defined similarly.
In this case, one may prove \eqref{eq:aux-2 Prop4.3(1)Case 1} when restricted to the $4 \times 1$ submatrix submatrix ${\bf M}'$ of  ${\bf M}$ given by
\begin{equation*}
{\bf M}'=
\begin{bmatrix}
x \\ y \\ z \\ w
\end{bmatrix}
=
\begin{bmatrix}
m_{i-1,1} \\ 
m_{i,1} \\ 
m_{i+1,1} \\ 
m_{i+2,1}
\end{bmatrix}.
\end{equation*}
We have
{\small{\allowdisplaybreaks
\begin{align*}
&e_0 \msf{B}_i {\bf M}' = 
e_0 \Bigg\{ \de_{y,1}\de_{z,0} 
\begin{bmatrix}
x \\ y-1 \\ z+1 \\ w
\end{bmatrix}   
+ \varsigma_iq^{z-y}\times \\
& \Big(\de_{x,0}\de_{y,0}\de_{z,1}\de_{w,1}-q^{-1}\de_{x,0}\de_{y,1}\de_{z,1}\de_{w,1}  -q^{-1}\de_{x,0}\de_{y,0}\de_{z,0}\de_{w,1}+q^{-2}\de_{x,0}\de_{y,1}\de_{z,0}\de_{w,1}\Big)
\begin{bmatrix}
x+1 \\ y \\ z \\ w-1
\end{bmatrix} 
\Bigg\} \\
&=\de_{y,1}\de_{z,0}\de_{w,1}q^{-2}\td{q}^{-2x-2y+2}
\begin{bmatrix}
x \\ y-1 \\ z \\ w-1
\end{bmatrix} 
+ \varsigma_i q^{z-y}\times \\ 
& \Big(\de_{x,0}\de_{y,0}\de_{z,1}\de_{w,1}-q^{-1}\de_{x,0}\de_{y,1}\de_{z,1}\de_{w,1}  -q^{-1}\de_{x,0}\de_{y,0}\de_{z,0}\de_{w,1}+q^{-2}\de_{x,0}\de_{y,1}\de_{z,0}\de_{w,1}\Big) \de_{y,1} \begin{bmatrix}
x \\ y-1 \\ z \\ w-1
\end{bmatrix} \\
&= \Bigg\{ 
\de_{y,1}\de_{z,0}\de_{w,1}q^{-2}\td{q}^{-2x-2y+2} + 
\varsigma_i q^{z-y} 
\Big(-q^{-1}\de_{x,0}\de_{y,1}\de_{z,1}\de_{w,1} +q^{-2}\de_{x,0}\de_{y,1}\de_{z,0}\de_{w,1}\Big) \Bigg\} 
\begin{bmatrix}
x \\ y-1 \\ z \\ w-1
\end{bmatrix},
\end{align*}}
and
\begin{align*}
\msf{B}_i e_0{\bf M}' &= \msf{B}_i \Bigg\{ \de_{x,1}\de_{y,1}
\begin{bmatrix}
x-1 \\ y-1 \\ z \\ w
\end{bmatrix} 
+ q^{-2}\td{q}^{-2x-2y}\de_{z,1}\de_{w,1}
\begin{bmatrix}
x \\ y \\ z-1 \\ w-1
\end{bmatrix} \Bigg\} \\
&=\Bigg\{\de_{x,1}\de_{y,1}\varsigma_iq^{z-y+1}(\de_{z,1}\de_{w,1}-q^{-1}\de_{z,0}\de_{w,1})  +q^{-2}\td{q}^{-2x-2y}\de_{y,1}\de_{z,1}\de_{w,1} \Bigg\}
\begin{bmatrix}
x \\ y-1 \\ z \\ w-1
\end{bmatrix}.
\end{align*}}
By comparing coefficients, we can check that $e_0\msf{B}_i {\bf M'}=\msf{B}_i e_0{\bf M}'$ holds if and only if $\varsigma_i = -q$.

\smallskip

\noindent{\em Step 3.} Finally, we can check \eqref{eq:aux-3 Prop4.3(1)Case 1} with $\msf{B}_i^{>2}$ replaced by $\msf{B}_i^{>1}$.
It can be also checked by tedious computations similar to but less involved than {\em Step 2}.

\medskip

\section*{List of Notations}
\begin{longtable}{ll}
$\U_X(\e)$ &: the generalized quantum group of type $X=C,D$ (see Section \ref{subsec:pre}) \\
$U_X(\e)$  &: the quantum superalgebra (see Section \ref{subsec:pre-2}) \\
$\U_A(\e)$ &: the subalgebra of $\U_X(\e)$ of type $A$ \\
$\mf{osp}_X(\e)$  &: an orthosymplectic Lie superalgebra (see Section \ref{subsec:classical limit})  \\
$\W$  &: the $q$-oscillator $\U_X(\e)$-module in Propositions \ref{prop:osc-D} and \ref{prop:osc-C} \\
$\W^2$  &: the $q$-oscillator $\U_X(\e)$-module in Propositions \ref{prop:osc-D-2} and \ref{prop:osc-C-2}  \\
$\ms{W}$ &: a $\U_X(\e)$-module $\W$ or $\W^2$ \\
$W$ &: the classical limit of $\ms{W}$ defined in \eqref{eq:classical limit of W} \\
$\ms{P}_\e$ &: the set of $(n_0|n_1)$-hook partitions (see Section \ref{subsec:classical Howe duality}) \\
$G_\ell$ &: $O_{\ell}$ ($\ell\ge 2$), $Sp_{2\ell}$ ($\ell\ge 1$)\\
$\mc{P}(G_\ell)$ &: the set of partitions defined in \eqref{eq:partitions for classical groups} with $\mc{P}(G_\ell)_\e=\mc{P}(G_\ell)\cap \cP_\e$ \\
$V_{G_\ell}(\la)$ &: a finite-dimensional irreducible $G_\ell$-module for $\la\in \mc{P}(G_\ell)$\\
$\Lambda_{\lambda, \e}$ &: a weight defined in \eqref{eq:highest weight La} for $\la\in \mc{P}(G_\ell)_\e$ \\
$V^\la$ &: the irreducible $\mf{osp}_X(\e)$--module with highest weight $\Lambda_{\lambda, \e}$  \\
$\V^\la$  &: the irreducible $\U_X(\e)$-module with highest weight $\Lambda_{\lambda, \e}$  \\
$(\mf{g},\mf{k})$ &: a symmetric pair $(\mf{g},\mf{k})=(\mf{sl}_\ell,\mf{so}_\ell)$ ($\ell\ge 2$), $(\mf{sl}_{2\ell},\mf{sp}_{2\ell})$ ($\ell\ge 1$) \\
$({\bf U}(\mf{g}),{\bf U}^\imath(\mf{k}))$ &: the quantum symmetric pair associated to $(\mf{g},\mf{k})$ \\
$\mc{V}^\imath_{\mf{k}}(\la)$ &: the ${\bf U}^\imath(\mf{k})$-module ${\rm Hom}_{\U_X(\e)}\left(\mc{V}^\la,\ms{W}^{\ot \ell}\right)$ \\
${V}_{\mf{k}}^\imath(\la)$ &: the classical limit of   $\mc{V}^\imath_{\mf{k}}(\la)$
\end{longtable}
\medskip

\noindent{\bf Acknowledgement} The second author would like to thank M. Okado from whom he learned about generalized quantum groups.

\noindent{\bf Funding} This work is supported by the National Research Foundation of Korea(NRF) grant funded by the Korea government(MSIT) (No.2020R1A5A1016126 and RS-2024-00342349).
\medskip

%\noindent{\bf Data availability} This paper has no associated data.
%\medskip

%\noindent{\bf Conflict of Interest} On behalf of all authors, the corresponding author states that there is no conflict of interest.

\end{document}